\documentclass[a4paper,12pt]{article}
\usepackage{amssymb}
\usepackage{latexsym}
\usepackage{amsmath}
\usepackage{amsthm}
\usepackage{amsfonts}
\usepackage{amssymb}
\usepackage{graphicx}

\newtheorem{Th}{Theorem}
\newtheorem{Prop}{Proposition}
\newtheorem{Co}{Corollary}
\newtheorem{Lm}{Lemma}

\newtheorem{Dfi}{Definition}
\newtheorem{Rm}{Remark}

\newcommand{\be}{\begin{equation}}
\newcommand{\ee}{\end{equation}}
\newcommand{\R}{\mathbb{R}}
\newcommand{\N}{\mathbb{N}}
\newcommand{\C}{\mathbb{C}}
\newcommand{\Z}{\mathbb{Z}}

\newcommand\res{\mathop{\hbox{\vrule height 7pt width .5pt depth 0pt
\vrule height .5pt width 6pt depth 0pt}}\nolimits}

\newcommand{\reset}{\setcounter{equation}{0}\setcounter{Th}{0}\setcounter{Prop}{0}\setcounter{Co}{0}
\setcounter{Lm}{0}\setcounter{Rm}{0}}

\def\ti{\tilde}
\def\lf{\left}
\def\rg{\right}

\def\al{\alpha}
\def\la{\lambda}

\def\ep{\varepsilon}
\def\ds{\displaystyle}
\def\ov{\overline}
\def\Om{\Omega}
\def\om{\omega}
\def\p{\partial}

\def\res{\mathop{\hbox{\vrule height 7pt width .5pt 
depth 0pt\vrule height .5pt width 6pt depth 0pt}}\nolimits}

\begin{document}
\title{The Variations of Yang-Mills Lagrangian.}
\author{ Tristan Rivi\`ere\footnote{Forschungsinstitut f\"ur Mathematik, ETH Zentrum,
CH-8093 Z\"urich, Switzerland.}}
\date{ }
\maketitle
\large
\section{Introduction}

Yang-Mills theory is growing at the interface between high energy physics and mathematics
It is well known that Yang-Mills theory and Gauge theory in general had a profound impact on
the development of modern differential and algebraic geometry. One could quote Donaldson invariants
 in four dimensional differential topology, Hitchin Kobayashi conjecture relating the existence of Hermitian-Einstein metric on holomorphic bundles
over K\"ahler manifolds and Mumford stability in complex geometry or also Gromov Witten invariants
 in symplectic geometry...etc. While the influence of Gauge theory in geometry is quite notorious, one
tends sometimes to forget that Yang-Mills theory has been also at the heart of fundamental progresses
in the non-linear analysis of Partial Differential Equations in the last decades.
The purpose of this mini-course is to present  the variations of this important lagrangian. We shall raise
analysis question such as existence and regularity of Yang-Mills minimizers or such as  the compactification of the 
moduli space of critical points to Yang-mills lagrangian
 in general.
 
 \medskip

\noindent{\bf Acknowledgments.} {\it These notes corresponds to a course given by the author at the 9th summer school in Differential Geometry at the
Korean Institute for Advanced Studies. The author is very grateful to the institute for it's hospitality and the excellent working conditions during his stay.}
\section{The Plateau Problem.}

Before to move to the Yang-Mills minimization problem we will first recall some fundamental facts regarding the minimization of the area
in the parametric approach and some elements of the resolution of the so called {\it Plateau problem}.

Let $\Gamma$ be a simple closed Jordan Curve in ${\R}^3$  : there exists $\gamma\in C^0(S^1,{\R}^3)$ such that $\Gamma=\gamma(S^1)$.

{\bf Plateau problem :} {\it Find a $C^1$ immersion $u$ of the two dimensional disc $D^2$ which is continuous up to the boundary, whose restriction to
$\p D^2$ is an homeomorphism and which minimizes the area
\[
\mbox{Area}(u)=\int_{D^2}|\p_{x_1}u\times\p_{x_2}u|\ dx_1\,dx_2
\]}
The area is a fairly degenerated functional :

\begin{itemize}
\item[i)] It has a huge invariance group : Diff$(D^2)$, the group of diffeomorphism of the disc.
 Let $u_n$ be a minimizing sequence of the {\it Plateau problem} above, then the composition of $u_n$ with any sequence of diffeomorphism $\Psi_n$ of $D^2$ is still a minimizing sequence. The sequence $\Psi_n$ can for instance degenerate so that $u_n\circ\Psi_n$ converges to a point !
 \item[ii)] The area of $u$ does not control the image $u(D^2)$ which could be uniformly bounded while $u(D^2)$ becomes dense in ${\R}^3$ !
 \end{itemize}

In order to solve the Plateau Problem J.Douglas and independently Rad\'o minimize instead the {\it Dirichlet energy}
\[
\mbox{Area}(u)\le E(u)=\frac{1}{2}\int_{D^2}[|\p_{x_1} u|^2+|\p_{x_2} u|^2]\ dx_1\,dx_2
\]
the inequality comes from the pointwise inequality
\[
2\, |\p_{x_1}u\times\p_{x_2}u|\le |\p_{x_1} u|^2+|\p_{x_2} u|^2
\]
and equality holds if and only if  
\[
H(u):=|\p_{x_1}u|^2-|\p_{x_2}u|^2-2\, i\ \p_{x_1}u\cdot\p_{x_2}u=0
\]
This condition is satisfied if and only if the differential $du$ preserves angle that is to say $u$ is {\it conformal}.

The {\it Dirichlet energy} $E$ has much better properties than the area : it is

\begin{itemize}
\item[i)] it is {\it coercive} :
\[
\forall\, u_n\quad\mbox{s.t.} \quad{\limsup}_{n\rightarrow+\infty}E(u_n)<+\infty\quad\exists\,u_{n'}\rightharpoonup u_\infty\quad\mbox{in }\quad W^{1,2}
\]
\item[ii)] {\it lower semi-continuous} 
\[
E(u_\infty)\le {\liminf}_{n\rightarrow+\infty}E(u_n)
\]
\item[iii)] it's {\it invariance group}
in the domain  is reduced to the 3-dimensional M\"obius group ${\mathcal M}(D^2)$
\[
{\mathcal M}(D^2):=\left\{\Psi(z)=e^{i\theta}\frac{z-a}{1-\ov{a}z}\quad \mbox{ s.t. } \theta\in {\R}\quad |a|<1\right\}
\]
\end{itemize}
Why minimizing $E$ instead of the area has any chance to give a solution to the Plateau problem ?

Let $u$ be an immersion of the disc $D^2$ and denote $g_u$ the pull back of the standard flat metric $g_{{\R}^3}$ in ${\R}^3$ : $g_u:=u^\ast g_{{\R}^3}$.
The {\it uniformization theorem} on $D^2$ gives the existence of a diffeomorphism  $\Psi$ of the disc such that
\[
\Psi^\ast g_u=e^{2\la}\ [dx_1^2+dx_2^2]\quad(\mbox{i.e.} \quad u\circ\Psi\mbox{ is conformal})
\]
for some function $\la$ from $D^2$ into ${\R}$. So if $u$ minimizes $E$ in the desired class we can replace it by $v:=u\circ \Psi$ which is conformal and for which 
\[
E(v)=\mbox{Area}(v)=\mbox{Area}(u)\le E(u)
\]
hence $u$ is conformal. Assuming now it is not a minimizer of the area, we would then find $w$ in the desired class such that 
\[
A(w)<A(u)=E(u)
\]
and taking again $\Psi'$ s.t. $E(w\circ \Psi')=A(w)$ we would contradict that $u$ minimizes $E$.

Hence the difficulty is to find a minimizer of $E$ in the class of $C^1$ immersions sending continuously and monotonically $\p D^2$ into $\Gamma$.

One introduces
\[
{\mathcal F}:=\lf\{
\begin{array}{l}
u\in W^{1,2}(D^2)\quad;\quad u\in C^0(\p D^2,\Gamma)\quad\\[5mm]
\mbox{ and } u\mbox{ is monotone from } \p D^2\simeq S^1\mbox{ into }\Gamma
\end{array}\rg\}
\]

Fixing the images of three distinct points in $\p D^2$ in order to kill the action of the remaining {\it gauge group} ${\mathcal M}(D^2)$ one proves the existence
of a minimizer of $E$ in ${\mathcal F}$ which happens to be a solution to the Plateau problem (A thorough analysis is required to prove that the minimizer 
is indeed a $C^1$ immersion).

\subsection{The conformal parametrization choice as a Coulomb gauge.}

As we have seen the conformal parametrizations of immersed discs play a central role in the resolution of the Plateau problem. In the present
subsection we establish a one to one correspondence between this choice of conformal parametrization and a {\it Coulomb gauge} choice.

Let $u$ be a conformal immersion of the disc $D^2$ (i.e. $H(u)\equiv 0$ on $D^2$). Let $\la\in{\R}$ such that
\[
e^{\la}:=|\p_{x_1}u|=|\p_{x_2}u|\quad .
\]
Introduce the {\it moving frame} associated to this parametrization
\[
\vec{e}_j:=e^{-\la}\, \p_{x_j}u\quad\quad\mbox{ for }j=1,2\quad.
\]
The family $(\vec{e}_1,\vec{e}_2)$ realizes an orthonormal basis of the tangent space of $u(D^2)$ at $u(x_1,x_2)$. This can be also interpreted as
a section of the {\it frame bundle} of $u(D^2)$ equipped with the induced metric $g_{{\R}^3}$.

A simple computation gives
\be
\label{II.1}
\begin{array}{l}
\ds\mbox{div}\lf(\vec{e}_1,\nabla\vec{e}_2\rg)=\p_{x_1}\lf[e^{-2\la}\ \p_{x_1}u\cdot\p^2_{x_1x_2}u\rg]+ \p_{x_2}\lf[e^{-2\la}\ \p_{x_1}u\cdot\p^2_{x^2_2}u\rg]\\[5mm]
\ds\quad\quad=2^{-1}\p_{x_1}\lf[e^{-2\la}\,\p_{x_2}|\p_{x_1}u|^2\rg]-2^{-1}\p_{x_2}\lf[e^{-2\la}\,\p_{x_1}|\p_{x_2}u|^2\rg]\\[5mm]
\ds\quad\quad=\p^2_{x_1x_2}\la-\p^2_{x_2x_1}\la=0
\end{array}
\ee
In other words, introducing the 1 form on $D^2$ given by $A:=\vec{e}_1\cdot d\vec{e}_2$, which is nothing but the {\it connection form} associated
to the {\it Levi-Civita connection} induced by $g_u=u^\ast g_{{\R}^3}$ on the corresponding frame bundle for the trivialization given by $(\vec{e}_1,\vec{e}_2)$,
equation (\ref{II.1}) becomes
\be
\label{II.2}
d^{\ast_g}A=d^{\ast_g} \lf(\vec{e}_1\cdot d\vec{e}_2\rg)=0
\ee
where $\ast_g$ is the {\it Hodge operator} associated to the induced metric $g_u$. The equation (\ref{II.2}) is known as being the {\it Coulomb condition}.
We will see again this condition in the following sections and it is playing a central role in the mini-course.

Vice versa one proves, see for instance \cite{He} or \cite{Ri1}, that for any immersion $u$, non necessarily conformal, any frame $(\vec{e}_1,\vec{e}_2)$ satisfying the {\it Coulomb condition} (\ref{II.2}) corresponds
to a conformal parametrization (i.e. $\exists\,\Psi\in\,$Diff$(D^2)$ s.t. $v:=u\circ\Psi$ is conformal and $\vec{e}_j=|\p_{x_j}v|^{-1}\,\p_{x_j}v$ . This observation is the basis of the {\it Chern method} for constructing conformal or isothermal coordinates.
\section{A Plateau type problem on the lack of integrability}
\reset

In the rest of the class $G$ denotes an arbitrary compact Lie group. We will sometime restrict to the case where $G$ is a special unitary group $SU(n)$ and we will
mention it explicitely.
The corresponding Lie algebra will be denoted by ${\mathcal G}$ and the neutral element of $G$ is denoted $e$.

\subsection{Horizontal equivariant plane distributions.}
\subsubsection{The definition.}

Consider the simplest principal fiber structure $P:=B^m\times G$ where $B^m$ is the unit $m-$dimensional ball of ${\R}^m$. Denote by  $\pi$ the projection map
which to $\xi=(x,h)\in P$ assigns the {\it base point } $x$. Denote by $Gr_m(TP)$ to be 
the {\it Grassmanian} of $m-$ dimensional subspaces of the tangent bundle to $P$.

 We define the notion
of {\it equivariant horizontal distribution of plane} to be a map
\[
\begin{array}{l}
\ds H\ :\ P=B^m\times G\ \longrightarrow Gr_m(TP)\\[5mm]
\ds\quad\quad\quad \xi=(x,h)\ \longrightarrow\ H_\xi
\end{array}
\]
satisfying the following 3 conditions

\begin{itemize}
\item[i)] the {\it bundle condition} :
\[
\forall\, \xi\in\, P\quad\quad H_{\xi}\in T_\xi P\quad,
\]
\item[ii)] the {\it horizontality condition}
\[
\forall\, \xi\in\, P\quad\quad \pi_\ast H_{\xi}=T_{\pi(\xi)}B^m\quad,
\]
\item[iii)] the {\it equivariance condition}
\[
\forall\, \xi\in\, P\quad\forall\,g\in G\quad\quad (R_g)_\ast H_\xi=H_{R_g(\xi)}
\]
where $R_g$ is the {\it right multiplication map} by $g$ on $P$ which to any $\xi=(x,h)$ assigns $R_g(\xi):=(x,h\,g)$.
\end{itemize}

\subsubsection{Characterizations of equivariant horizontal distribution of plane by $1-$forms on $B^m$ taking values into ${\mathcal G}$.}

Let $H$ be an equivariant horizontal distribution of plane in $P=B^m\times G$. Clearly the following holds
\[
\forall\, \xi=(x,h)\in P\quad\forall X\,\in T_xB^m\quad\exists\, ! \, X^H(\xi)\in T_{\xi}P\quad\mbox{ s.t. }\pi_\ast X^H=X
\]
The vector $X^H(\xi)$ is called the {\it horizontal lifting} of $X$ at $\xi$.

At the point $(x,e)$ (recall that $e$ denotes the neutral element of $G$) we identify $T_{(x,e)}P\simeq T_xM\oplus{\mathcal G}$. Using this identification we deduce the existence of $A_{x}\cdot X$ such that
\[
X^H(x,e)=(X,-\,A_x\cdot X)
\]
The one form $A$ is called {\it connection 1-form} associated to $H$. The linearity of $A_x$ with respect to $X$ is a straightforward consequence of the definition of $X^H$ and therefore $A$ defines a 1-form on $B^m$ taking values into ${\mathcal G}$.

For any element $B\in {\mathcal G}\simeq T_eG$ we denote by $B^\ast$ the unique vector-field on $G$ satisfying
\[
B^\ast(e)=B\quad\mbox{and }\quad\forall g\in G \quad B^\ast(g):=(R_g)_\ast B
\]
and by an abuse of notation $B^\ast(g)$ is simply denoted $B\,g$. Using this notation we have
\[
\begin{array}{l}
\forall\, \xi=(x,h)\in P\quad\forall X\in T_xB^m\quad \\[5mm]
\quad X^H(\xi)=(X,-(A_x\cdot X)^\ast(h))=(R_h)_\ast X^H(x,e)
\end{array}
\]
At any point $\xi\in P$ Any vector $Z\in T_\xi P$ admits a decomposition according to $H$ : we denote by $Z^V$ the projection parallel to $H_\xi$ onto the tangent
plane to the vertical fiber  given by the kernel of $\pi_\ast$ :
\[
Z^V:=Z-(\pi_\ast Z)^H\quad.
\]
\subsection{The lack of integrability of equivariant horizontal distribution of planes.}

A $m-$dimensional plane distribution $H$ is said to be {\it integrable} if it identifies at every point with the tangent space to a $m-$dimensional {\it foliation}.

We aim to ''measure'' the {\it lack of integrability} of an equivariant horizontal distribution of planes. To that aim we recall the following classical result
\begin{Th} [Frobenius]
An $m-$dimensional plane distribution $H$ is integrable if and only if for any pair of vector fields $Y$ and $Z$ contained in $H$ at every
point the bracket $[Y,Z]$ is still included in $H$ at every point.
\hfill$\Box$
\end{Th}
In the particular case of equivariant horizontal distribution of planes in $P=B^m\times G$ we have that $H$ is integrable if and only if
\[
\forall\, X,Y\mbox{ vector-fields on }B^m\quad\quad [X^H,Y^H]^V\equiv 0\quad.
\]
We shall now compute $[X^H,Y^H]^V$ in terms of the one form $A$.

We write
\be
\label{III.1}
\begin{array}{l}
[X^H,Y^H]_{(x,e)}=[(X, -(A\cdot X)^\ast),(Y, -(A\cdot Y)^\ast)]_{(x,e)}\\[5mm]
\quad= [(X,0),(Y,0)]_{(x,e)}+[(X,0),(0,-(A\cdot Y)^\ast)]_{(x,e)}\\[5mm]
\quad+[(0,-(A\cdot X)^\ast),(Y,0)]_{(x,e)}+[(0,(A\cdot X)^\ast),(0,(A\cdot Y)^\ast)]_{(x,e)}
\end{array}
\ee
The definition of the Bracket operation on the Lie algebra ${\mathcal G}$ together with the commutation of the 
vector-field bracket operation with the push-forward operation of the right multiplication map gives that
\be
\label{III.2}
[(A\cdot X)^\ast),(A\cdot Y)^\ast]=([A\cdot X,A\cdot Y])^\ast
\ee
where the brackets in the r.h.s. of the identity is the Lie algebra bracket operation. The definition of the Lie bracket of vector fields gives
\be
\label{III.3}
[(X,0),(0,-(A\cdot Y)^\ast)]_{(x,e)}=(0,-\, d(A\cdot Y)\cdot X)
\ee
Finally we write
\be
\label{III.4}
[(X,0),(Y,0)]=([X,Y],0)=([X,Y],-A\cdot[X,Y])+(0,A\cdot[X,Y])
\ee
Combining (\ref{III.1}), (\ref{III.2}), (\ref{III.3}) and (\ref{III.4}) gives
\[
[X^H,Y^H]^V=d(A\cdot X)\cdot Y-d(A\cdot Y)\cdot X+A([X,Y])+[A\cdot X,A\cdot Y]
\]
and using Cartan formula on the expression of the exterior derivative of a one form we obtain
\be
\label{III.5}
[X^H,Y^H]^V=dA(X,Y)+[A\cdot X,A\cdot Y]\quad.
\ee
The two form we obtained
\be
\label{III.6}
F_A(X,Y):=dA(X,Y)+[A\cdot X,A\cdot Y]
\ee
is the so called {\it curvature} of the plane distribution $H$ and ''measures'' the lack of integrability of $H$. It will be conventionally denoted
\[
F_A=dA+\frac{1}{2}[A\wedge A]\quad \mbox{ or simply }\quad F_A=dA+A\wedge A\quad.
\]
The Lie algebra, and hence the compact Lie group, is equipped with the Killing form associated to a finite dimensional representation, hence unitary, for which
the form defines a scalar product invariant under adjoint action. For instance ${\mathcal G}=o(n)$ or ${\mathcal G}=u(n)$
\[
<B,C>= -\ \mbox{Tr}(B\,C)\quad.
\]
If the Lie algebra ${\mathcal G}$ is semi-simple : it is a direct sum of Lie algebras with no non trivial ideal, the Lie algebra is equipped with the {\it Killing scalar product} :
\[
<B,C>:=-\mbox{Tr}(\mbox{ad}(B)\, \mbox{ad}(C))
\]
where $\mbox{ad}(B)$ is the following endomorphism of ${\mathcal G}$ : $\mbox{ad}(B)(D):=[B,D]$. 

\medskip

The Lagrangian we are considering for measuring the lack of integrability of the plane distribution $H$ is just the $L^2$ norm of the curvature
\[
\int_{B^m}\sum_{i<j}\lf||[\p_{x_i}^H,\p_{x_j}^H]^V\rg|^2\ dx^m=\int_{B^m}|F_A|^2\ dx^m\quad,
\]
where $dx^m$ is the canonical volume form on $B^m$.  The $L^2$ norm of the curvature is also known as being the {\it Yang-Mills energy}
of the connection form $A$ and is denoted $YM(A)$. We can now state the main problem this mini-course is addressing.

{\bf Yang-Mills Plateau Problem :} {\it Let $\eta$ be a 1-form on $\p B^m$ taking values into a Lie algebra ${\mathcal G}$ of a compact Lie group
$G$ does there exists a 1-form $A$ into ${\mathcal G}$ realizing
\[
\inf\lf\{YM(A)=\int_{B^m}|dA+A\wedge A|^2\ dx^m\quad;\quad \iota_{\p B^m}^\ast A=\eta\rg\}
\]
where $\iota_{\p B^m}$ is the canonical inclusion of the boundary $\p B^m$ into ${\R}^m$.}

In other words we are asking the following question : {\it given an equivariant horizontal plane distribution over the boundary of the unit ball in ${\R}^m$, can one extend it inside the ball in an optimal way with respect to the
$L^2$ norm of the integrability defect.}

\medskip

In order to study this variational problem we first have to identify it's invariance group.

\subsection{The Gauge Invariance.}

In this subsection we identify the group of the Yang-Mill Plateau problem corresponding to the diffeomorphism group of the disc for the area
in the classical Plateau problem.

Let $g$ be a map from $B^m$ into $G$. We denote by $L_{g^{-1}}$ the left multiplication by ${g^{-1}}$ defined as follows
\[
\begin{array}{l}
\ds L_{g^{-1}}\ ;\ P=B^m\times G\ \longrightarrow P\\[5mm]
\ds \quad \xi=(x,h)\ \longrightarrow\ (x,g^{-1}\, h)
\end{array}
\]
Let $H$ be an {\it equivariant horizontal distribution of planes} on $P$ we observe that the push-forward by $L_{g^{-1}}$ of $H$, $(L_{g^{-1}})_\ast H$,
is still an {\it equivariant horizontal distribution of planes}. We now compute the connection 1-form associated to this new distribution.

Let $X$ be a vector of $T_x B^m$ and $x(t)$ a path in $B^m$ such that $\dot{x}(0)=X$. Let $h(t)\in G$ such that $\xi(t):=(x(t),h(t))$ is the horizontal lifting of $x(t)$
starting at the neutral element $e$ of $G$ (i.e. $\dot{\xi}=(\dot{x})^H(\xi(t))$ and $\xi(0)=(x,e)$). Since $\dot{xi}^V=0$ we have in particular
\[
\frac{d h}{dt}(0)=-A\cdot X
\]
The push forward by $L_{g^{-1}}$ of the horizontal vector $X^H(x,e)$ is the horizontal lifting of $X$ at $(x,g^{-1})$ for the distribution $(L_{g^{-1}})_\ast H$ :
$X^{(L_{g^{-1}})_\ast H}(x,g^{-1})$. Hence we have
\[
\begin{array}{l}
\ds X^{(L_{g^{-1}})_\ast H}(x,g^{-1})=(L_{g^{-1}})_\ast X^H(x,e)=(L_{g^{-1}})_\ast (X, - A\cdot X)\\[5mm]
\ds\quad=\frac{d}{dt}\lf(x(t), g^{-1} h(t)\rg)= (X, dg^{-1}\cdot X-g^{-1}\, A\cdot X)\\[5mm]
\ds\quad=\lf(X, -\lf(g^{-1}dg\cdot X+g^{-1}\, A\cdot X\,g\rg) g^{-1}\rg)\\[5mm]
\ds\quad=\lf(X, -\lf(g^{-1}dg\cdot X+g^{-1}\, A\cdot X\,g\rg)^\ast (g^{-1})\rg)
\end{array}
\]
Hence we have proved that the horizontal lift at $(x,e)$ for the new plane distribution $(L_{g^{-1}})_\ast H$ is $\lf(X, -\lf(g^{-1}dg\cdot X+g^{-1}\, A\cdot X\,g\rg)\rg)$ and the associated connection one form associated to the distribution $(L_{g^{-1}})_\ast H$ is
\[
A^g=g^{-1}dg\cdot X+g^{-1}\, A\cdot X\,g
\]
The curvature associated to this new distribution is given by
\[
F_{A^g}(X,Y)=dA^g(X,Y)+[A^g(X),A^g(Y)]
\]
We have in one hand
\[
dA^g(X,Y)=dg^{-1}\wedge dg(X,Y)+[dg^{-1}\,g\wedge g^{-1}Ag](X,Y)
\]
and in the other hand
\[
[A^g(X),A^g(Y)]=[g^{-1} dg\wedge g^{-1} A g](X,Y)+g^{-1}[A(X),A(Y)] g\quad.
\]
summing the two last identities and using the fact that $g^{-1}dg+dg^{-1}g=0$ gives finally
\[
F_{A^g}=g^{-1}\, F_A\, g
\]
Since the Killing scalar product on ${\mathcal G}$ is invariant under the adjoint action of $G$ we have
\[
YM(A^g)=YM(A)
\]
The action of $L_{g^{-1}}$ on the plane distribution $H$ leaves invariant it's Yang Mills energy and realizes therefore a ''huge'' invariance
group of the {\it Yang-Mills Plateau Problem} which is called the {\it Gauge group} of the problem.

\medskip

Exactly as for the classical Plateau problem we discussed in the first part of the mini-course the task for solving the Yang-Mills Plateau problem  will be to ''kill'' this gauge invariance and, here again, the Coulomb Gauge choices will be of great help.

\subsection{The Coulomb Gauges.}

We first start with the simplest group, the abelian group $G=S^1$. The Yang-Mills Plateau problem in this case becomes : 

{\it Find a minimizer of
\[
\inf\lf\{YM(A)=\int_{B^m}|dA|^2\ dx^m\quad;\quad \iota_{\p B^m}^\ast A=\eta\rg\}
\]
where $\eta$ is some given 1-form on the boundary $\p B^m$. }

In a reminiscent way to the classical Plateau problem our starting functional is degenerate and we shall replace it by a more coercive one 
\[
YM(A)\le E(A)=\int_{B^m} |dA|^2+|d^\ast A|^2\ dx^m
\]
with equality if and only if $d^\ast A=0$ (i.e. $A$ satisfies the Coulomb condition).

The following coercivity inequality holds
\be
\label{III.7}
\begin{array}{l}
\ds\forall\, A\in W^{1,2}(B^m,{\mathcal G})\quad\mbox{ s.t. } \iota_{\p B^m}^\ast A=\eta\quad\\[5mm]
\ds\quad \int_{B^m}|A|^2+\sum_{i,j=1}^m|\p_{x_i}A_j|^2\ dx^m\le C \lf[E(A)+\|\eta\|^2_{H^{1/2}(\p B^m)}\rg]
\end{array}
\ee
for some fixed constant independent of $\eta$ and $A$,  where $H^{1/2}(\p B^m)$ is the fractional
trace space of $W^{1,2}(B^m)$.  The convexity of $E$ in $W^{1,2}_\eta(\wedge^1B^m,{\mathcal G})$ together with the previous coercivity inequality
implies that the following problem admits a unique minimizer
\[
\min_{A\in W^{1,2}_\eta(\wedge^1B^m,{\mathcal G})} E(A)
\]
and it is the unique solution of the following system
\[
\lf\{
\begin{array}{l}
d^\ast dA_0=0\quad\quad\mbox{ in }{\mathcal D}'(B^m)\\[5mm]
d^\ast A_0=0\quad\quad\mbox{ in }{\mathcal D}'(B^m)\\[5mm]
\iota_{\p B^m}^\ast A_0=\eta
\end{array}
\rg.
\]
The components of $A$ are harmonic in $B^m$ and are therefore smooth moreover
 we have
 \[
 YM(A_0)=E(A_0)
 \]
Let now $B$ in $W^{1,2}_\eta(\wedge B^m,{\mathcal G})$ we claim that there exists a gauge change $g$ such that $YM(B)=YM(B^g)=E(B^g)$. This can be seen as follows. Let $\varphi$ be the solution of
\[
\lf\{
\begin{array}{l}
-\Delta \varphi=d^\ast B\quad\quad\mbox{ in }{\mathcal D}'(B^m)\\[5mm]
\varphi=0\quad\quad\mbox{ on }\p B^m
\end{array}
\rg.
\]
Hence we have
\[
\lf\{
\begin{array}{l}
d(B+d\varphi)= dB\quad\quad\mbox{ in }{\mathcal D}'(B^m)\\[5mm]
d^\ast (B+d\varphi)=0\quad\quad\mbox{ in }{\mathcal D}'(B^m)\\[5mm]
\iota_{\p B^m}(B+d\varphi)=\eta
\end{array}
\rg.
\]
Taking $g:=exp(i\varphi)$ we have $YM(B)=YM(B^g)=E(B^g)$.

Hence $A_0$ realizes
\[
\min_{A\in W^{1,2}_\eta(\wedge^1B^m,{\mathcal G})} YM(A)
\]
Indeed if there would be $B\in W^{1,2}_\eta(\wedge^1B^m,{\mathcal G})$ such that $YM(B)<YM(A_0)$ we choose $g$ such that
$YM(B^g)=E(B^g)$ and we would contradict the fact that $A_0$ minimizes $E$.

\medskip

Taking now a general compact Lie group $G$ we would also like to propose to minimize $E$ instead of $YM$ but we need first ensure 
that a Coulomb gauge always exists. We have the following lemma which answers positively to this last question

\begin{Lm}
\label{lm-III.1}
Let $A\in L^2(\wedge^1B^m,{\mathcal G})$. The following variational problem
\[
\inf_{g\in W^{1,2}_e(B^m,G)}\int_{B^m}|g^{-1}dg+g^{-1} A g|^2\ dx^m
\]
is achieved and each minimizer satisfies the Coulomb condition
\[
d^\ast(g^{-1} dg+g^{-1} A\, g)=0
\]
\hfill $\Box$
\end{Lm}
\noindent{\bf Proof of Lemma~\ref{lm-III.1}.} Let $g_k$ be a minimizing sequence. Since the group is compact we have that 
\[
\limsup_{k\rightarrow +\infty}\int_{B^m}|dg_k|^2\ dx^m<+\infty
\]
Hence, using Rellich Kondrachov's theorem there exists a subsequence $g_{k'}$ converging weakly in $W^{1,2}(B^m,G)$ to $g_\infty$
and strongly in every $L^p(B^m,G)$ space for any $p<+\infty$ hence $g_\infty\in W^{1,2}_e(B^m,G)$. The same holds for $g_k^{-1}$ and 
it's weak limit in $W^{1,2}(B^m,G)$ is the inverse of $g_\infty$. Hence we have
\[
g_k^{-1}dg_k+g_k^{-1}A\,g_k\rightharpoonup g_\infty^{-1} dg_\infty+g_\infty^{-1} A g_\infty\quad\quad\mbox{ in }{\mathcal D}'(B^m)
\]
The lower semi continuity of the $L^2$ norm implies that $g_\infty$ is a minimizer of (\ref{III.6}).

For any $U\in C^\infty_0(B^m,{\mathcal G})$ we introduce
\[
g_\infty(t):=g_\infty\exp(t\,U)\quad.
\]
We have
\[
\begin{array}{l}
g_\infty^{-1}(t) dg_\infty^{-1}(t) +g_\infty^{-1}(t)A g_\infty(t)= \exp(-t\,U)\,d\exp(t\,U)\\[5mm]
\quad\quad+\exp(-t\,U)\lf[g_\infty^{-1} dg_\infty+g_\infty^{-1} A g_\infty\rg]\exp(t\,U)
\end{array}
\]
Hence
\[
\frac{d}{dt}\lf[g_\infty^{-1}(t) dg_\infty^{-1}(t) +g_\infty^{-1}(t)A g_\infty(t)\rg]=dU-[U,A^{g_\infty}]
\]
Since $g_\infty$ is a minimizer  we have
\[
0=\frac{d}{dt}\int_{B^m}|A^{g_\infty(t)}|^2\ dx^m=2\int_{B^m}\lf<\lf(dU-[U,A^{g_\infty}]\rg)\cdot A^{g_\infty}\rg>\ dx^m
\]
We use the identity $<[U,V],W>=<U,[V,W]>$ to deduce that $$[U,A^{g_\infty}]\cdot A^{g_\infty}=0$$ and then we have proved that
for any $U\in C^\infty_0(B^m,{\mathcal G})$ we have
\[
0=\int_{B^m}dU\cdot A^{g_\infty}\ dx^m=\int_{B^m}<U, d^\ast A^{g_\infty}>\ dx^m
\]
This finishes the proof of the lemma.\hfill $\Box$

Since every connection form posses a Coulomb Gauge representative it is then tempting to minimize $E$ instead of $YM$ following the main
lines of the abelian case. However due to the non-linearity $A\wedge A$ in $F_A$ it is not clear whether the $E$ energy controls the $W^{1,2}$ norm of $A$
in a similar way of (\ref{III.7}) in the general case.

In fact the answer to that question is ''no'' as we can see in the following example. We take $G=SU(2)$ and we identify $su(2)$ with the imaginary quaternions.
On ${\R}^4$ we identify canonically the point of coordinates $(x_0,x_1,x_2,x_3)$ with the quaternion $\mathbf{x}:=x_0 + x_1\,\mathbf{i}+x_2\,\mathbf{j}+x_3\,\mathbf{k}$. For a quaternion $\mathbf{y}=y_0 + y_1\,\mathbf{i}+y_2\,\mathbf{j}+y_3\,\mathbf{k}$ we denote by $\Im(\mathbf{y})$ the element in  $su(2)$
given by $y_1\,\sigma_1+y_2\,\sigma_2+y_3\,\sigma_3$ where $\sigma_i$ are the Pauli matrices to which we identify ${\mathbf i}$, ${\mathbf j}$ and ${\mathbf k}$
\[
{\mathbf i}\leftrightarrow\sigma_1=\lf(\begin{array}{cc}
i& 0\\[2mm]
0 &-i
\end{array}\rg)\quad\quad{\mathbf j}\leftrightarrow\sigma_2=\lf(\begin{array}{cc}
0& 1\\[2mm]
-1 &0
\end{array}\rg)\quad\quad{\mathbf k}\leftrightarrow\sigma_3=\lf(\begin{array}{cc}
0& i\\[2mm]
i&0
\end{array}\rg)
\] forming an orthogonal basis of $su(2)$ with norms $\sqrt{2}$ for each vector of the basis. 

On ${B}^4$, for $\la\in {\R}^\ast_+$,
we consider the family of one forms into $su(2)$ given by
\[
A_\la:=\la^2\ \frac{\Im\lf(\mathbf{x}\,d\overline{\mathbf{x}}\rg)}{1+\la^2\,|\mathbf{x}|^2}
\]
The corresponding curvature is given by
\[
F_{A_\la}=\la^2\,\frac{d\mathbf{x}\wedge d\overline{\mathbf{x}}}{(1+\la^2\,|\mathbf{x}|^2)^2}
\]
One easily verifies that\footnote{ The somehow surprizing factor 48 comes from the fact that there are 6 curvature coordinates and each curvature coordinate has the form
\[
|(F_{A_1})_{ij}|^2=\frac{8}{(1+|x|^2)^2}
\] where we have used that the square of the norm of each Pauli matrix  is 2.}
\[
\lim_{\la\rightarrow +\infty}\int_{B^4}|F_{A_\la}|^2\ dx^4=\int_{{\R}^4}|F_{A_1}|^2\ dx^4=\int_{\R^4}48\ \frac{dx^4}{(1+|x|^2)^4}=8\pi^2<+\infty
\]
but one verifies also that
\[
\lim_{\la\rightarrow +\infty}\int_{B^4}\sum_{i,j=1^4}|\p_{x_i}(A_\la)_j|^2\ dx^4=+\infty
\]
we might then think that by changing the gauge we can avoid this blow up of the $W^{1,2}$ norm of the connection form
but, as we see now, this cannot be the case.
Consider the {\it second Chern form} $Tr(F_{A_\la}\wedge F_{A_\la})$, it satisfies
\be
\label{III.8}
\lim_{\la\rightarrow +\infty} \mbox{Tr}(F_{A_\la}\wedge F_{A_\la})=8\pi^2\,\delta_0\ dx^4
\ee
The {\it second Chern form} is invariant under gauge transformation and for any choice of Gauge $g$ this $4-$form, which has to be closed, is on $B^4$
the exterior derivative of the {\it transgression form } known as the {\it Chern-Simon 3-form} :
\[
\begin{array}{l}
\ds\forall \, g\ : B^4\ \rightarrow SU(2)\  \\[5mm]
\ds \mbox{Tr}(F_{A_\la}\wedge F_{A_\la})=d\lf[ \mbox{Tr}\lf(A^g_\la\wedge dA^g_\la+\frac{1}{3} A^g_\la\wedge[A^g_\la,A^g_\la]\rg)\rg]
\end{array}
\]
Assume now there would have been a gauge $g_\la$ s.t. 
\[
\liminf_{\la\rightarrow +\infty}\|A_\la^{g_\la}\|_{W^{1,2}(B^4)}<+\infty
\]
Then for some sequence $\la_k\rightarrow +\infty$, using Rellich Kondrachov theorem $A^k:=A_{\la_k}^{g_k}$ would weakly converge to some limit
$A^\infty$ in $W^{1,2}(\wedge^1B^4, su(2))$ and strongly in $L^p(\wedge^1B^4, su(2))$ for any $p<4$. Hence
\[
\mbox{Tr}\lf(A^k\wedge dA^k+\frac{1}{3} A^k\wedge[A^k,A^k]\rg)\rightharpoonup \mbox{Tr}\lf(A^\infty\wedge dA^\infty+\frac{1}{3} A^\infty\wedge[A^\infty,A^\infty]\rg)\]
in ${\mathcal D}'(B^4)$. Taking now the exterior derivative and using again the gauge invariance of the second Chern form we obtain 
\[
\mbox{Tr}(F_{A_{\la_k}}\wedge F_{A_{\la_k}})\rightharpoonup \mbox{Tr}(F_{A^\infty}\wedge F_{A^\infty})\quad\mbox{ in }{\mathcal D}'(B^4)
\]
Since $A^\infty$ is in $W^{1,2}$ the 4-form $\mbox{Tr}(F_{A^\infty}\wedge F_{A^\infty})$ is an $L^1$ function, however, from (\ref{III.8}) it is equal to the Dirac mass. This gives a contradiction and we have proved the following proposition
\begin{Prop}
\label{pr-III.1}
There exists $A^k\in W^{1,2}(B^4,su(2))$ such that
\[
\limsup_{k\rightarrow +\infty}\int_{B^4}|F_{A^k}|^2\ dx^4<+\infty
\]
but
\[
\liminf_{k\rightarrow +\infty}\inf\lf\{\int_{B^4}\sum_{i,j=1}^4|\p_{x_i}(A^k)^g_j|^2\ dx^4\ ;\ g\in W^{2,2}(B^4,SU(2))\rg\}=+\infty
\]
\hfill$\Box$
\end{Prop}
Hence by minimizing $E$ instead of $YM$ we don't get enough control on the minimizing sequence $A^k$ in order to extract a converging subsequence
to a solution to the Yang-Mills Plateau problem.

The situation would have been much better in dimension less than $4$ where a $W^{1,2}$ control of $A$ in terms of $E(A)$ do exist. In dimension
equal to 4, despite proposition~\ref{pr-III.1}, there is still a positive result in that line which says roughly that such a control do exist for some gauge
provided the {\it Yang-Mills energy} stays below some positive threshold. The following section is devoted to the proof of this result by K.Uhlenbeck.

\section{Uhlenbeck's Coulomb Gauge Extraction Method}
\reset

\subsection{Uhlenbeck's construction.}
We have seen that in dimension 4 - and higher of course- there is no hope to control the  $W^{1,2}$ norm of sequences of connection forms from the $E$ energy. The fact that the dimension $4$ is critical for this phenomenon comes form the optimal {\it Sobolev embedding}
\[
W^{1,2}(B^4)\hookrightarrow L^4(B^4)
\]
which does not hold in higher dimension. 

\begin{Th}
\label{th-III.2}
Let  $m\le 4$ and $G$ be a compact Lie group. There exists $\ep_{G}>0$ and $C_{G}>0$ such that for any $A\in W^{1,2}(B^m,{\mathcal G})$ satisfying
\[
\int_{B^m}|dA+A\wedge A|^2\ dx^m<\ep_G\quad,
\]
there exists $g\in W^{2,2}(B^m,G)$ such that
\be
\label{III.9}
\lf\{
\begin{array}{l}
\ds\int_{B^m}|A^g|^2+\sum_{i,j=1}^4|\p_{x_i}A^g_j|^2\, dx^m\le C_G\, \int_{B^m}|dA+A\wedge A|^2\, dx^m\\[8mm]
\ds d^\ast A^g=0\quad\quad\mbox{ in }B^m\\[5mm]
\ds\iota_{\p B^m}^\ast(\ast A^g)=0
\end{array}
\rg.
\ee
where $A^g=g^{-1}dg+g^{-1} A g$ and $\iota_{\p B^m}$ is the canonical inclusion map of the boundary of the unit ball into ${\R}^m$.\hfill$\Box$
\end{Th}

For $m<4$ the non linearity $A\wedge A$ is a compact perturbation of $dA$ and the problem is a perturbation to a simple linear one
that we solved in the abelian case. We will then restrict the presentation to the case $m=4$. We assume that the compact Lie group is represented by a subgroup of invertible matrices in $R^n$ for some $n\in {\N}^\ast$ which gives an isometric embedding of $G$ in an euclidian space.

We aim to solve the Coulomb equation $d^\ast A^g=0$ keeping this time a control of $A^g$ in $W^{1,2}$ (lemma~\ref{III.1} was only giving an $L^2$ control).
Since we have little energy the hope is to use a fixed point argument close to the zero connection and for $g$ close to the identity. The linearization of the {\it Coulomb non-linear elliptic PDE}
\[
\lf\{
\begin{array}{l}
\ds d^\ast \lf[g^{-1}\, dg\rg]=-\,d^\ast\lf[g^{-1}A\, g\rg]\quad\mbox{ in }B^4\\[5mm]
\ds \p_rg\, g^{-1}=-<A,\p_r>\quad\mbox{ on }\p B^4
\end{array}
\rg.
\]
for $A^t=t\, \omega$ and $g^t=exp(t\ U)$ gives
\[
\lf\{
\begin{array}{l}
\ds\Delta U=-d^\ast dU=d^\ast\omega\quad\mbox{ in }B^4\\[5mm]
\ds \p_r U=-<\omega,\p_r>\quad\mbox{ on }\p B^4
\end{array}
\rg.
\]
This linearized problem is  solvable for any map $\omega\in W^{1,2}(B^4,{\mathcal G})$ and we get a unique $U\in W^{2,2}(B^4,{\mathcal G})$
solving the previous linear equation.

However, in order to be able to apply the implicit function theorem we need  the following non-linear mapping to be \underbar{smooth}

\medskip

\[
\begin{array}{l}
\ds{\mathcal N}^0\ :\ W^{1,2}(B^4,{\mathcal G})\times W^{2,2}(B^4,{\mathcal G})\, \longrightarrow \, W^{1,2}(B^4,{\mathcal G})\times H^{1/2}(\p B^4,{\mathcal G})\\[7mm]
\ds\quad (\omega,U)\ \longrightarrow\ \lf(d^\ast \lf[g_U^{-1}\, dg_U+g_U^{-1}\omega\, g_U\rg],\p_r g_U\,g_U^{-1}-<\om,\p_r>\rg)
\end{array}
\]

\medskip

where $g_U:=\exp(U)$.  This is however \underbar{not the case} in dimension $4$. This is due to the fact
that $W^{2,2}$ does not embed in $C^0$ in 4 dimension but only in the {\it Bounded Mean Oscillation} space $BMO(B^4)$ hence 
simple algebraic operations such as the multiplication of two $G$ valued $W^{2,2}$ maps is not continuous in 4 dimension.

If one replaces however $W^{1,2}\times W^{2,2}$ by a ''slightly smaller space'' $W^{1,p}\times W^{2,p}$ for any $p>2$ ($p$ being as close as we want from $2$) then the space $W^{2,p}$ embeds continuously (and compactly) in $C^0$ and the map ${\mathcal N}^0$ becomes sudently \underbar{smooth} ! and a fixed point argument is conceivable in this smaller space.

Uhlenbeck's strategy consists in combining a fixed point argument in smaller spaces  - in which the problem is invertible - together with a continuity argument.

This method is rather generic in the sense that it can be applied to critical {\it extensions} or {\it lifting problems}  of maps  in the Sobolev space
$W^{1,m}(M^m)$ which misses to embed in $C^0$ but for which however the notion of homotopy class is well defined (see \cite{SU} and \cite{Wh}) and prevents to find global extensions or liftings when the norm of the map is too high. As an illustration we shall give two results. 
\begin{Th}
\label{th-III.1-a} 
For any $m\ge 1$,  and any compact Lie group $G$ there exists $\ep_{m,G}>0$ and $C_{m,G}>0$ such that for any map $g\in W^{1,m}(S^m,G)$ satisfying
\[
\int_{S^m}|dg|^m\ dvol_{S^m}<\ep_m
\]
there exists an extension $\ti{g}\in W^{1,{m+1}}(B^{m+1},G)$ equal to $g$ on $\p B^{m+1}$ such that
\[
\int_{B^{m+1}}|d\ti{g}|^{m+1}\ dx^{m+1}\le C_m\ \int_{S^m}|dg|^m\ dvol_{S^m}\quad.
\]
\hfill $\Box$
\end{Th}
\begin{Rm}
\label{rm-III.1-a}
The existence of such an extension $\ti{g}\in W^{1,{m+1}}(B^{m+1},G)$ is clearly not true for general $g\in W^{1,m}(S^m,S^m)$. Indeed,  consider for instance
$m=3$ and $G=SU(2)\simeq S^3$, if such an extension would exists one would have using Stokes theorem. 
\[
0=\int_{B^{4}}\ti{g}^\ast dx^{4}=\frac{1}{4}\int_{S^3}g^\ast dvol_{SU(2)}=\frac{|S^3|}{4}\mbox{ deg}(g)\quad,
\]
where deg$(g)$ is the topological degree of the map $g$ which is not necessarily zero. \hfill $\Box$
\end{Rm}
Theorem~\ref{th-III.1-a} is proved in \cite{PR2} using
 Uhlenbeck's method. The second example is the following one
\begin{Th}
\label{th-III.1-b}
Let $P$ be a $G$ principal bundle over $S^m$ where $G$ is a compact Lie group and where $\pi$ is the projection associated to this bundle. There exists $\ep_{m,G}>0$ and $C_{m,G}>0$ such that for any $g\in W^{1,m}(S^m,S^m)$ satisfying
\[
\int_{S^m}|d g|^m\ dvol_{S^m}<\ep_{m,G}
\]
there exists $v\in W^{1,m}(S^m,P)$ such that
\[
\int_{S^m}|dv|^m\ dvol_{S^m}\le C_{m,G}\ \int_{S^m}|d g|^m\ dvol_{S^m}\quad,
\]
and
\[
\pi\circ v=g\quad.
\]
\hfill $\Box$
\end{Th}
\medskip

\noindent{\bf Proof of theorem~\ref{th-III.2}.}

Fix some $2<p<4$. For any $\ep>0$ we introduce
\[
{\mathcal U}^\ep:=\lf\{A\in W^{1,p}(B^4,{\mathcal G})\quad\mbox{s.t. }\quad\int_{B^4}|F_A|^2\ dx^4<\ep\rg\}
\]
and for any $\ep>0$ and $C>0$ and  we consider
\[
{\mathcal V}^\ep_C:=\lf\{ 
\begin{array}{l}
\ds A\in {\mathcal U}^\ep\quad\mbox{ s. t. }\quad\exists\, g\in W^{2,p}(B^4,G)\\[5mm]
\ds \int_{B^4}|dA^g|_{g_{S^4}}^p\, dx^4\le C \int_{B^4}|F_A|^p\ dx^4\\[5mm]
\ds \int_{B^4}|dA^g|_{g_{S^4}}^2\, dx^4\le C \int_{B^4}|F_A|^2\ dx^4\\[5mm]
\ds\quad d^\ast A^g=0\quad\mbox{ and }\quad \iota_{\p B^4}\ast\,A^g=0
\end{array}
\rg\}
\]
The first goal is to show the following

\medskip

\noindent{\bf Claim}
\[
{ \exists\, \ep>0\quad C>0\quad\mbox{ s.t. }\quad{\mathcal V}^\ep_C={\mathcal U}^\ep}
\] 
In order to prove the claim we shall establish successively

\begin{itemize}
\item[1  -] The set ${\mathcal U}^\ep$ is path connected.
\item[2 -] The set ${\mathcal V}^\ep_C$ is closed in ${\mathcal U}^\ep$ for the $W^{1,p}-$topology
\item[3 -] For $\ep>0$ chosen small enough and $C>0$ large enough the set  ${\mathcal V}^\ep_C$ is open
for the $W^{1,p}-$topology
\end{itemize}
\medskip
Since ${\mathcal V}^\ep_C$ is non empty, this will imply the claim 1 for this choice of $\ep$ and $C$.
\medskip

\noindent{\it Proof of the path connectedness of ${\mathcal U}^\ep$.} 
For $A$ in ${\mathcal U}^\ep$ and $t\in [0,1]$ we define the connection form $A^t$ to be the image of $A$ by the dilation of rate $t^{-1}$ : $A^t=t\, \sum_{j=1}A_j(t\,x)\ dx_j$. We have in particular
\[
F_{A^t}= t^2\, \sum_{i,j=1}^4 (F_{A})_{ij}(x)\ dx_i\,dx_j
\]
hence
\[
\int_{B^4}|F_{A^t}|^2\ dx^4=\int_{B_{t}^4}|F_A|^2\ dx^4<\ep 
\]
 and 
\be
\label{III.9-a}
 \sum_{i,j=1}^4\int_{B^4}|\p_{x_i}(A_t)_{x_j}|^p\ dx^4= O(t^{2p-4})
\ee
This shows that  in one hand $A^t\in{\mathcal U}^\ep$ for any $t\in [0,1]$ and that in the other hand  $A^t\rightarrow 0$ strongly in $W^{1,p}(B^4)$. Hence $A^t$ is a continuous path contained in ${\mathcal U}^\ep$ connecting $A$ and $0$ which prove the path connectedness of ${\mathcal U}^\ep$.

\medskip

\noindent{\it Proof of the closeness of ${\mathcal V}^\ep_C$  in ${\mathcal U}^\ep$.}

Let $A^k\in {\mathcal V}^\ep_C$ and assume $A^k$ converges strongly in $W^{1,p}$ to some limit $A^\infty\in{\mathcal U}^\ep$. We claim that
$A^\infty\in {\mathcal V}^\ep_C$.

\medskip

Since $A^k\rightarrow A^\infty$ strongly in $W^{1,p}$, $dA^k \rightarrow dA^\infty$ strongly in $L^p$ and, using Sobolev embedding,
$A^k\rightarrow A^\infty$ in $L^{4p/4-p}$. Hence, due to the later, $$A^k\wedge A^k\longrightarrow A^\infty \wedge A^\infty\mbox{ strongly in }L^{2p/4-p}(B^4)\quad.$$
We have chosen $2<p<4$ in such a way that $p<2p/4-p$. Hence we deduce that
\be
\label{III.10}
F_{A^k}\longrightarrow F_{A^\infty}\quad\quad\mbox{ strongly in }L^p(B^4)\quad.
\ee
Let $g^k$ be a  sequence such that 
\be
\label{III.11}
 \int_{B^4}|d(A^k)^{g^k}|_{g_{B^4}}^q\, dvol_{S^4}\le C \int_{B^4}|F_{A^k}|^q\ dvol_{S^4}
\ee
for $q=2,p$ and
\be
\label{III.12}
d^\ast(A^k)^{g^k}=0\quad\mbox{ and }\quad \iota_{\p B^4}^\ast \ast\,(A^k)^{g^k}=0\quad.
\ee
Since both $d(A^k)^{g^k}$ and $d^\ast(A^k)^{g^k}$ are uniformly bounded in $L^p$ and since there is no harmonic 3-form\footnote{Assume $B$ is a 3 form satisfying $dB=0$ and $d^\ast B$ on $B^4$ and $\iota_{\p B^4}^\ast B=0$ then $\ast B=d\varphi$, $B=d\beta$ and $\iota_{\p B^4}^\ast\beta=dC$. With these notations we have
\[
\int_{B^4}|B|^2=\int_{B^4} d\beta\wedge d\varphi=\int_{\p B^4}dC\wedge d\varphi=0
\] 
hence $B=0$.} on $B^4$ whose restriction 
on the boundary is zero, classical $L^p$-Hodge theory  (see for instance \cite{IM})  infers that $(A^k)^{g^k}$ is uniformly bounded in $W^{1,p}$. So, using Sobolev embedding, it is in particular bounded in $L^{4p/4-p}$ and since $(A^k)^{g^k}=(g^k)^{-1}dg^k+(g^k)^{-1}A^k\,g^k$, using that $A^k$ is bounded in $W^{1,p}$ and hence in $L^{4p/4-p}$, we deduce that $dg^k$
is bounded also in $L^{4p/4-p}$. Since $4<4p/4-p$, there exists a subsequence $g^{k'}$  converging strongly to some limit $g^\infty$ in $C^0$. 
Going back to the weak convergence of $(A^k)^{g^k}$ in $W^{1,p}$, the strong convergence of $A^k$ in $W^{1,p}$ and  the weak convergence of $g^k$
in $W^{1,4p/4-p}$ we deduce that
\[
dg^k=g^k\, (A^k)^{g^k}-A^k\,g^k
\]
is uniformly bounded in ${W^{1,p}}$ and therefore $ g^{k'}\rightharpoonup g^\infty$ weakly in $W^{2,p}(B^4,G)$.
Thus we deduce that the following weak convergence in $W^{1,p}(B^4)$ holds
\[
(g^{k'})^{-1}dg^{k'}+(g^{k'})^{-1}A^{k'}\,g^{k'}\rightharpoonup (g^{\infty})^{-1}dg^{\infty}+(g^{\infty})^{-1}A^{\infty}\,g^{\infty}\quad.
\]
 Combining (\ref{III.10}) and the latest weak convergence 
we deduce that
\be
\label{III.13}
 \int_{B^4}|d(A^\infty)^{g^\infty}|^q\, dx^4\le C \int_{B^4}|F_{A^\infty}|^q\ dx^4
\ee
for $q=2,p$ and, using the  following continuous embedding of $$W^{1,p}(B^4)\hookrightarrow W^{1-1/p,p}(\p B^4)$$
we have
\be
\label{III.14}
d^\ast(A^\infty)^{g^\infty}=0\quad \mbox{ and }\quad \iota_{\p B^4}^\ast (A^\infty)^{g^\infty} 
\ee
 So we have proved that $A^\infty$ fulfills all the conditions for being in ${\mathcal V}^\ep_C$.

\medskip

\noindent{\it Proof of the openness of ${\mathcal V}^\ep_C$.}

Let $A$ be an element of ${\mathcal V}^\ep_C$. It is clear that if we find in ${\mathcal V}_C^\ep$ an open neighborhood for the $W^{1,p}-$topology of the $W^{1,p}$  Coulomb gauge $A^g$, then $A$ posses also such a neighborhood. So we can assume right away that $d^\ast A=0$ and $\iota_{\p B^4} \ast A=0$.

We are looking for the existence of $\delta$ sufficiently small - possibly depending on $A$ - such that for any $\omega$ satisfying $\|\om\|_{W^{1,p}}<\delta$ there exists $g$ close to
the identity in $W^{2,p}-$norm such that
\[
d^\ast\lf[g^{-1} dg+g^{-1} (A+\omega)\, g\rg]=0\quad\mbox{ and }\quad \iota_{\p B^4}\ast(A+\omega)^g=0
\]
To that purpose we introduce the map
\be
\label{III.14-a}
\begin{array}{l}
\ds{\mathcal N}^A\ :\ W^{1,p}(B^4,{\mathcal G})\times W^{2,p}(B^4,{\mathcal G})\, \longrightarrow \, L^{p}(B^4,{\mathcal G})\times W^{1-1/p,p}(\p B^4, {\mathcal G})\\[5mm]
\ds\quad (\omega,U)\ \longrightarrow\ \lf(d^\ast \lf[g_U^{-1}\, dg_U+g_U^{-1}(A+\omega)\, g_U\rg],\ \iota_{\p B^4}\ast(A+\omega)^{g_U}\rg)
\end{array}
\ee
where $g_U:=\exp(U)$.  We have seen that this map is smooth. 

The derivative of ${\mathcal N}^A$ along the $U$ direction at $(0,0)$ gives
\[
\p_U{\mathcal N}^A(0,0)\cdot V=\lf(-\Delta V+[A,dV], \p_r V\rg)\quad .
\]
where $\Delta=\sum_{k=1}^4 \p^2_{x_k^2}$. 

Using Calderon Zygmund $L^p$ theory (see for instance \cite{St} or \cite{GT}) we have the following a-priori estimate for any $V$ satisfying $\int_{B^4}V\ dx^4=0$
\[
\begin{array}{l}
\ds\|V\|_{W^{2,p}(B^4)}\le c\, \lf[\|\Delta V\|_{L^p(B^4)}+\|\p_r V\|_{W^{1-1/p,p}(\p B^4)}\rg]\\[5mm]
\ds\quad\le c\, \lf[\|\p_U{\mathcal N}^A(0,0)\cdot V\|_{{\mathcal F}}+\|[A,dV]\|_{L^p(B^4)}\rg]\\[5mm]
\ds\quad\le  c\, \lf[\|\p_U{\mathcal N}^A(0,0)\cdot V\|_{{\mathcal F}}+c\,\|A\|_{L^4(B^4)}\ \|dV\|_{L^{4p/4-p}(B^4)}\rg]
\end{array}
\]
where ${\mathcal F}$ is the hyperplane of $L^{p}(B^4,{\mathcal G})\times W^{1-1/p,p}(\p B^4, {\mathcal G})$ made of couples $(f,g)$ such that
\[
\int_{B^4}f(x)\ dx^4=-\int_{\p B^4} g(y)\ dvol_{\p B^4}
\]
From the fact that $A\in {\mathcal V}_C^\ep$ we deduce that  $\|A\|_{L^4}\le C_4\, \sqrt{C\ \ep}$ where $C_4$ is the Sobolev constant coming from the embedding
 into $L^4(B^4)$ of closed 3 forms on $B^4$ with adjoint exterior derivative in $L^2$ and whose restriction to $\p B^4$ is zero. Hence for any $V$ with
 average zero on $B^4$ we have 
\[
\|V\|_{W^{2,p}(B^4)}\le c\, [\|\p_U{\mathcal N}^A(0,0)\cdot V\|_{{\mathcal F}}+c\, C_4\, \sqrt{C\ \ep}\ \|dV\|_{L^{4p/4-p}}
\]

\noindent Using again the embedding of $W^{1,p}(\wedge^1B^4,{\mathcal G})$ into $L^{4p/4-p}(\wedge^1 B^4,{\mathcal G})$ and denoting $C_p$ the corresponding constant, we have then
\[
[1- c\, C_4\, \sqrt{C\ \ep}\ C_p]\ \|V\|_{W^{2,p}(B^4)}\le c\, \|\p_U{\mathcal N}^A(0,0)\cdot V\|_{{\mathcal F}}
\]
Having chosen $\ep$ such that  $c\, C_4\, \sqrt{C\ \ep}\ C_p<1/2$ we have that $\p_U{\mathcal N}^A(0,0)$ has zero Kernel. A classical result from Calderon-Zymund theory (see \cite{GT}) asserts that
\[
\begin{array}{l}
\ds{\mathcal L}\ :\ W^{2,p}(B^4,{\mathcal G})\longrightarrow {\mathcal F}\\[5mm]
\ds \quad V\longrightarrow (-\Delta V,\p_rV)
\end{array}
\]
is invertible and hence has zero index. By continuity of the index the maps
\[
\begin{array}{l}
\ds{\mathcal L}_t\ :\ W^{2,p}(B^4,{\mathcal G})\longrightarrow {\mathcal F}\\[5mm]
\ds \quad V\longrightarrow (-\Delta V+t[A,dV],\p_rV)
\end{array}
\]
has also zero index and since ${\mathcal L}_1=\p_U{\mathcal N}^A(0,0)$ has trivial kernel
it is invertible. So we can apply the implicit function  theorem and there exists $\delta>0$ together with an open neighborhood ${\mathcal O}$ of 0 in the subspace  of
 $W^{2,p}_0(S^4,{\mathcal G})$ with average 0 on $B^4$ such that 
 \[
 \begin{array}{l}
 \ds\forall \, \om\in W^{1,p}(B^4,{\mathcal G}) \quad\mbox{ satisfying }\quad  \|\om\|_{W^{1,p}(S^4,{\mathcal G})}<\delta\\[5mm]
 \ds\quad\exists\ !\ V_\omega\in {\mathcal O}\quad \mbox{ s. t.}\quad {\mathcal N}^A(V_\om,\om)=0\quad\mbox{ and}\quad\int_{B^4}V_\omega=0
  \end{array}
 \]
 and ${\mathcal O}$ can be taken smaller and smaller as $\delta$ tends to zero. We denote $g_\om:=\exp(V_\om)$.

 It remains to establish the control
 of the $L^p$ norm (resp. $L^2$ norm) of $d(A+\om)^{g_\om}$ in terms of the $L^p$ norm (resp. $L^2$ norm) of $F_{A+\om}$. 
 
 The Coulomb gauge
 $(A+\om)^{g_\om}$ satisfies for $q=2$ and $q=p$
 \[
 \|d(A+\om)^{g_\om}\|_{L^q}\le \|F_{A+\om}\|_{L^p}+ \|(A+\om)^{g_\om}\wedge(A+\om)^{g_\om}\|_{L^p}\\[5mm]
 \]
 We have $$\|(A+\om)^{g_\om}\|_{L^4}\le C_4\ \|dA\|_{L^2}+\|\om\|_{L^4}+\|dg^\om\|_{L^4}\quad.$$ Using the fact that $A$
 is the Coulomb gauge whose $W^{1,2}$ norm is controlled by the $L^2$ norm of $F_A$ - which is assumed itself to be less than $\ep$ - by taking $\ep$ small
 enough - independently of $A$ - by taking  $\delta$ small enough - depending possibly on $A$ - which ensures in particular that $\|dg_\om\|_{L^4}$ is sufficiently small, and by using the embedding of closed forms with $L^q$ exterior co-derivative  and whose restriction to $\p B^4$ is zero into $L^{4q/4-q}$ since $d^\ast(A+\om)^{g_\om}=0$ and $\iota_{\p B^4}^\ast (A+\om)^g=0$ we have established that
 \[
  \|d(A+\om)^{g_\om}\|_{L^q}\le \|F_{A+\om}\|_{L^q}+2^{-1} \|d(A+\om)^{g_\om}\|_{L^q}
 \]
 This implies that $A+\om$ fulfills the conditions for being in ${\mathcal V}^\ep_{C}$ for and $\omega$ satisfying $\|\om\|_{W^{1,p}}<\delta$ where $\delta$ has also been taken small enough in such a way that $\|F_{A+\om}\|<\ep$. This concludes the proof of the openness of ${\mathcal V}^\ep_{C}$ with respect to the $W^{1,p}-$topology for well chosen constants $\ep>0$ and $C>0$ and this concludes the proof of the {\bf claim}.
 
\medskip

\noindent{\it End of the proof of theorem~\ref{th-III.2}.}

\medskip

With the claim 1 at hand now, we are going to conclude the proof of theorem~\ref{th-III.2}

Let $A\in W^{1,2}(\wedge^1B^4,{\mathcal G})$ such that $\int_{S^4}|F_A|^2\ dvol_{S^4}<\ep$. Since $C^\infty$ is dense in $W^{1,2}$ there exists $A^t$  a family of smooth $1-$form on $S^4$ into ${\mathcal G}$
converging strongly to $A$ in $W^{1,2}$ as $t$ goes to zero. Using again the embedding of $W^{1,2}$ into $L^4$ we have the existence of $t_0>0$ such that
\[
\forall\ t<t_0\quad\quad\int_{B^4}|F_{A^t}|^2\ dvol_{B^4}<\ep\quad.
\]
Thus $A^t$ is in ${\mathcal U}^\ep$ and, due to the claim 1, it is also in ${\mathcal V}^\ep_C$. Let $g^t$ such that $d^\ast(A^t)^{g^t}=0$ with
\[
\ds \int_{B^4}\lf|d\lf((A^t)^{g^t}\rg)\rg|^2\, dx^4\le C \int_{B^4}|F_{A^t}|^2\ dx^4
\]
Again, since there is no non-trivial closed and co-closed $3-$form on $B^4$ the previous identity implies that $(A^t)^{g^t}$  is bounded in $W^{1,2}$ and then in $L^4$ too.The approximating connection one forms $A^t$ are converging strongly to $A$ in $W^{1,2}$ and hence in $L^4$ thus $d(g^t)$ is bounded in $L^4$.
We then deduce the existence of a sequence $t_k\rightarrow 0$ such that $g^{t_k}$ converges weakly in $W^{1,4}(B^4,G)$ to some limit $g^0$. 
Using Rellich Kondrachov compactness theorem $g^{t_k}$ converges strongly to $g^0$ in $L^p(B^4)$ for any $p<+\infty$ and hence $g^0$ is
also taking values in $G$ since we have almost everywhere convergence of the sequence.
 Using the previous convergences we deduce first
that
\[
(g^{t_k})^{-1}\, dg^{t_k}+(g^{t_k})^{-1} \, A\, g^{t_k}\rightharpoonup (g^0)^{-1}\, dg^0+(g^0)^{-1}\, A\, g^0\quad\mbox{ in }{\mathcal D}'(S^4)
\]
This implies that $d^\ast((A)^{g^0})=0$.  Since both $A^t$ and $(A^t)^{g^t}$ are bounded in $W^{1,2}$ and since $g^t$ is bounded in $L^4$, using 
\[
dg^t=g^t\, (A^t)^{g^t}-A^t\,g^t
\] 
we deduce that $g^t$ is bounded in $W^{2,2}$ and hence the trace of $g^{t_k}$ weakly converges to the trace of $g^0$ in $H^{3/2}(\p B^4,G)$. So we can pass to the limit in the equation $\iota_{\p B^4}^\ast\ast(A^{t_k})^{g^{t_k}}=0$ and we obtain
\[
\iota_{\p B^4}^\ast\ast(A^{0})^{g^0}=0
\]
Finally since $F_{A^t}$ is strongly converging to $F_A$ in $L^2$, using also the lower semi-continuity
of the $L^2$ norm together with the weak convergence of $(A^{t_k})^{g_{t_k}}$ towards $A^{g^0}$ we have
\[
\ds \int_{B^4}|d(A^{g^0})|_{g_{S^4}}^2\, dvol_{S^4}\le C \int_{S^4}|F_{A}|^2\ dvol_{S^4}\quad.
\]
This concludes the proof of theorem~\ref{th-III.2}.

\subsection{A refinement of Uhlenbeck's Coulomb Gauge extraction theorem.}
\medskip

This part can be skipped in a first reading.

\medskip

We have seen that Uhlenbeck's result is optimal in the sense that without assuming anything about the smallness of the Yang-Mills energy there
is no hope to obtain a gauge of $W^{1,2}$ controlled energy. One might wonder however if the smallness of the $L^2$ norm of the curvature
is the ultimate criterium for ensuring the existence of controlled Coulomb Gauges. The answer is ''no'' and one can very slightly reduce this 
requirement. Recall the notion of weak $L^2$ quasi-norm. We say that a measurable function $f$ on $B^m$ is in the {\it weak $L^2$ space} if
\[
|f|_{2,\infty}:=\lf[\sup_{\al>0}\ \al^2\ \lf|\lf\{x\in B^m\ ; \ |f(x)|>\al\rg\}\rg|\rg]^{1/2}<+\infty
\]
where $|\cdot |$ denotes the Lebesgue measure on $B^m$. This quantity defines a quasi-norm which is equivalent to a norm (see for instance \cite{Gr1}) that we denote $\|\cdot\|_{2,\infty}$. The {\it weak $L^2$ space} equipped with $\|\cdot\|_{2,\infty}$ is complete and define then a {\it Banach space} denoted $L^{2,\infty}$ called also {\it Marcinkiewicz weak $L^2$ space} or also {\it Lorentz weak $L^2$ space}. It is  larger than $L^2$. Indeed, for any function $f\in L^2$
we have
\[
\|f\|_{2,\infty}\le \sup_{\al>0}\int_{x\ ;|f|(x)>\al}|f|^2(x) dx\le \int |f|^2(x)\ dx=\|f\|_2^2
\]
It is strictly larger than $L^2$ : the function $f(x):=|x|^{-m/2}$ is in $L^{2,\infty}(B^m)$ but not in $L^2(B^m)$. It is also not difficult to see that
$L^{2,\infty}(B^m)\hookrightarrow L^p(B^m)$ for any $1\le p<2$.  
More generally we define the $L^{q,\infty}$ space of measurable functions $f$ satisfying
\[
|f|_{q,\infty}:=\lf[\sup_{\al>0}\ \al^q\ \lf|\lf\{x\in B^m\ ; \ |f(x)|>\al\rg\}\rg|\rg]^{1/q}<+\infty
\]
This defines again a quasi-norm equivalent to a norm\footnote{This is not true for $q=1$, the space $L^{1}-$weak cannot be made equivalent to a normed space - unfortunately, otherwise the analysis could make the economy of Calderon Zygmund theory and a major part of harmonic analysis that would sudently become trivial....!} if $q>1$. So it is a space which ''sits'' between $L^q(B^m)$ and all the $L^p(B^m)$ spaces for any $p<q$.
This is a space which has the same scaling properties as $L^q$ but  has however the big advantage of containing the Riesz functions $|x|^{-m/q}$ which play a central role in the theory of elliptic PDE. As we will see later the space $L^{q,\infty}(B^m)$ has also the advantage of being the dual of a Banach space, the {\it Lorentz space} $L^{q',1}(B^m)$ of measurable functions $f$ satisfying 
\be
\label{III.15}
\int_0^{+\infty}\lf|\lf\{x\ ;\ |f|(x)>\al\rg\}\rg|^{1/q'}\ d\al<+\infty
\ee
where $1/q'=1-1/q$ (see \cite{Gr1}). This later space has very interesting ''geometric'' properties that will be useful for the analysis of Yang-Mills Lagrangian as we will see below.

\medskip

We have the following theorem
\begin{Th}
\label{th-III.3}
Let  $m\le 4$ and $G$ be a compact Lie group. There exists $\ep_{G}>0$ and $C_{G}>0$ such that for any $A\in W^{1,2}(B^m,{\mathcal G})$ satisfying
\be
\label{III.8b}
\sup_{\al>0}\ \al^2\ \lf|\lf\{x\in B^m\ ; \ |F_A(x)|>\al\rg\}\rg|<\ep_G\quad,
\ee
there exists $g\in W^{2,2}(B^m,G)$ such that
\be
\label{III.9a}
\lf\{
\begin{array}{l}
\ds\int_{B^m}|A^g|^2+\sum_{i,j=1}^4|\p_{x_i}A^g_j|^2\, dx^m\le C_G\, \int_{B^m}|dA+A\wedge A|^2\, dx^m\\[8mm]
\ds d^\ast A^g=0\quad\quad\mbox{ in }B^m\\[5mm]
\ds\iota_{\p B^m}^\ast(\ast A^g)=0
\end{array}
\rg.
\ee
where $A^g=g^{-1}dg+g^{-1} A g$ and $\iota_{\p B^m}$ is the canonical inclusion map of the boundary of the unit ball into ${\R}^m$.
Moreover we have also
\[
\sum_{i,j=1}^4\|\p_{x_i}A^g_j\|^2_{2,\infty}\le C_G \ \|F_A\|_{2,\infty}^2\quad.
\]
\hfill$\Box$
\end{Th}
The weakening of the smallness criterium by replacing {\it small $L^2$}  by the less restrictive {\it small $L^{2,\infty}$}  condition for the existence of a controlled Coulomb gauge has been first
observed in \cite{BR}. This was a  very precious observation for the control of the loss of energies in so called {\it neck annular regions} in the study of conformally invariant problems such as {\it Willmore surfaces} or also {\it Yang-Mills Fields} as we will see below. The estimate (\ref{III.8b}) comes naturally
from the {\it $\epsilon-$regularity property} which holds in neck regions.

\medskip

\noindent{\bf Proof of theorem~\ref{th-III.3}.}

It follows exactly the same scheme as the proof of theorem~\ref{th-III.3} but we will need to use interpolation spaces between $L^{q,\infty}$ and $L^{q,1}$, the
{\it Lorentz spaces } $L^{q,s}$ and some of their properties.

Let $2<p<4$ and 
\[
\hat{\mathcal U}^\ep:=\lf\{A\in W^{1,p}(B^4,{\mathcal G})\quad\mbox{s.t. }\quad|F_{A}|^2_{L^{2,\infty}(B^4)}<\ep\rg\}
\]
and for any $\ep>0$ and $C>0$ and  we consider
\[
\hat{\mathcal V}^\ep_C:=\lf\{ 
\begin{array}{l}
\ds A\in \hat{\mathcal U}^\ep\quad\mbox{ s. t. }\quad\exists\, g\in W^{2,p}(B^4,G)\\[5mm]
\ds \int_{B^4}|dA^g|^p\, dx^4\le C \int_{B^4}|F_A|^p\ dx^4\\[5mm]
\ds \int_{B^4}|dA^g|^2\, dx^4\le C \int_{B^4}|F_A|^2\ dx^4\\[5mm]
\ds \|dA^g\|^2_{L^{2,\infty}(B^4)}\le C\ \|F_{A}\|^2_{L^{2,\infty}(B^4)}\\[5mm]
\ds\quad d^\ast A^g=0\quad\mbox{ and }\quad \iota_{\p B^4}\ast\,A^g=0
\end{array}
\rg\}
\]
The first goal is to show the following

\medskip

\noindent{\bf Claim}
\[
{ \exists\, \ep>0\quad C>0\quad\mbox{ s.t. }\quad\hat{\mathcal V}^\ep_C=\hat{\mathcal U}^\ep}
\] 

The proof of the claim is again divided in 3 steps.

\medskip

\noindent{\it Proof of the path connectedness of $\hat{\mathcal U}^\ep$.} 
For $A$ in ${\mathcal U}^\ep$ and $t\in [0,1]$ we define the connection form $A^t$ to be the image of $A$ by the dilation of rate $t^{-1}$ : $A^t=t\, \sum_{j=1}A_j(t\,x)\ dx_j$. We have in particular
\[
F_{A^t}= t^2\, \sum_{i,j=1}^4 (F_{A})_{ij}(x)\ dx_i\,dx_j
\]
hence $|F_{A^t}|(x)=t^2\ |F_{A}|(t\,x)$ and
\[
|F_{A_t}|_{L^{2,\infty}(B^4)}=|F_{A}|_{L^{2,\infty}(B_t^4)}\le |F_{A}|_{L^{2,\infty}(B^4)}<\ep
\]
and \footnote{This last inequality illustrates what we meant at the beginning of this subsection by {\it $L^{2,\infty}$  has the same scaling properties as $L^2$}.}
this path connects $A$ to $0$ in the $W^{1,p}$ topology due to (\ref{III.9-a}). Hence this concludes the proof of the path connectedness of $\hat{\mathcal U}^\ep$.

\medskip

\noindent The {\it proof of the closeness of $\hat{\mathcal V}^\ep_C$  in $\hat{\mathcal U}^\ep$} is identical to the proof of the closeness of ${\mathcal V}^\ep_C$  in ${\mathcal U}^\ep$.

\medskip

\noindent{\it Proof of the openness of $\hat{\mathcal V}^\ep_C$  in $\hat{\mathcal U}^\ep$.}

We consider the  map ${\mathcal N}^A$ defined by (\ref{III.14-a}). We recall\footnote{For a more thorough presentation of the Lorentz spaces and it's interaction with Calderon Zygmund theory in particular the reader is invited to consult the first chapter of \cite{Gr1} as well as \cite{SW} or \cite{Ta}.} the definition of the space $L^{q,s}(B^m)$ where $1<q<\infty$ and $1\le s<+\infty$.
A measurable function $f$ on $B^m$ belongs to $L^{q,s}(B^m)$ if
\be
\label{III.16}
|f|_{q,s}:=\lf[\int_0^\infty t^{\frac{s}{q}}\ f^\ast(t)\ \frac{dt}{t}\rg]^{1/s}<+\infty
\ee
where $f^\ast(t)$ is the decreasing rearrangement  function associated to $f$, defined on ${\R}_+$, and satisfying
\[
\forall \al>0   \quad\quad \lf|\lf\{t>0\ ;\ f^\ast(t)>\al\rg\}\rg|=\lf|\lf\{x\in B^m\ ;\ |f|(x)>\al\rg\}\rg|\quad.
\]
This defines again a quasi-norm equivalent to a norm for which the space is complete (see \cite{Gr1}). One verifies that the space $L^{q,1}(B^m)$ defined by (\ref{III.15}) coincides with the space given by (\ref{III.16}) for $s=1$. One verifies also
that $L^{q,q}(B^m)=L^q(B^m)$ and that for any $q\in (1,+\infty)$ and any $1\le s<\sigma\le +\infty$ we have $L^{q,s}(B^m)\hookrightarrow L^{q,\sigma}(B^m)$.
We have also that $\forall q<r$ and $\forall t,s\in [1,\infty]$ the following continuous embedding holds  $L^{p,s}(B^m)\hookrightarrow L^{q,t}(B^m)$.
The following multiplication rules holds and are continuous bilinear mappings in the corresponding spaces with the corresponding estimates
\be
\label{III.17}
L^{p,s}\cdot L^{q,t}\hookrightarrow L^{r,\sigma}
\ee
where $r^{-1}=p^{-1}+q^{-1}$ and $s^{-1}+t^{-1}=\sigma^{-1}$ and where $1<p,q<+\infty$ such that $r\ge 1$ and $1\le s,t\le \infty$ such that $1\le\sigma\le+\infty$. In particular we have for any $2\le p<4$
\be
\label{III.18}
L^{4,\infty}\cdot L^{\frac{4p}{4-p},p}\hookrightarrow L^p
\ee
Before to move on with the proof of theorem~\ref{th-III.3} we shall need a last tool from function theory : the {\it improved Sobolev embeddings} (see \cite{Ta}).
For $1\le p<m$ the following embedding is continuous
\be
\label{III.19}
W^{1,p}(B^m)\hookrightarrow L^{\frac{m\,p}{m-p},p}(B^m)
\ee
and more generally for any $t\in[1,+\infty]$
\be
\label{III.20}
W^{1,(p,t)}(B^m)\hookrightarrow L^{\frac{m\,p}{m-p},t}(B^m)
\ee
where $W^{1,(p,t)}(B^m)$ denotes the space of measurable functions on $B^m$ with distributional derivative in the Lorentz space $L^{p,t}(B^m)$.

\medskip

\noindent{\it Proof of openness of $\hat{\mathcal V}^\ep_C$ continued.} Using Calderon Zygmund $L^{q,t}$ theory we have the following bound
\[
\begin{array}{l}
\ds\|V\|_{W^{2,p}(B^4)}\le c\, \lf[\|\Delta V\|_{L^p(B^4)}+\|\p_r V\|_{W^{1-1/p,p}(\p B^4)}\rg]\\[5mm]
\ds\quad\le c\, \lf[\|\p_U{\mathcal N}^A(0,0)\cdot V\|_{{\mathcal F}}+\|[A,dV]\|_{L^p(B^4)}\rg]\\[5mm]
\ds\quad\le  c\, \lf[\|\p_U{\mathcal N}^A(0,0)\cdot V\|_{{\mathcal F}}+c\,\|A\|_{L^{4,\infty}(B^4)}\ \|dV\|_{L^{4p/4-p,p}(B^4)}\rg]
\end{array}
\]
where ${\mathcal F}:=W^{1,p}(B^4,{\mathcal G})\times W^{1-1/p,p}(\p B^4, {\mathcal G})$. From the fact that $A\in {\mathcal V}_C^\ep$ we deduce that  $\|A\|_{L^{4,\infty}}\le C_4\, \sqrt{C\ \ep}$ where $C_4$ is the Sobolev constant coming from the embedding
 into $L^4(B^4)$ of closed 3 forms on $B^4$ with adjoint exterior derivative in $L^2$ and whose restriction to $\p B^4$ is zero. Hence for any $V$ with
 average zero on $B^4$ we have 
\[
\|V\|_{W^{2,p}(B^4)}\le c\, [\|\p_U{\mathcal N}^A(0,0)\cdot V\|_{{\mathcal F}}+c\, C_4\, \sqrt{C\ \ep}\ \|dV\|_{L^{4p/4-p,p}}
\]

\noindent Using again the embedding (\ref{III.19}) and denoting $C_p$ the corresponding constant, we have then
\[
[1- c\, C_4\, \sqrt{C\ \ep}\ C_p]\ \|V\|_{W^{2,p}(B^4)}\le c\, \|\p_U{\mathcal N}^A(0,0)\cdot V\|_{{\mathcal F}}
\]
Having chosen $\ep$ such that  $c\, C_4\, \sqrt{C\ \ep}\ C_p<1/2$ we have that $\p_U{\mathcal N}^A(0,0)$, which is again a Fredholm operator of index zero, has a trivial kernel and is hence invertible.

The rest of the proof is completed by easily transposing to our present setting each argument of the case of small $L^2$ Yang-Mills energy which was detailed in the previous subsection.

\subsection{Controlled  gauges without small energy assumption.}

One might wonder why a $W^{1,2}$ control is wished and why one could not give up a bit our requirements and look for some
control of a ''weaker norm''.  
This is indeed possible, together with Mircea Petrache \cite{PR2}, the author proved the existence of global gauges $A^g$ whose $L^{4,\infty}$ norm is controlled by the Yang-Mills energy which is not necessarily small. Precisely we have. 
\begin{Th}
Let $(M^4,g)$ be a riemannian $4$-manifold. There exists a function $f:\mathbb R^+\to\mathbb R^+$ with the following properties.\\
 Let $\nabla$ be a $W^{1,2}$ connection over an $SU(2)$-bundle over $M$. Then there exists a \underbar{global}
$W^{1,(4,\infty)}$ section of the bundle (possibly allowing singularities) over the whole $M^4$ such that in the corresponding trivialization $\nabla$ is given by $d+A$ with the following bound.
\begin{equation*}
\|A\|_{L^{(4,\infty)}(M^4)}\leq f\left(\|F_{\nabla}\|_{L^2(M^4)}\right),
\end{equation*}
where $F_{\nabla}$ is the curvature form of $\nabla$.

\hfill $\Box$
\end{Th}

\section{The resolution of the Yang-Mills Plateau problem in the critical dimension.}
\reset
\subsection{The small energy case.}
We first present the resolution of the {\it Yang-Mills Plateau} problem in the case where the given connection at the boundary
has a small trace norm. Precisely we shall prove the following result.
\begin{Th}
\label{th-V.0}
Let $G$ be a compact Lie group and $m\le 4$. There exists $\delta_G>0$ such that for any 1-form $\eta\in H^{1/2}(\wedge^1\p B^m,{\mathcal G})$ satisfying
\be
\label{V.0}
\|\eta\|_{H^{1/2}(\p B^m)}<\delta_G\quad,
\ee
then the following minimization problem is achieved by a 1-form $A^0\in W^{1,2}(\wedge^1B^m,{\mathcal G})$
\[
\inf\lf\{ YM(A)=\int_{B^m}|dA+A\wedge A|^2\ dx^m\quad;\quad \iota_{\p B^m}^\ast A=\eta    \rg\}\quad.
\]
\hfill $\Box$
\end{Th}
The previous theorem is a corollary of the following weak closure theorem
\begin{Th}
\label{th-V.1}
Let $G$ be a compact Lie group and $m\le 4$. There exists $\delta_G>0$ such that for any 1-form $\eta\in H^{1/2}(\wedge^1\p B^m,{\mathcal G})$ satisfying
\be
\label{V.0-a}
\|\eta\|_{H^{1/2}(\p B^m)}<\delta_G\quad,
\ee
then for any $A^k\in W^{1,2}(\wedge^1B^m,{\mathcal G})$  satisfying
\[
\limsup_{k\rightarrow+\infty}YM(A^k)=\int_{B^m}|dA^k+A^k\wedge A^k|^2\ dx^m<+\infty\quad\mbox{ and }\quad\iota_{\p B^m}^\ast A^k=\eta
\]
there exists a subsequence $A^{k'}$ and a Sobolev connection $A^\infty\in W^{1,2}(\wedge^1B^m,{\mathcal G})$ such that
\[
D(A^{k'},A^\infty):=\inf_{g\in W^{2,2}(B^4,G)}\int_{B^m}|A^{k'}-(A^\infty)^g)|^2\ dx^m\longrightarrow 0
\]
moreover
\[
YM(A^\infty)\le\liminf_{k'\rightarrow 0} YM(A^{k'})\quad\mbox{ and }\quad\iota_{\p B^m}^\ast A^\infty=\eta\quad.
\]
\hfill $\Box$
\end{Th}
\noindent{\bf Proof of theorem~\ref{th-V.1}.}
We present the proof in the critical case $m=4$. The case $m<4$ being almost like the abelian linear case treated  . Let $B$ be the minimizer of $E$ in
 $W^{1,2}_\eta(\wedge^1B^4,{\mathcal G})$ and using (\ref{III.7}) and the Sobolev embedding $W^{1,2}(B^4)$ into $L^4(B^4)$   we have
 \be
\label{V.1}
\begin{array}{l}
\ds\quad \lf[\int_{B^4}|B|^4\ dx^4\rg]^{\frac{1}{4}}+\sum_{i,j=1}^m\int_{B^4}|\p_{x_i}B_j|^2\ dx^4\le C \lf[E(B)+\|\eta\|^2_{H^{1/2}}\rg]
\end{array}
\ee
The one form $B$ is the harmonic extension of $\eta$ and classical elliptic estimate gives
\be
\label{V.2}
E(B)\le C\,\|\eta\|_{H^{1/2}(\p B^4)}^2
\ee
Combining (\ref{V.1}) and (\ref{V.2})  we obtain the existence of a constant $C$ independent of $\eta$ such that
\[
\int_{B^4}|F_B|^2\ dx^4\le C\ \lf[\|\eta\|^2_{H^{1/2}(\p B^4)}+\|\eta\|^4_{H^{1/2}(\p B^4)}\rg]\quad.
\]
We choose first $\delta_G>0$ such that $C\ [\delta_G^2+\delta_G^4]< \ep_{G}$ in such a way that we can apply theorem~\ref{th-III.2} and we have
the existence of a minimizing sequence $A^k$ of $YM$ in $W^{1,2}_\eta(\wedge^1B^4,{\mathcal G})$ with a   Coulomb gauge $(A^k)^{g^k}$ controlled in $W^{1,2}$ :
\be
\label{V.3}
\|(A^k)^{g^k}\|_{W^{1,2}(B^4)}\le\ C\ \|F_{A^k}\|_{L^2(B^4)}\le C\ [\delta_G+\delta_G^2]
\ee
Without loss of generality
we can assume that
\[
(A^k)^{g^k}\rightharpoonup \hat{A}^\infty\quad\mbox{ weakly in }W^{1,2}(\wedge^1B^4,{\mathcal G})
\]
for some 1-form $\hat{A}^\infty$ which satisfies the Coulomb condition $d^\ast \hat{A}^\infty=0$ and for which
\be
\label{V.3-b}
\int_{B^4}|F_{\hat{A}^\infty}|^2\ dx^4\le \liminf_{k\rightarrow +\infty}\int_{B^4}|F_{(A^k)^{g^k}}|^2\ dx^4=\liminf_{k\rightarrow +\infty}\int_{B^4}|F_{A^k}|^2\ dx^4
\ee
 We claim that the restriction of $\hat{A}^\infty$ to $\p B^4$
is gauge equivalent to $\eta$. Because of the weak convergence of $(A^k)^{g^k}$ to $\hat{A}^\infty$  weakly in $W^{1,2}$, by continuity of the trace
operation from $W^{1,2}$ into $H^{1/2}$ we have
\be
\label{V.4}
\iota^\ast_{\p B^4}(A^k)^{g^k}= (g^k)^{-1}\, \iota_{\p B^4}^\ast dg^k +(g^k)^{-1}\,\eta\, g^k\rightharpoonup \iota_{\p B^4}^\ast\hat{A}^\infty\quad
\ee weakly in $H^{1/2}(\wedge^1\p B^4,{\mathcal G})$. Using the continuous embedding
\[
H^{1/2}(\p B^4)\hookrightarrow L^3(\p B^4)
\]
we have that the restriction of $g^k$ to $\p B^4$ converges weakly to some limit $g^\infty$ in $W^{1,3}(\p B^4)$ and we have, using (\ref{V.3}),
\be
\label{V.5}
\begin{array}{l}
\ds\|dg^\infty\|_{L^3(\p B^4)}\le \liminf_{k\rightarrow +\infty}\|dg^k\|_{L^3(\p B^4)}\\[5mm]
\quad\ds\le C\ \lf[\|\eta\|_{H^{1/2}(\p B^4)}+\liminf_{k\rightarrow +\infty}\|(A^k)^{g^k}\|_{H^{1/2}(\p B^4)}\rg]\\[5mm]
\quad\le C\ [\delta_G+\delta_G^2]
\end{array}
\ee
 Using now Rellich Kondrachov theorem (see for instance \cite{Br}), this convergence is strong in $L^q$ for any $q<+\infty$ which implies that $g^\infty$ takes values almost everywhere in $G$ and $g^\infty\in W^{1,3}(\p B^4,G)$.
We have moreover
\[
(g^\infty)^{-1}\,  dg^\infty +(g^\infty)^{-1}\,\eta\, g^\infty= \iota_{\p B^4}^\ast\hat{A}^\infty
\]
Using the continuous embedding 
\[
L^\infty\cap W^{1,3}(\p B^4)\cdot H^{1/2}(\p B^4)\hookrightarrow H^{1/2}(\p B^4)
\]
(the proof of this continuous embedding is also similar to the one of Lemma B1 in \cite{PR2}), we have that
\be
\label{V.6}
\begin{array}{l}
\ds\|d g^\infty\|_{H^{1/2}(\p B^4)}\\[5mm]
\ds\ \le C\ \lf[\|g^\infty\|_\infty+\|g^\infty\|_{W^{1,3}(\p B^4)}\rg] [\|\eta\|_{H^{1/2}(\p B^4)}+\|\iota^\ast_{\p B^4}\hat{A}^\infty\|_{H^{1/2}(\p B^4)}]\\[5mm]
\ds\ \le C\ [\delta_G+\delta_G^2]
\end{array}
\ee
We shall now make use of the following theorem which, as for theorem~\ref{th-III.1-a}
 can be proved following Uhlenbeck's Coulomb gauge extraction method.
 \begin{Th}
 \label{th-V.1-a}
 Let $G$ be a compact Lie group. There exists $\ep_{G}>0$ such that for any $g\in H^{3/2}(\p B^3,G)$ satisfying
 \[
 \|g\|_{H^{3/2}(\p B^4,G)}<\ep_G
 \]
 there exists an extension $\ti{g}\in W^{2,2}(B^4,G)$ of $g$ satisfying
 \[
 \|\ti{g}\|_{W^{2,2}(B^4,G)}\le C\ \|g\|_{H^{3/2}(\p B^4,G)}\quad.
 \]
 \hfill $\Box$
 \end{Th}
\noindent{\bf End of the proof of theorem~\ref{th-V.1}.} 
We choose $\delta_G$ small enough such that the r.h.s. of (\ref{V.6}) $C\ [\delta_G+\delta_G^2]$ is smaller than $\ep_G$ given by the previous theorem.
Let $\ti{g}^\infty\in W^{2,2}(B^4,G)$ be an extension of $g^\infty$ given by theorem~\ref{th-V.1-a}. Then 
$$
A^\infty:=(\hat{A}^\infty)^{(\ti{g}^\infty)^{-1}}\in W^{1,2}_\eta(\wedge^1 B^4,{\mathcal G})
$$
and we have using (\ref{V.3-b})
\[
\int_{B^4}|F_{A^\infty}|^2\ dx^4\le \liminf_{k\rightarrow +\infty}\int_{B^4}|F_{A^k}|^2\ dx^4\quad.
\]
Since $A^k$ is a minimizing sequence of the Yang-Mills Plateau problem in $W^{1,2}_\eta(\wedge^1 B^4,{\mathcal G})$, the connection form
$A^\infty$ is a solution to this problem and theorem~\ref{th-V.1} is proved.\hfill $\Box$
 
 \subsection{The general case and the point removability result for $W^{1,2}$ Sobolev connections. }

The theorem~\ref{th-V.1} as it is stated does not hold without the small norm assumption (\ref{V.0}) this is due to the fact that the theorem~\ref{th-V.1-a}
and similar results such as theorem~\ref{th-III.1-a}
do not hold for general data without smallness assumption (see again remark~\ref{rm-III.1-a}). We shall instead prove the following result where the boundary
condition is relaxed to a constrained trace modulo gauge action.
\begin{Th}
\label{th-V.2}
Let $G$ be a compact Lie group and $m\le 4$. For any 1-form $\eta\in H^{1/2}(\wedge^1\p B^m,{\mathcal G})$  the following minimization problem 
\be
\label{V.7}
\inf\lf\{\int_{B^m}|F_A|^2\ dx^m\ ;\ \iota_{\p B^m}^\ast A=\eta^g\quad\mbox{ for some }g\in H^{3/2}(\p B^4,G)    \rg\}
\ee
is achieved by a 1-form $A^0\in W^{1,2}(\wedge^1B^m,{\mathcal G})$.\hfill $\Box$
\end{Th}
In fact theorem~\ref{th-V.2} is a corollary of a general closure result.
\begin{Th}
\label{th-V.3}
For any compact Lie group $G$ and any dimension $m\le 4$, the space of Sobolev connections
\[
{\mathfrak A}_\eta(B^m):=\lf\{ A\in W^{1,2}(B^4,{\mathcal G})\ ; \ \iota_{\p B^m}^\ast A=\eta^g\quad\mbox{ for some }g\in H^{3/2}(\p B^4,G)    \rg\}
\] 
is weakly sequentially closed for sequences of controlled Yang-Mills energy. Precisely, for any $A^k\in  {\mathfrak A}_\eta(B^m)$ satisfying
\[
\limsup_{k\rightarrow+\infty}YM(A^k)=\int_{B^m}|dA^k+A^k\wedge A^k|^2\ dx^m<+\infty
\]
there exists a subsequence $A^{k'}$ and a Sobolev connection $A^\infty\in  {\mathfrak A}_\eta(B^m)$ such that
\[
d(A^{k'},A^\infty):=\inf_{g\in W^{1,2}(B^4,G)}\int_{B^m}|A^{k'}-(A^\infty)^g|^2\ dx^m\longrightarrow 0
\]
moreover
\[
YM(A^\infty)\le\liminf_{k'\rightarrow 0} YM(A^{k'})\quad.
\]
\hfill $\Box$
\end{Th}
\noindent{\bf Proof of theorem~\ref{th-V.3}.} Here again we restrict to the most delicate case : $m=4$. 

Let $A^k$ be a  sequence of ${\mathcal G}-$valued 1-forms and denote by $\ep_G$ the positive constant in Uhlenbeck's
theorem~\ref{th-III.2}\footnote{We choose in fact $\ep_G$ small enough for the controlled Gauge Uhlenbeck theorem to be valid for this constant
when the domain is any intersection of $B^4$ with a ball $B_\rho(y)$ for $y\in \ov{B^4}$ and  $0<\rho<1$}. A straightforward covering argument combined by some induction procedure gives the existence of a subsequence that we still denote $A^k$ and $N$
points $p_1\cdots p_N$ in $\ov{B^4}$ such that
\[
\begin{array}{l}
\ds\forall\, \delta>0 \quad\exists\,\rho_\delta>0\\[5mm]
\ds\sup_{k\in {\N}}\sup\lf\{\int_{B_{\rho_\delta}(y)}|F_{A^k}|^2\ dx^4\ ;\ y\in B^4\setminus\cup_{l=1}^NB_\delta(p_l)\rg\}<\ep_G
\end{array}
\]
\noindent{\it The case without concentration : $\{p_1\cdots p_N\}=\emptyset$.}

Let $\rho>0$ such that
\[
\sup_{k\in {\N}}\sup_{y\in B^4}\lf\{\int_{B_{\rho}(y)}|F_{A^k}|^2\ dx^4\rg\}<\ep_G
\]
We fix a finite {\it good covering~\footnote{The word ''good'' means that any intersections of elements of the covering is either empty or diffeomorphic to $B^4$ (see \cite{BT}).} } of $B^4$ by balls of  radius $\rho/2$. Denote $\{B_{\rho/2}(x_i)\}_{i\in I}$ this covering.
On each of the rice larger ball $B_{\rho}(x_i)$ for any $k\in {\N}$ we take a controlled Coulomb gauge $(A^k)^{g^k_i}$ such that
\be
\label{V.8}
\begin{array}{l}
\ds \lf[\int_{B_\rho^4(x_i)}|(A^k)^{g^k_i}|^4\ dx^4\rg]^{\frac{1}{4}}+\sum_{l,j=1}^m\int_{B_\rho^4(x_i)}|\p_{x_l}((A)^{g^k_i})_j|^2\ dx^4\\[5mm]
\ds\quad\quad\le\ C\int_{B_{\rho}(x_i)}|F_{A^k}|^2\ dx^4
\end{array}
\ee
and
\be
\label{V.9}
d^\ast(A^k)^{g^k_i}=0\quad.
\ee
For any pair $i\ne j$ in $I$ such that $B_\rho(x_i)\cap B_\rho(x_i)\ne \emptyset$ we denote $$g^k_{ij}:=g^k_i(g^k_j)^{-1}\in W^{2,2}(B_\rho(x_i)\cap B_\rho(x_i),G)\quad,$$
and we have in particular
\be
\label{V.10}
(A^k)^{g^k_j}=(g^k_{ij})^{-1}dg^k_{ij}+(g^k_{ij})^{-1}\, (A^k)^{g^k_{i}}\,g^k_{ij}\quad .
\ee
Observe that for any triplet $i\ne j$, $j\ne l$ and $i\ne l$ such that $B_\rho(x_i)\cap B_\rho(x_i)\cap B_\rho(x_l)\ne \emptyset$ we have the co-cycle
condition 
\be
\label{V.11}
\forall\ k\in {\N}\quad\quad g^k_{ij}\ g^k_{jl}=g^k_{il}\quad.
\ee
Combining (\ref{V.8}) and (\ref{V.10}) together with the improved Sobolev embedding $W^{1,2}(B^4)\hookrightarrow L^{4,2}(B^4)$ where $L^{4,2}$ is the Lorentz interpolation space given by (\ref{III.16})
we obtain that for any pair $i\ne j$ such that $B_\rho(x_i)\cap B_\rho(x_j)\ne \emptyset$
\be
\label{V.12}
\|dg^k_{ij}\|^2_{L^{4,2}( B_\rho(x_i)\cap B_\rho(x_j))}\le C\ \int_{B_{\rho}(x_i)\cup B_\rho(x_j)}|F_{A^k}|^2\ dx^4
\ee
From (\ref{V.10}) we have
\be
\label{V.13}
-\Delta g^k_{ij}=(A^k)^{g^k_{i}}\cdot dg^k_{ij}-dg^k_{ij}\cdot(A^k)^{g^k_j}
\ee
Using again the improved Sobolev embedding $W^{1,2}(B^4)\hookrightarrow L^{4,2}(B^4)$, inequalities (\ref{V.8}) and (\ref{V.12}) together with the continuous embedding 
\[
L^{4,2}\cdot L^{4,2}\hookrightarrow L^{2,1}
\]
we obtain
\be
\label{V.14}
\|\Delta g^k_{ij}\|_{L^{2,1}(B_\rho(x_i)\cap B_\rho(x_i))}\le C\ \int_{B_{\rho}(x_i)\cup B_\rho(x_j)}|F_{A^k}|^2\ dx^4
\ee
Using Calderon Zygmund theory in Lorentz interpolation spaces (see \cite{SW}) we obtain that $g^k_{ij}\in W^{2, (2,1)}_{loc}(B_\rho(x_i)\cap B_\rho(x_i))$ where $W^{2,(2,1)}$ denotes the space of functions
with two derivatives in $L^{2,1}$ and using (\ref{V.12}) together with (\ref{V.14}) we obtain the following estimate
\be
\label{V.15}
\|\nabla^2 g^k_{ij}\|_{L^{2,1}(B_{3\rho/4}(x_i)\cap B_{3\rho/4}(x_i))}\le C\ \lf[\int_{B_{\rho}(x_i)\cup B_\rho(x_j)}|F_{A^k}|^2\ dx^4\rg]^{1/2}
\ee
We can then extract a subsequence such that 
\[
\lf\{
\begin{array}{l}
\forall i\in I\quad (A^k)^{g^k_i}\rightharpoonup A^{i,\infty}\quad\mbox{ weakly in }W^{1,2}(B_\rho(x_i))\\[5mm]
\forall i\ne j\quad\quad g^k_{ij}\rightharpoonup g^\infty_{ij}\quad\mbox{ weakly in }W^{2,(2,1)}(B_{3\rho/4}(x_i)\cap B_{3\rho/4}(x_i)))
\end{array}
\rg.
\]
moreover $A^{i,\infty}$ and $g^\infty_{ij}$ satisfy the following identities
\be
\label{V.15a}
\lf\{
\begin{array}{l}
\ds \forall i\ne j\quad A^{j,\infty}=(g^\infty_{ij})^{-1}dg^\infty_{ij}+(g^\infty_{ij})^{-1}\, A^{i,\infty}\,g^\infty_{ij}\quad\\[5mm]
\ds  \forall i,j,l \quad g^\infty_{ij}\ g^\infty_{jl}=g^\infty_{il}\
\end{array}
\rg.
\ee
and we have the following estimate
\be
\label{V.16}
\|\nabla^2 g^\infty_{ij}\|^2_{L^{2,1}(B_{3\rho/4}(x_i)\cap B_{3\rho/4}(x_j))}\le C\ \liminf_{k\rightarrow +\infty}\int_{B_{\rho}(x_i)\cup B_\rho(x_j)}|F_{A^k}|^2\ dx^4
\ee
 It is proved in \cite{Ri0} that 
\[
W^{2, (2,1)}(B^4)\hookrightarrow C^0(B^4)
\]
hence we deduce that $g^\infty_{ij}\in C^0\cap W^{2,(2,1)}(B_{3\rho/4}(x_i)\cap B_{3\rho/4}(x_i))$ and for any $i\ne j$ there exists $\ov{g^\infty_{ij}}\in G$ such that
\be
\label{V.17}
\begin{array}{l}
\ds\|g^\infty_{ij}-\ov{g^\infty_{ij}}\|^2_{L^{\infty}(B_{3\rho/4}(x_i)\cap B_{3\rho/4}(x_j))}\\[5mm]
\ds\quad\quad\le C\ \liminf_{k\rightarrow +\infty}\int_{B_{\rho}(x_i)\cup B_\rho(x_j)}|F_{A^k}|^2\ dx^4<2\ C\ \ep_G
\end{array}
\ee
Taking $\ep_G$ small enough there exists a unique lifting $$U^\infty_{ij}\in W^{2,(2,1)}(B_{3\rho/4}(x_i)\cap B_{3\rho/4}(x_j))$$ such that
\[
\forall i\ne j\quad\quad g^\infty_{ij}=\ov{g^\infty_{ij}}\ \exp(U^\infty_{ij})
\]
and 
\[
\|U_{ij}^\infty\|_\infty\le C\ \ep_G
\]
for some constant $C$ depending only on $G$. Following an induction argument\footnote{A co-cycle smoothing argument by induction argument is also proposed in \cite{Is1} under the weaker hypothesis that the co-cycles $g^\infty_{ij}$ are $W^{1,4}$ in 4 dimension. This is made possible  due to the fact that $C^\infty(B^4,G)$ is dense in $W^{1,4}(B^4,G)$ (see \cite{SU}). The works of Takeshi Isobe \cite{Is1}, \cite{Is2} are proposing a framework for studying the analysis of gauge theory in conformal and super-critical dimension.} such as the one followed in \cite{MW} for the proof of theorem II.11, we can smooth the $U^\infty_{ij}$ in order to produce a sequence $$g^\infty_{ij}(t)\in C^\infty(B_{3\rho/4}(x_i)\cap B_{3\rho/4}(x_j),G)$$ satisfying
\[
g^\infty_{ij}(t)\longrightarrow g^\infty_{ij}\quad\quad\mbox{ strongly in }W^{2,(2,1)}(B_{3\rho/4}(x_i)\cap B_{3\rho/4}(x_j))\quad \mbox{ as } t\rightarrow 0
\]
and
\[
\forall\ t \quad\forall i,j,l\quad \quad g^\infty_{ij}(t)\ g^\infty_{jl}(t)=g^\infty_{il}(t)
\]
Since the ball $B^4$ is topologically trivial, the previous {\it cocycle condition} defines a trivial {\it \v{C}ech smooth co-chain} for the {\it presheaf} of $G-$valued smooth functions (see for instance \cite{BT} section 10 chapter II)
and for any $i\in I$ and any $t>0$ there exists $\rho_i(t)\in C^\infty(B_{3\rho/4}(x_i),G)$ such that
\be
\label{V.18}
g^\infty_{ij}(t)=\rho_i(t)\,\rho_j(t)^{-1}
\ee
We shall now make use of the following technical lemma which is proved in \cite{Uh2}.
\begin{Lm}
\label{lm-V.1}
Let $G$ be a compact Lie group and $\{U_i\}_{i\in I}$ be a good covering of $B^4$. There exists $\delta>0$ such that for any pair of co-chains
\[
\forall i\ne j\quad\quad h_{ij},g_{ij}\in W^{2,2}\cap C^0(U_i\cap U_j,G) 
\]
satisfying 
\[
\forall i,j,l\quad \quad g_{ij}\ g_{jl}=g_{il}\quad\mbox{ and }\quad h_{ij}\ h_{jl}=h_{il}\quad\mbox{ in }U_i\cap U_j\cap U_l.
\]
Assume
\[
\forall i\ne j\quad\quad \|g_{ij}^{-1}\,h_{ij}- e\|_{L^\infty(U_i\cap U_j)}<\delta
\]
where $e$ is the constant map equal to the neutral element of $G$, then, for any strictly smaller good covering of $B^4$ $\{V_i\}_{i\in I}$ satisfying $\ov{V_i}\subset U_i$, there exists a family of maps $\sigma_i\in W^{2,2}\cap C^0(U_i\cap U_j,G) $ such that
\[
\forall i\ne j\quad\quad h_{ij}=(\sigma_i)^{-1} g_{ij}\,\sigma_j\quad\quad\mbox{ in }V_i\cap V_j
\]
\hfill$\Box$
\end{Lm}
We apply the previous lemma to $h_{ij}:=g^\infty_{ij}(t)$ and $g_ij:=g^\infty_{ij}$ for $t$ small enough and we deduce the existence of 
\[
\sigma_i(t)\in W^{2,2}\cap C^0(B_{\rho/2}(x_i)\cap B_{\rho/2}(x_j))
\]
such that
\be
\label{V.19}
\forall i\ne j\quad\quad g^\infty_{ij}(t)=\sigma_{i}(t)^{-1}\, g^\infty_{ij}\,\sigma_j(t)\quad\quad\mbox{ in }B_{\rho/2}(x_i)\cap B_{\rho/2}(x_j)
\ee
Combining (\ref{V.18}) and (\ref{V.19}) we have
\[
\forall i\ne j\quad\quad g^\infty_{ij}=\sigma_i\rho_i\, (\sigma_j\rho_j)^{-1}\quad\quad\mbox{ in }B_{\rho/2}(x_i)\cap B_{\rho/2}(x_j)
\]
Combining this identity with (\ref{V.15a}) we set
\[
A^0:=(\sigma_i\rho_i)^{-1}\,  d(\sigma_i\rho_i)
+ (\sigma_i\rho_i)^{-1}\, A^{i,\infty}\, (\sigma_i\rho_i) \quad\mbox{ in }B_{\rho/2}(x_i)
\]
Clearly $A^0$ extends to a $W^{1,2}$ ${\mathcal G}-$valued 1-form in $B^4$, moreover, following the arguments in the proof of theorem~\ref{th-V.1}, the restriction of $A^0$ to $\p B^4$ is gauge equivalent to $\eta$. This concludes the proof of theorem~\ref{th-V.3} in the absence of concentration points.

\medskip

\noindent{\it The general case with possible concentration : $\{p_1\cdots p_N\}\ne\emptyset$.}

Following the arguments in the previous case, for any $\delta>0$ we exhibit a subsequence $A^{k'}$, a covering by balls $B_{\rho_\delta}(x_i)$
of $B^4\setminus \cup_{l=1}^NB^4_\delta(p_l)$ and a family of gauge changes $g_i^k$  such that
\[
\lf\{
\begin{array}{l}
\forall i\in I\quad (A^k)^{g^k_i}\rightharpoonup A^{i,\infty}\quad\mbox{ weakly in }W^{1,2}(B_\rho(x_i))\\[5mm]
\forall i\ne j\quad\quad g^k_{ij}\rightharpoonup g^\infty_{ij}\quad\mbox{ weakly in }W^{2,(2,1)}(B_{3\rho/4}(x_i)\cap B_{3\rho/4}(x_i)))
\end{array}
\rg.
\]
The family $g^\infty_{ij}$ defines again a $W^{2,(2,1)}-$co-chain that we can approximate in $C^0\cap W^{2,2}$ by a smooth one $g^\infty_{ij}(t)$. Using the fact that the second homotopy group of the compact Lie group is trivial $\pi_2(G)=0$ (see for instance \cite{BrDi} chapter V proposition 7.5) we deduce that the co-chain $g^\infty_{ij}(t)$, defined on a covering of $B^4\setminus \cup_{l=1}^NB^4_\delta(p_l)$. is trivial for the \v{C}ech cohomology for the co-chains on the pre-sheaf of smooth $G-$valued functions. Following each step of the above argument we construct a  $W^{1,2}$ ${\mathcal G}-$valued 1-form $A^0$ in $B^4\setminus \cup_{l=1}^NB^4_\delta(p_l)$ which is gauge equivalent to  $A^{i,\infty}$ in $B^4_{\rho/2}(x_i)$ for each $i\in I$ and whose restriction on $\p B^4\setminus \cup_{l=1}^NB^4_\delta(p_l)$ is also gauge equivalent
to $\eta$. Moreover we have
\[
\int_{B^4\setminus \cup_{l=1}^NB^4_\delta(p_l)}|F_{A^0}|^2\ dx^4\le \liminf_{k\rightarrow +\infty}\int_{B^4}|F_{A^k}|^2\ dx^4\quad.
\]
Using a diagonal argument with $\delta\rightarrow 0$ we can extend  $A^0$ as a ${\mathcal G}-$valued 1-form in $W^{1,2}_{loc}(B^4\setminus\{p_1\cdots p_N\})$ and still satisfying
\be
\label{V.19a}
\int_{B^4}|F_{A^0}|^2\ dx^4\le \liminf_{k\rightarrow +\infty}\int_{B^4}|F_{A^k}|^2\ dx^4\quad.
\ee
We conclude the proof of theorem~\ref{th-V.3} by changing the gauge of $A^0$ in the neighborhood of each blow up point $p_l$ making use of the following theorem~\ref{V.4}, known as {\it point removability theorem}, which gives the existence of a change of gauge $g$ in order to extend our connection 1-form $(A^0)^g$ as a $W^{1,2}$ ${\mathcal G}$ valued $1-$form in the neighborhood of each $p_l$. We then paste together these $W^{1,2}-$ gauges by using the same technique as the one we used in the case without blow up points
in order to get a global $W^{1,2}$ representative of $A^0$ on $B^4$ gauge equivalent to $\eta$ on $\p B^4$ and satisfying (\ref{V.19a}). This concludes the proof of theorem~\ref{th-V.3}. \hfill $\Box$
\begin{Th}
\label{th-V.4}
{\bf[Point removability]}
Let $G$ be a compact Lie group and $A\in W^{1,2}_{loc}(\wedge^1B^4,{\mathcal G})$ such that
\[
\int_{B^4}|dA+A\wedge A|^2\ dx^4<+\infty
\]
then there exists a gauge change $g\in W^{2,2}_{loc}(B^4,G)$ such that
\[
A^g\in W^{1,2}(\wedge^1B^4,{\mathcal G})\quad.
\]
\hfill $\Box$
\end{Th}
\begin{Rm}
\label{rm-V.5}
Point removability results play an important r\^ole in the analysis of conformally invariant variational problems. This is a natural consequence of due to the existence of point concentration which is inherent to the conformal invariance. These results are often formulated for the critical points of conformally invariant lagrangians 
and in the present case it has been first proved by K.Uhlenbeck for Yang-Mills fields (see \cite{Uh1}). Observe that here we are not assuming that $A$ is satisfying a particular equation. \hfill $\Box$
\end{Rm}
\begin{Rm}
\label{rm-V.6}
Beyond geometric analysis, point removability results play also an important r\^ole in complex geometry. One could for instance quote the work of Bando 
\cite{Ban} about the possibility to extend an hermitian holomorphic structure $F_{A}^{0,2}=0$ with $L^2$ bounded curvature on a the punctured ball $B^4\setminus\{0\}$ as a smooth holomorphic bundle throughout the origin. Beside the holomorphicity condition $F^{0,2}_A=0$ no further ''equation'' is assumed and in particular the Einstein equation $\omega\cdot F_A^{1,1}=c\, I$ is not assumed and  the connection form is not necessarily a Yang-Mills field.\hfill $\Box$
\end{Rm}
\noindent{\bf Proof of theorem~\ref{th-V.4}.}
Without loss of generality we can assume that
\[
\int_{B^4}|F_A|^2\ dx^4<\delta
\]
where $\delta>0$ will be fixed later on in the proof. Denote for $i\ge 2$ $$T_i:=B^4_{2^{-i+2}}(0)\setminus B^4_{2^{-i-2}}(0)\quad.$$
 From theorem~\ref{th-III.2}     there exists $\delta>0$ such that, on each annulus $T_i$ there exists a change of gauge $g_i$ such that 
there exists $g_i\in W^{2,2}(T_i,G)$ such that
\be
\label{V.20}
\lf\{
\begin{array}{l}
\ds\int_{T_i}2^{2\,i}|A^{g_i}|^2+\sum_{k,l=1}^4|\p_{x_k}A^{g_i}_l|^2\, dx^m\le C_G\, \int_{T_i}|dA+A\wedge A|^2\, dx^4\\[8mm]
\ds d^\ast A^{g_i}=0\quad\quad\mbox{ in }T_i\\[5mm]
\ds\iota_{\p B^m}^\ast(\ast A^{g_i})=0
\end{array}
\rg.
\ee
On $T_i\cap T_{i+1}=B^4_{2^{-i+1}}\setminus B^4_{2^{-i-2}}$ the transition function $g_{i \,i+1}=g_i(g_{i+1})^{-1}$ satisfy
\be
\label{V.21}
A^{g_{i+1}}=(g_{i\,i+1})^{-1}\, dg_{i\,i+1}+(g_{i\,i+1})^{-1}\, A^{g_{i}}\,g_{i\,i+1}
\ee
Hence (\ref{V.20}) imply 
\be
\label{V.22}
2^{2\,i}\int_{B^4_{2^{-i+1}}\setminus B^4_{2^{-i-2}}}|dg_{i\, i+1}|^2\ dx^4\le C\  \int_{B_{2^{-i+2}}\setminus B_{2^{-i-3}}}|dA+A\wedge A|^2\, dx^4
\ee
Taking the adjoint of the covariant derivative of equation (\ref{V.21})
\[
-\Delta g_{i\, i+1}=A^{g_i}\cdot dg_{i\,i+1}-dg_{i\,i+1}\cdot A^{g_{i+1}}
\]
and, arguing as in the first part of the proof of theorem~\ref{th-V.3}, we deduce the existence of $\ov{g_{i\, i+1}}\in G$ such that
\be
\label{V.23}
\begin{array}{l}
\ds\|g_{i\, i+1}-\ov{g_{i\, i+1}}\|_{L^\infty(T_i\cap T_{i+1})}\\[5mm]
\ds\le C\ 2^{i}\|dg_{i\,i+1}\|_{L^2(T_i\cap T_{i+1})}+C\,\sum_{k,l=1}^4\|\p^2_{x_k x_l} g_{i\, i+1}\|_{L^{2,1}(T_i\cap T_{i+1})}\\[5mm]
\ds\le C\ \lf[  \int_{T_i\cup T_{i+1}}|dA+A\wedge A|^2\, dx^4 \rg]^{1/2}\le C\ \sqrt{\delta}
\end{array}
\ee
We now modify the gauge change $g_i$ as follows. Precisely, for any $i\in {\N}$, we denote
\[
\ov{\sigma_i}:=\ov{g_{1\,2}}\  \ov{g_{2\, 3}}\, \cdots \ov{g_{i-1\, i}}\in G\quad.
\] 
Observe that 
\be
\label{V.23aa}
A^{g_i\,\ov{\sigma_i}^{-1}}=\ov{\sigma_i}\, A^{g_i}\, \ov{\sigma_i}^{-1}
\ee
Hence $A^{g_i\,\ov{\sigma_i}^{-1}}$ is still a Coulomb gauge satisfying 
\be
\label{V.23a}
\lf\{
\begin{array}{l}
\ds\int_{T_i}2^{2\,i}|A^{g_i\,\ov{\sigma_i}^{-1}}|^2+\sum_{k,l=1}^4|\p_{x_k}A^{g_i\,\ov{\sigma_i}^{-1}}_l|^2\, dx^m\le C_G\, \int_{T_i}|dA+A\wedge A|^2\, dx^4\\[8mm]
\ds d^\ast A^{g_i\,\ov{\sigma_i}^{-1}}=0\quad\quad\mbox{ in }T_i\\[5mm]
\ds\iota_{\p B^m}^\ast(\ast A^{g_i\,\ov{\sigma_i}^{-1}})=0
\end{array}
\rg.
\ee
Denote $h_i:=g_i\,\ov{\sigma_i}^{-1}$ the transition functions on $T_i\cap T_{i+1}$ for these new gauges are given by
\[
h_{i\,i+1}:=g_i\,\ov{\sigma_i}^{-1}\,\ov{\sigma_{i+1}}\, (g_{i+1})^{-1}=g_i\ \ov{g_{i\,i+1}}\, (g_{i+1})^{-1}
\]
Using (\ref{V.23}) we have
\be
\label{V.23b}
\|h_{i\,i+1}-e\|_{L^{\infty}(T_i\cap T_{i+1})}\le C\ \lf[  \int_{T_i\cup T_{i+1}}|dA+A\wedge A|^2\, dx^4 \rg]^{1/2}\le C\ \sqrt{\delta}
\ee
Exactly as for $g_i$, using the identity 
\be
\label{V.21a}
A^{h_{i+1}}=(h_{i\,i+1})^{-1}\, dh_{i\,i+1}+(h_{i\,i+1})^{-1}\, A^{h_{i}}\,h_{i\,i+1}
\ee
together with (\ref{V.23aa}) and (\ref{V.23a}) we obtain
\be
\label{V.25}
\begin{array}{l}
 \ds 2^{i}\|dh_{i\,i+1}\|_{L^2(T_i\cap T_{i+1})}+\,\sum_{k,l=1}^4\|\p^2_{x_k x_l} h_{i\, i+1}\|_{L^{2,1}(T_i\cap T_{i+1})}\\[5mm]
\ds\le C\ \lf[  \int_{T_i\cup T_{i+1}}|dA+A\wedge A|^2\, dx^4 \rg]^{1/2}\le C\ \sqrt{\delta}
\end{array}
\ee
where $e$ is the content function on $T_i\cap T_{i+1}$ equal to the neutral element of $G$. Having chosen $\delta$ small enough we ensure that the transition functions of this new set of trivialization
are contained in a neighborhood of the neutral element into which the exponential map defines a diffeomorphism and there exist $U_{i\, i+1}$ such that $h_{i\,i+1}= \exp(U_{i\, i+1})$ and
\be
\label{V.26}
\begin{array}{l}
\ds \|U_{i\,i+1}\|_{L^\infty(T_i\cap T_{i+1})}+2^{i}\|dU_{i\,i+1}\|_{L^2(T_i\cap T_{i+1})}+\,\sum_{k,l=1}^4\|\p^2_{x_k x_l} U_{i\, i+1}\|_{L^{2,1}(T_i\cap T_{i+1})}\\[5mm]
\ds\le C\ \lf[  \int_{T_i\cup T_{i+1}}|dA+A\wedge A|^2\, dx^4 \rg]^{1/2}\le C\ \sqrt{\delta}
\end{array}
\ee
Let $\rho$ be a smooth function on ${\R}_+$ identically equal to 1 between 0 and $\sqrt{2}$ and compactly supported in $[0,2]$. On $B^4$ we define
\[
\rho_i(x):=\rho(|x|\, 2^{i})\quad V_i:=B_{2^{-i+3/2}}\setminus B_{2^{-i}}\quad\mbox{ and }\tau_i:=\exp(\rho_i\, U_{i\, i+1})
\]
With these notations we have
\[
\mbox{ on }\quad V_{i+1}\cap V_i=B_{2^{-i+1/2}}\setminus B_{2^{-i}}\quad\mbox{ we have }\tau_i=h_{i\, i+1}\mbox{ and }\tau_{i+1}=e
\]
Hence on $V_{i+1}\cap V_i$ we have
\be
\label{V.27}
\begin{array}{l}
\ds A^{h_i\,\tau_i}=\tau_i^{-1}\, d\tau_i+\tau_i^{-1}\,A^{h^i}\,\tau_i\\[5mm]
\ds\quad=(h_{i\,i+1})^{-1}\, dh_{i\,i+1}+(h_{i\,i+1})^{-1}\, A^{h_{i}}\,h_{i\,i+1}=A^{h_{i+1}}=A^{h_{i+1}\,\tau_{i+1}}
\end{array}
\ee
and the 1-form $\hat{A}$ equal to $A^{h_i\,\tau_i}$ on each annulus $V_i$  defines a global $W^{1,2}_{loc}$ connection 1-form on $B^4\setminus \{0\}$ gauge equivalent to $A$. Clearly, for $k=1,2$, we have the pointwise estimate
\[
|d^k\tau_i|\le C\ \sum_{l=0^k} 2^{i\, l}\ |d^{k-l}U_{i\,i+1}|\quad\mbox{ on }V_i
\]
Combining this fact together with (\ref{V.23a}), (\ref{V.26}) and (\ref{V.27}) we obtain
\[
\int_{V_i}2^{2\,i}|A^{h_i\,\tau_i}|^2+\sum_{k,l=1}^4|\p_{x_k}A^{h_i\,\tau_i}_l|^2\, dx^4\le C_G\, \int_{T_i\cup T_{i+1}}|dA+A\wedge A|^2\, dx^4
\]
Summing over $i$ gives
\[
\int_{B^4}|x|^{-2}\ |\hat{A}|^2+\sum_{k,l=1}^4|\p_{x_k}\hat{A}_l|^2\, dx^4\le C_G\, \int_{B^4}|dA+A\wedge A|^2\, dx^4
\]
$\hat{A}$ is then in $W^{1,2}(\wedge^1B^4,{\mathcal G})$ and this concludes the proof of theorem~\ref{th-V.4}.\hfill $\Box$

\section{The Yang-Mills equation in sub-critical and critical dimensions.}

\subsection{Yang-Mills fields.}
Until now we have produced solutions to the Yang-Mills Plateau problem  in dimension less or equal to 4 but we have not addressed issues related to the special properties that should be satisfied by these solutions. Maybe one
of the first question that should be looked at is whether these minima define smooth {\it equivariant horizontal plane distributions} or not.

In order to study the regularity of solutions to the Yang-Mills Plateau problem we have first to produce the {\it Euler Lagrange equation} attached to this variational problem. This is the so called {\it Yang-Mills equation}.
\begin{Dfi}
\label{df-VI.1}
Let $G$ be a compact Lie group and $A$ be an $L^2$ connection 1-form on $B^m$ into the Lie Algebra ${\mathcal G}$ of $G$. Assume that
\[
\int_{B^m}|dA+A\wedge A|^2\ dx^m<+\infty
\]
we say that $A$ is a Yang-Mills field if
\[
\begin{array}{l}
\ds\forall\,\xi\in C^\infty_0(\wedge^1 B^m,{\mathcal G})\\[5mm]
\ds\quad\quad\frac{d}{dt}\int_{B^m}|d(A+t\xi)+(A+t\xi)\wedge (A+t\xi)|^2\ dx^m\lf|_{t=0}\rg.=0\quad.
\end{array}
\]
\hfill $\Box$
\end{Dfi}
Observe that this definition makes sense for any $A\in L^2$ such that $F_A\in L^2$, indeed we have for any $\xi$ in $C^\infty_0(\wedge^1 B^m,{\mathcal G})$
\[
F_{A+t\xi}=F_A+t\ \lf(d\xi+A\wedge\xi+\xi\wedge A\rg)+t^2\ \xi\wedge\wedge\xi\quad\in L^2(\wedge^2B^m,{\mathcal G})
\]
For such a $A\in L^2$ and for any $\xi$ in $C^\infty_0(\wedge^1 B^m,{\mathcal G})$ we denote by $d_A\xi$ the following 2-form
\[
d_A\xi(X,Y):=d\xi(X,Y)+[A(X),\xi(Y)]+[\xi(X),A(Y)]
\]
So we have for instance
\[
d_A\xi(\p_{x_i},\p_{x_j})=\p_{x_i}\xi_j-\p_{x_j}\xi_i+[A_i,\xi_j]+[\xi_i,A_j]
\]
We have then the following proposition.
\begin{Prop}
\label{pr-VI.1}{\bf[Yang-Mills Equation]}
Let $A\in L^2(\wedge^1 B^m,{\mathcal G})$ such that $F_A\in L^2(\wedge^2 B^m,{\mathcal G})$. The connection 1-form $A$ is a Yang-Mills field if
\be
\label{VI.1}
\begin{array}{l}
\ds\forall\,\xi\in C^\infty_0(\wedge^1 B^m,{\mathcal G})\\[5mm]
\ds\quad\quad\int_{B^4}d_A\xi\cdot F_A=0
\end{array}
\ee
which is equivalent to
\be
\label{VI.2}
d_{A}^\ast F_A=0\quad\mbox{ in }{\mathcal D}'(B^m)
\ee
In coordinates this reads
\be
\label{VI.3}
\forall\ i=1\cdots m\quad\sum_{j=1}^m\p_{x_j}(F_{A})_{ij}+[A_j,(F_{A})_{ij}]=0\quad.
\ee
\hfill $\Box$
\end{Prop}
The Yang-Mills equation (\ref{VI.1}) is also written symbolically as follows
\[
d^\ast F_A+[A,\res F_{A}]=0
\]
where $\res$ is referring to the contraction operation between tensors with respect to the flat metric on $B^m$.
The proof of the last statement of the proposition goes as follows (\ref{VI.1}) in coordinates is equivalent to
\[
\sum_{i,j=1}^m\int_{B^m}\lf<\p_{x_i}\xi_j-\p_{x_j}\xi_i+[A_i,\xi_j]+[\xi_i,A_j],(F_A)_{ij}\rg>\ dx^m=0
\]
using integration by parts and the fact that the Killing metric, invariant under adjoint action, satisfies $<U,[V,W]>=<W,[U,V]>$ we obtain
\[
\sum_{i,j=1}^m\int_{B^m}\lf<-\p_{x_i}(F_A)_{ij}+[(F_A)_{ij},A_i],\xi_j\rg>+\lf<\p_{x_j}(F_A)_{ij}+[A_j,(F_A)_{ij}],\xi_i\rg> dx^m=0
\]
which implies (\ref{VI.3}).

\medskip

The gauge invariance of the integrant of Yang-Mills lagrangian implies that (\ref{VI.3}) is solved for $A$ if and only if it is solved
for any gauge transformation $A^g$ of $A$. More generally we have the following 
\be
\label{VI.3-sst}
\lf\{
\begin{array}{l}
\ds\forall\,A\in W^{1,2}(\wedge^1B^m,{\mathcal G})\quad\quad \forall\ g\in W^{2,2}(B^m,G)\\[5mm]
\ds d_{(A^g)}^\ast F_{A^g}=g^{-1}\,d_A^\ast F_A\, g
\end{array}
\rg.
\ee
where we recall that $A^g:=g^{-1}dg+g^{-1}A\, g$.

\medskip

The Yang-Mills equation (\ref{VI.2}) has to be compared with the {\it Bianchi identity} to which it is a kind of ''dual equation''. This is a {\it structure equation} which holds for any
connection 1-form.

\begin{Prop}
\label{pr-VI.2}{\bf [Bianchi identity]}
For any $A\in L^2(\wedge^1 B^m,{\mathcal G})$ such that $F_A\in L^2(\wedge^2 B^m,{\mathcal G})$ the following identity holds
\[
d_AF_A=0
\] 
where $d_AF_A$ is the 3-form given by
\[
\begin{array}{l}
\ds d_AF_A(X,Y,Z):=dF_A(X,Y,Z)\\[5mm]
\ds\quad\quad+[A(X),F_A(Y,Z)]+[A(Y),F_A(Z,X)]+[A(Z),F_A(X,Y)]
\end{array}
\]
\hfill $\Box$
\end{Prop}
The proof of {\it Bianchi identity} goes as follows. We have
\[
\begin{array}{l}
\ds dF_A(X,Y,Z)=d(A\wedge A)(X,Y,Z)\\[5mm]
\ds\quad=[dA(X,Y), A(Z)]+[dA(Y,Z), A(X)]+[dA(Z,X), A(Y)]\\[5mm]
\ds\quad=[F_{A}(X,Y),A(Z)]+[F_A(Y,Z), A(X)]+[F_A(Z,X), A(Y)]
\end{array}
\]
where we have used the {\it Jacobi identity}
\[
[[A(X),A(Y)],A(Z)]+[[A(Y),A(Z)],A(X)]+[[A(Z),A(X)],A(Y)] =0
\]
This concludes the proof of {\it Bianchi identity}.

\medskip

In the particular case where $G$ is {\bf abelian}, Yang-Mills equation together with Bianchi identity is equivalent to the harmonic map form system
\[
\lf\{
\begin{array}{l}
d^\ast F_A=0\quad\quad\mbox{\bf Yang-Mills }\\[5mm]
d F_A=0\quad\quad\mbox{\bf Bianchi }
\end{array}
\rg.
\]
whose solutions are known to be analytic in every dimension. We are now asking about the same regularity issue in the \underbar{non abelian} case.

\subsection{The regularity of $W^{1,2}$ Yang-Mills Fields in  sub-critical and critical dimensions.}

The Yang-Mills equation (\ref{VI.3})  is impossible to exploit for proving any regularity. Indeed, Since the associated lagrangian is invariant under gauge transformation
the equation is also invariant under this action of this huge group and assuming any regularity for $A$ would be established then, tacking any arbitrary other non smooth gauge $g$
there is no reason why $A^g$ should be again smooth though it still solves Yang-Mills. here again breaking the gauge invariance plays a condimental role and we are proving the following
result
\begin{Th}
\label{th-VI.1}
Let $G$ be a compact Lie group and $m\le 4$ and let $A\in W^{1,2}(\wedge^1 B^m,{\mathcal G})$ satisfying the Yang-Mills equation (\ref{VI.2}) then for any $W^{1,2}$
Coulomb gauge $A^g$ (i.e. satisfying $d^\ast A^g=0$), $A^g$ is $C^\infty$.\hfill $\Box$
\end{Th}
\noindent{\bf Proof of theorem~\ref{th-VI.1}.}

We assume that $A$ is Coulomb and satisfy the Yang-Mills equation (\ref{VI.2}). Hence we have
\[
d^\ast dA+d^\ast(A\wedge A)+[A,\res dA]+[A,\res A\wedge A]=0\quad.
\]
Using the fact that $d^\ast A=0$ the Yang-Mills equation in this Coulomb gauge reads then
\be
\label{VI.3aaa-0}
\Delta A=d^\ast(A\wedge A)+[A,\res dA]+[A,\res (A\wedge A)]\quad,
\ee
and theorem~\ref{th-VI.1} is now the direct consequence of the following result.\hfill $\Box$
\begin{Th}
\label{th-VI.2}
Let $m\le 4$ and $N\in {\N}^\ast$. Let $f\in C^\infty({\R}^N\times({\R}^m\otimes{\R}^N),{\R}^N)$ and let $g\in C^\infty({\R}^N,{\R}^N)$ such that there exists $C>0$
satisfying
\be
\label{VI.3aa}
|f(\xi,\Xi)|\le C\ |\xi|\ |\Xi|\quad\quad\mbox{ and }\quad\quad|g(\xi)|\le C\ |\xi|^3\quad.
\ee
Let $u\in W^{1,2}(B^m,{\R}^N)$ satisfying
\be
\label{VI.4}
\Delta u=f(u,\nabla u)+g(u)\quad,
\ee
then $u$ is $C^\infty$.\hfill $\Box$
\end{Th}
\noindent{\bf Proof of theorem~\ref{th-VI.2}.}
First we observe that in dimension 4 only the theorem is not a straightforward consequence of a classical {\it bootstrap argument} in elliptic PDE. Indeed, if $m<4$, using the embedding 
\[
W^{1,2}(B^m)\hookrightarrow L^q(B^m)
\]
for $q=6$ in 3 dimension and any $q<+\infty$ in 2 dimensions, we obtain from the equation (\ref{VI.4}) that $\Delta u\in L^{3/2}$ in 3 dimension and $\Delta u\in L^p$ for any $p<2$. Using classical Calderon-Zygmund theory
this last information gives that $\nabla u\in L^3_{loc}(B^3)$ in 3 dimension and respectively $\nabla u\in L^q(B^2)$ for any $q<+\infty$ in 2 D. So we have gained informations about the regularity of $u$. Injecting again
this improved regularity for $u$ in the equation we obtain that $u$ is bounded in $L^\infty$ and, after one more bootstrap that $u$ is lipshitz. The proof of the full regularity is established then by induction taking more and more derivatives
of the PDE and arguing as we just did with the successive derivatives of $u$.

\medskip

We shall now concentrate on the case $m=4$ which is critical and for which the direct bootstrap argument we just used in the sub-critical dimensions 2 and 3 is not offering any more information on the regularity of $u$.
In 4 dimension we have the embedding
\[
W^{1,2}(B^4)\hookrightarrow L^4(B^4)
\]
We claim that there exists $\al>0$ such that
\be
\label{VI.5}
\sup_{x_0\in B^4_{1/2}(0)\ ;\ 0<\rho<1/4}\rho^{-4\,\al}\int_{B_\rho(x_0)}|u|^4(x)\ dx^4<+\infty\quad.
\ee
Let $\epsilon>0$ to be fixed later.
There exists $\rho_0>0$ such that
\be
\label{VI.6}
\sup_{x_0\in B^4_{1/2}(0)\ ;\ 0<\rho<\rho_0}\int_{B^4_\rho(x_0)}[|u|^4(x)+|\nabla u|^2(x)] dx<\ep\quad.
\ee
Let now $x_0\in B^4_{1/2}(0)$ and $\rho<\rho_0$ arbitrary. On $B_\rho(x_0)$ we consider $\varphi$ to be the solution of
\be
\label{VI.7}
\lf\{
\begin{array}{l}
\ds\Delta \varphi=f(u,\nabla u)+g(u)\quad\quad\mbox{ in }B^4_\rho(x_0)\quad\\[5mm]
\ds\varphi=0\quad\quad\mbox{ on }\quad\p B^4_\rho(x_0)\quad
\end{array}
\rg.
\ee
Classical elliptic estimates (see \cite{GT}) gives the existence of a constant independent of $\rho$ such that
\be
\label{VI.8}
\begin{array}{l}
\ds\|\varphi\|_{L^4(B^4_\rho(x_0))}\le C\ \|f(u,\nabla u)+g(u)\|_{L^{4/3}(B^4_\rho(x_0))}\\[5mm]
\ds\quad\quad\quad\quad\le C\ \|u\|_{L^4(B^4_\rho(x_0))}\ \|\nabla u\|_{L^2(B^4_\rho(x_0))}+C\ \|u\|^3_{L^4(B^4_\rho(x_0))}\quad.
\end{array}
\ee
The difference $v:=u-\varphi$ is harmonic on $B^4_\rho(x_0)$. Hence $|v|^4$ is subharmonic 
\[
\Delta|v|^4=12\, |v|^2\ |\nabla v|^2\ge 0
\]
This gives that
\[
\forall\ r<\rho\quad\quad \int_{\p B_r(x_0)}\frac{\p |v|^4}{\p r}\ge 0
\]
which implies that
\[
\frac{d}{dr}\lf[\frac{1}{r^4}\int_{B^4_r(x_0)}|v|^4(x) \ dx^4\rg]\ge 0
\]
So we have in particular
\be
\label{VI.9}
\int_{B_{\rho/4}(x_0)}|v|^4(x)\ dx^4\le4^{-4} \int_{B_{\rho}(x_0)}|v|^4(x)\ dx^4
\ee
From this inequality we deduce
\be
\label{VI.10}
\begin{array}{l}
\ds\int_{B_{\rho/4}(x_0)}|u|^4(x)\ dx^4\le 8\ \int_{B_{\rho/4}(x_0)}[|v|^4+|\varphi|^4]\ dx^4\\[5mm]
\ds\le 2^{-5} \int_{B_{\rho}(x_0)}|v|^4(x)\ dx^4+ 8\ \int_{B_{\rho}(x_0)}|\varphi|^4\ dx^4\\[5mm]
\ds\le 2^{-2} \int_{B_{\rho}(x_0)}|u|^4(x)\ dx^4+16\ \int_{B_{\rho}(x_0)}|\varphi|^4\ dx^4
\end{array}
\ee
Combining (\ref{VI.8}) and (\ref{VI.10}) we then have
\be
\label{VI.11}
\begin{array}{l}
\ds\int_{B_{\rho/4}(x_0)}|u|^4(x)\ dx^4\\[5mm]
\ds\quad\le \lf[2^{-2}+C_0\ [\|\nabla u\|_{L^2(B_\rho)}^4 + \|u\|_{L^4(B_\rho)}^8]\rg]\int_{B_{\rho}(x_0)}|u|^4(x)\ dx^4
\end{array}
\ee
We choose $\ep>0$ such that $ C_0\ep^2\le 2^{-1}$ and we have then established that for any $\rho<\rho_0$
\be
\label{VI.12}
\int_{B_{\rho/4}(x_0)}|u|^4(x)\ dx^4\le \frac{1}{2}\int_{B_{\rho}(x_0)}|u|^4(x)\ dx^4\quad
\ee
Iterating this inequality gives (\ref{VI.5}). Inserting the {\it Morrey bound} (\ref{VI.5}) into the equation (\ref{VI.4}) gives
\be
\label{VI.13}
\sup_{x_0\in B^4_{1/2}(0)\ ;\ 0<\rho<1/4}\rho^{-4\,\al/3}\int_{B_\rho(x_0)}|\Delta u|^{4/3}(x)\ dx^4<+\infty\quad.
\ee
The {\it Adams-Sobolev} embeddings (see \cite{Ad}) give then the existence of $p>2$ such that $\nabla u\in L^p_{loc}(B^4_{1/2}(0))$. It is easy then to see that the PDE (\ref{VI.4}) becomes
sub-critical for $W^{1,p}$ (with $p>2$) in 4 dimension and we can apply a similar bootstrap arguments as we did in dimensions 2 and 3 in order to obtain the desired regularity for $u$. This concludes
the proof of theorem~\ref{th-VI.2}.\hfill $\Box$.

\medskip

\begin{Rm}
\label{rm-VI.22}
The proof of the regularity of Yang-Mills fields in the critical 4 dimension is ''soft'' in comparison with the proof of the regularity of the ''cousin problem'' : the harmonic maps between a surface and a manifold.
Both equations are critical respectively in 4 and 2 dimensions but the analysis of the harmonic map equation is made more delicate by the fact that the non-linearity in the harmonic map equation is in $L^1$
which is a space which does not behave ''nicely'' with respect to Calderon-Zygmund operations. There is no such a difficulty for Yang-Mills. What is delicate however is to construct a ''good gauge'' in which Yang-Mills equation becomes elliptic.
In a somewhat parallel way the difficulty posed by the harmonic maps  equation was overcome by the author by solving a gauge problem (see \cite{Ri2}). 
\end{Rm}

\medskip
One consequence of the previous regularity result and the point removability result~\ref{th-V.4} is the following {\it point removability} theorem for Yang-Mills fields in 4 dimension
\begin{Th}
\label{th-V.6}{\bf[Point Removability for Yang-Mills in conformal dimension.]}
Let $A$ be a weak solution in $W^{1,2}_{loc}(\wedge^1B^4,{\mathcal G})$ to Yang-Mils equation
\[
d^\ast_AF_A= d^\ast F_A+[A,\res F_{A}]=0\quad\quad\mbox{in }{\mathcal D}'(B^4\setminus \{0\})\quad.
\]
Assume
\[
\int_{B^4}|dA+A\wedge A|^2\ dx^4<+\infty
\]
then there exists  $g\in W^{2,2}_{loc}(B^4\setminus\{0\})$ such that 
\[
A^g\in C^\infty(B^4)
\]
and $A^g$ solves the Yang-Mills equation strongly in the whole ball $B^4$.\hfill $\Box$
\end{Th}
\noindent{\bf Proof of theorem~\ref{th-V.6}.}
The point removability result~\ref{th-V.2}  qives the existence of a $W^{2,2}_{loc}-$gauge such that $A^g\in W^{1,2}(\wedge^1 B^4,{\mathcal G})$. Using
Uhlenbeck's Coulomb gauge~\ref{th-III.2} extraction theorem we can assume that $A^g$ satisfies the Coulomb condition $d^\ast A^g=0$ in ${\mathcal D}'(B^4)$. So $A^g$ is a $W^{1,2}-$solution of a system of the form 
\[
\Delta u=f(u,\nabla u)+g(u)\quad\mbox{in }{\mathcal D}'(B^4\setminus\{0\})
\]
where $f$ and $g$ are smooth maps satisfying (\ref{VI.3aa}). The distribution $\Delta u-f(u,\nabla u)-g(u)$ is supported in $\{0\}$. Hence by a classical
result in distribution theory this distribution is a \underbar{finite} linear combination of derivatives of Dirac masses :
\be
\label{VI.13a}
-\Delta u+f(u,\nabla u)+g(u)=\sum_{|\al\le N} C_\al\ \p_{\al}\delta_0
\ee
where $N\in {\N}$, $\al=(\al_1\cdots\al_4)\in {\N}^4$, $|\al|:=|\al_1|+\cdots+|\al_4|$,  $C_\al\in {\R}^4$ and $\p_\al$ denotes the partial derivative
 $$
 \p_\al:=\frac{\p^{\al_1}}{\p{x_1^{\al_1}}}\ \frac{\p^{\al_2}}{\p{x_2^{\al_2}}}\ \frac{\p^{\al_3}}{\p{x_3^{\al_3}}}\ \frac{\p^{\al_4}}{\p{x_4^{\al_4}}}\
$$
Let $\chi$ be an arbitrary smooth compactly supported function in $B^4_1(0)$. Denote $\chi_\ep(x):=\chi(x/\ep)$, multiply equation (\ref{VI.13a}) by this
function and integrate over $B^4$ gives
\[
\sum_{|\al|\le N}C_\al \frac{\p^\al\chi(0)}{\ep^{|\al|}}=\frac{1}{\ep}\int_{B^4_\ep}\nabla\chi(x/\ep)\cdot\nabla u+\int_{B^4_\ep}\chi_\ep\ \lf[f(u,\nabla u)+g(u)\rg]
\]
Hence, using $\nabla u\in L^2$ and $u\in L^4$ we have
\[
\lf| \sum_{|\al|\le N}C_\al \frac{\p^\al\chi(0)}{\ep^{|\al|}}\rg|=o(\ep)
\]
Since $\p^\al\chi(0)$ are arbitrary, this implies that $C_\al=0$ for any $\al$. So the equation $-\Delta u+f(u,\nabla u)+g(u)=0$ holds on the
whole ball and we can apply theorem~\ref{th-VI.2} to $u=A^g$ and obtain that it is $C^\infty$ which concludes the proof of theorem~\ref{th-V.6}.
\hfill $\Box$

\medskip

It is clear that the solutions to the Yang-Mills Plateau problems satisfy the Yang-Mills equation and hence we have the following corollary.

\begin{Co}
\label{co-V.2}
Let $G$ be a compact Lie group and $m\le 4$. For any 1-form $\eta\in H^{1/2}(\wedge^1\p B^m,{\mathcal G})$  the following minimization problem 
\be
\label{V.7-aa}
\inf\lf\{\int_{B^m}|F_A|^2\ dx^m\ ;\ \iota_{\p B^m}^\ast A=\eta^g\ \mbox{ for some }g\in H^{3/2}(\p B^4,G)    \rg\}
\ee
is achieved by a 1-form $A^0\in W^{1,2}(\wedge^1B^m,{\mathcal G})$ which is $C^\infty$ in any local $W^{1,2}-$Coulomb gauge inside the ball $B^m$.\hfill $\Box$
\end{Co}

\section{Concentration compactness and energy quantization for Yang-Mills Fields in critical dimension.}
\reset

The goal of this section is to study establish the behavior of sequences of Yang-Mills fields of uniformly bounded energy in critical dimension 4. There are three main problematics attached
to this study

\medskip

$\bullet$ Modulo extraction of subsequence, do we have strong converge to a limiting Yang-Mills ?

\medskip 

$\bullet$ If the strong convergence does not hold where is located the lack of strong convergence in the base ?

\medskip

$\bullet$ How much Yang-Mills energy is lost at the limit ?

\medskip

We have already several tools and results at hand that we established in the previous sections in order to provide a relatively
precise answer to these three questions. The proof or our main result in this section is based in particular on the following ''quantitative reformulation''
of the regularity theorem~\ref{th-VI.1} which belongs to the family of the so called $\epsilon-$regularity results for {\it conformally invariant problems}.
\begin{Th}
\label{th-X.1-a} 
{\bf[$\epsilon-$regularity for Sobolev solutions to Yang-Mills in conformal dimension]}
Let $G$ be a compact Lie group, there exists $\ep_{G,4}>0$ such that for any  ${\mathcal G}-$valued 1-forms $A$ in $W^{1,2}(\wedge^1B^4_1(0),{\mathcal G})$   satisfying the Yang-Mills equation
\[
d_{A}^\ast F_{A}= d^\ast F_A+[A,\res F_{A}]=0\quad\mbox{ in }\quad{\mathcal D}'(B^4_1(0)).
\]
and the small energy condition
\[
\int_{B^4_1(0)}|F_A|^2\ dx^4<\ep_{G,4}
\]
then there exists a gauge $g$ in which the following estimates holds : for any $l\in {\N}$ there exists $C_l>0$ such that
\be
\label{X.1-0}
\|\nabla^l(A)^g\|^2_{L^\infty(B_{1/2}(0))}\le C_l\ \int_{B^m_1(0)}|F_A|^2\ dx^4\quad.
\ee
\hfill $\Box$
\end{Th}
\noindent{\bf Proof of theorem~\ref{th-X.1-a}.}
We choose $\ep_{G,4}>0$ that will be definitively fixed a bit later in the proof to be at least smaller than  the $\epsilon_G>0$ of the {\it Coulomb gauge} extraction result theorem~{th-III.2}. We now work in this {\it Coulomb gauge} and we omit to mention the superscript $g$. So, from now on until the end of the proof, we are then assuming that
we have
\be
\label{X.1-1}
\int_{B^4_1}|A|^2\ dx^4+\sum_{i,j=1}^4|\p_{x_i}A_j|^2\, dx^4\le C_G\, \int_{B^4_1}|F_A|^2\, dx^4<\ep_{G,4}
\ee
from which we deduce in particular
\be
\label{X.1-2}
\int_{B^4_1}|A|^4\ dx^4\le C_0\ \lf[\int_{B^4_1}|F_A|^2\, dx^4\rg]^2
\ee
for some constant $C_0>0$. Recall that in the Coulomb gauge we are choosing, the connection form $A$ satisfies the elliptic system (\ref{VI.3aaa-0}) to which we can apply the arguments of
the proof of theorem~\ref{th-VI.2} that we are going to follow closely keeping track this time of each estimate. In particular, having chosen $\ep_{G,4}$ small enough we have inequality (\ref{VI.12}) which holds for $u=A$ and for any 
$B_\rho(x_0)\subset B_1(0)$ and we deduce 
\be
\label{X.1-3}
\begin{array}{l}
\ds\forall\, x_0\ \in B_{3/4}(0)\quad\forall \, \rho<1/4\\[5mm]
\ds\int_{B_\rho(x_0)}|A|^4\ dx^4\le 2\ \rho^\al\ \int_{B^4_1(0)}|A|^4\ dx^4\\[5mm]
\ds\quad\quad\quad\le 2\  C_0\ \rho^\al\  \lf[\int_{B^4_1(0)}|F_A|^2\, dx^4\rg]^2
\end{array}
\ee
where $\al=\log 2/\log 4$. Inserting this inequality in the Yang-Mills  PDE in Coulomb gauge (\ref{VI.3aaa-0}) we obtain the existence of a constant $C_1>0$
such that
\be
\label{X.1-4}
\begin{array}{l}
\ds\forall\, x_0\ \in B_{3/4}(0)\quad\forall \, \rho<1/4\\[5mm]
\ds\int_{B_\rho(x_0)}|\Delta A|^{4/3}\ dx^4\le C_1\ \rho^{\al/3}\ \lf[\int_{B^4_1(0)}|F_A|^2\, dx^4\rg]^{4/3}
\end{array}
\ee
Combining (\ref{X.1-3}) and (\ref{X.1-4}) we deduce from {\it Adams-Morrey inequalities} (see \cite{Ad})
\be
\label{X.1-5}
\begin{array}{l}
\ds\|\nabla A\|_{L^p(B_{3/5}(0))}\le C\ \sup_{x_0\in B_{3/4}(0)\,; \, \rho<1/4}\lf[\rho^{-\al/3}\int_{B_\rho(x_0)}|\Delta A|^{4/3}\ dx^4  \rg]^{3/4}\\[5mm]
\ds\quad\quad+C\ \|A\|_{L^2(B_1(0))}
\end{array}
\ee
where
\[
p=\frac{16-4\,\al/3}{8-\al}>2\quad.
\]
Hence we have for this $p>2$
\[
\|\nabla A\|_{L^p(B_{3/5}(0))}\le C\ \lf[\int_{B^4_1(0)}|F_A|^2\, dx^4\rg]^{1/2}\quad.
\]
Since $p>2$ the non-linear elliptic system (\ref{VI.3aaa-0}) becomes sub-critical  in 4 dimensions and a standard bootstrap argument gives (\ref{X.1-0}). This concludes the
proof of theorem~\ref{th-X.1-a}.\hfill $\Box$
The previous $\epsilon-$regularity result is the main step for proving the following {\it concentration compactness} theorem for sequences of Yang-Mills
fields.
\begin{Th}
\label{th-X.1-b}{\bf[Concentration compactness for Yang-Mills Fields in conformal dimension]}
Let $(M^4,h)$ be a closed 4 dimensional riemannian manifold and $P$ a principal smooth $G$ bundle over $M^4$. Let $\nabla^k$ be a sequence of Yang-Mills
connections satisfying
\[
\limsup_{k\rightarrow +\infty}\int_{M^4}|F_{\nabla^k}|^2_h\ dvol_{h}<+\infty
\]
Then there exists a subsequence $\nabla^{k'}$, a smooth $G-$bundle $P^\infty$ over $(M^4,h)$, a smooth Yang-Mills connection $\nabla^\infty$ of $P^\infty$
and finitely many points $\{p_1\cdots p_N\}$ in $M^4$ such that for any contractible open set $D^4\subset M^4\setminus \{p_1\cdots p_N\}$ there exists
a sequence of trivialization of $P$ over $D^4$ for which
\[
A^{k'}\rightarrow A^\infty\quad\quad\mbox{ strongly in }C^l(D^4)\quad\forall\, l\in {\N}
\]
where $A^{k'}$ (resp. $A^\infty$) is the connection 1-form associated to $\nabla^{k'}$ (resp. $\nabla^\infty$) in this sequence of local  trivializations of $P$  (resp. $P^\infty$) over $D^4$. Moreover we have the following weak convergence in Radon measure
\be
\label{X.1-7-a}
\mu^{k'}:=|F_{\nabla^{k'}}|^2_h\ dvol_h\rightharpoonup \mu^\infty:=|F_{\nabla^{\infty}}|^2_h\ dvol_h+\nu
\ee
where $\nu$ is a non negative atomic measure supported by the points $p_j$
\be
\label{X.1-6}
\nu:=\sum_{j=1}^N f_j\ \delta_{p_j}\quad.
\ee
\hfill $\Box$
\end{Th}
{\bf Proof of theorem~\ref{th-X.1-b}}. We follow step by step the proof of theorem~\ref{th-V.3} replacing for the choice of the covering the Uhlenbeck Coulomb Gauge threshold $\ep_G(M^4,h)$
by the smaller positive constant $\ep_{G,4}(M^4,h)$ given by the $\epsilon-$regularity result~\ref{th-X.1-a} on the manifold $(M^4,h)$. Observe that, because
of the epsilon regularity result, the Coulomb gauges $(A^k)^{g^k_i}$ are pre-compact for any $C^l-$topology on each ball $B_\rho(x_i)$. This gives also the pre-compactness of the transition functions $g^k_{ij}$ in any of the $C^l$ topologies. Hence, the  co-cycles $g_{ij}^k$ converge in any of these topology to the limiting (now smooth) co-cycle $g^\infty_{ij}$ which defines a smooth $G-$bundle $P^\infty$ over $M^4\setminus \{p_1\cdots p_N\}$. Moreover the limiting collection of 1-forms $A^{j,\infty}$ defines a connection $\nabla^\infty$ on $P^\infty$ satisfying also the Yang-Mills equation which passes obviously to the limit under $C^\infty$ convergence. The gauge invariant quantities such as $|F_{\nabla^k}|^2$ converge also to the corresponding limiting quantities and we have then, modulo extraction of a further  subsequence, the existence of a limiting radon measure $\nu$ supported on the points $p_j$ exclusively such that
(\ref{X.1-7}) holds. Finally applying the point removability theorem and once again the $\epsilon-$regularity we extend $\nabla^\infty$ globally on $M^4$ as a Yang-Mills 
smooth connection of the bundle $P^\infty$ which also obviously extends throughout the $p_j$ as a smooth bundle. This concludes the proof of theorem~\ref{th-X.1-b}. \hfill $\Box$

\medskip

Finally, we identify the {\it concentration atomic measure} by proving that the weights $f_j$ in front of the Dirac masses $\delta_{p_j}$ are the sums of Yang-Mills energies of Yang-Mills fields over $S^4$, the so called ''bubbles''. Precisely we have the following {\it energy identity result} which was first established for instantons in \cite{Tau} and for Yang-Mills fields in general in \cite{Ri3}. The proof we present below is using the {\it interpolation Lorentz spaces} following a technic introduced in \cite{LR1} and \cite{LR2}. 

\begin{Th}
\label{th-X.1-c}{\bf [Energy quantization for Yang-Mills Fields in conformal dimension.]}
Let $(M^4,h)$ be a closed 4 dimensional riemannian manifold and $P$ a principal smooth $G$ bundle over $M^4$. Let $\nabla^k$ be a sequence of Yang-Mills connections of uniformly bounded Yang-Mills energy converging strongly away from finitely many points $\{p_1\cdots\p_N\}$ to a limiting Yang-Mills connection $\nabla^\infty$ as described in theorem~\ref{th-X.1-b}. Let $\nu$ be the atomic
concentration measure 
\[
\nu:=\sum_{j=1}^N f_j\ \delta_{p_j}
\]
satisfying
\[
\mu^{k'}:=|F_{\nabla^{k'}}|^2_h\ dvol_h\rightharpoonup \mu^\infty:=|F_{\nabla^{\infty}}|^2_h\ dvol_h+\nu\quad.
\]
Then for each $j=1\cdots N$ there exists finitely many $G-$Yang-Mills connections $(D_j^i)_{i=1\cdots N_j}$ over $S^4$ such that
\be
\label{X.1-7}
\forall\, j=1\cdots N\quad f_j=\sum_{i=1}^{N_j}\int_{S^4}|F_{D_j^i}|^2\ dvol_{S^4}\quad.
\ee
\hfill $\Box$
\end{Th}
\noindent{\bf Proof of theorem~\ref{th-X.1-c}.} Since the result is local, the metric in the domain does not play much role and we will present the proof for $M^4=B^4_1(0)$ equipped with the flat metric and assuming moreover that there is exactly one limiting blow-up point, $N=1$, which coincide with the origin, $p_1=0$. We also express the connection $\nabla^k$ in s 

\medskip

Recall that we denote by $\ep_{G,4}$ the positive constant given by the $\epsilon-$regularity theorem~\ref{th-X.1-a}. We detect the ''most concentrated bubble''
precisely let
\[
\rho^k:=\inf\lf\{\rho\ ;\ \exists x\in B_1^4(0)\quad\mbox{ s. t. } \int_{B_\rho(x)}|F_{A^k}|^2\ dx^4=\frac{\ep_{G,4}}{2}\rg\}\quad.
\]
Since we are assuming that blow-up is happening exactly at the origin we have that 
\[
\lf\{
\begin{array}{l}
\ds\rho^k\rightarrow 0\quad\mbox{ and }\\[5mm]
\ds\exists\, x^k\rightarrow 0\quad\mbox{ s. t. } \int_{B_{\rho^k}(x^k)}|F_{A^k}|^2\ dx^4=\frac{\ep_{G,4}}{2}
\end{array}
\rg.
\]
We choose a sequence $x^k$ that we call {\it center of the first bubble} and $\rho^k$ is called the {\it critical radius of the first bubble}. Let
\[
\hat{A}^k(y):=\rho^k\ \sum_{i=1}^4 A^k_i(\rho^k\, y+x^k)\ dy_i\quad.
\]
Due to the scaling invariance of the Yang-Mills lagrangian in 4-dimensions, $\hat{A}^k$, which is the pull-back of $A^k$ by the dilation map
$D^k(y):=\rho^k\, y+x^k$, is a Yang-Mills fields moreover
\be
\label{X.1-8}
\max_{y\in B_{1/(2\rho^k)}(0)}\int_{B_1^4(y)}|F_{\hat{A}^k}|^2\ dy^4=\int_{B_1^4(0)}|F_{\hat{A}^k}|^2\ dy^4=\frac{\ep_{G,4}}{2}\quad.
\ee
Applying the $\epsilon-$regularity theorem~\ref{th-X.1-a} we deduce that on  {\it Uhlenbeck's Coulomb Gauges} $g$ which
exists on each unit ball $B_1(y)$ since the $\epsilon_{G,4}$ has been taken smaller than $\epsilon_G$ from theorem~\ref{th-III.2}
\[
\forall\, l\in {\N}\quad\sup_{y\in B_{1/(2\rho^k)}(0)}\|\nabla^l(\hat{A})^g\|_{L^\infty(B_{1/2})}\le C_l
\]
where $C_l$ is independent of $k$. Hence, locally in Coulomb gauge, modulo extraction of a subsequence, the sequence $\hat{A}^k$ converges strongly in any $C^l_{loc}$ topology on ${\R}^4$ to a limiting Yang-Mills connection $\hat{A}^\infty$ satisfying
\[
\int_{B_1^4(0)}|F_{\hat{A}^k}|^2\ dy^4=\frac{\ep_{G,4}}{2}
\]
and which is therefore non trivial. Let now $\pi$ be the stereographic projection with respect to the north pole, due to the {\it conformal invariance} of Yang-Mills energy $\tilde{A}:=\pi^\ast \hat{A}^\infty$ is a non trivial Yang-Mills Field on $S^4\setminus\{\mbox{south p\^ole}\}$. Since $\tilde{A}$ is a smooth ${\mathcal G}-$valued 1-form with finite Yang-Mills energy and satisfying the Yang-Mills equation we can apply the {\it point removability} result for Yang-Mills fields, theorem~\ref{th-V.6} and conclude that $\tilde{A}$ extends to a global smooth Yang-Mills $G-$connection $D_1^1$   over the whole $S^4$  which is our {\it first bubble} and using again the conformal invariance of Yang-Mills energy we have
\be
\label{X.1-9}
\frac{\ep_{G,4}}{2}\le\int_{S^4}|F_{D^1_1}|^2\ dvol_{S^4}=\lim_{R\rightarrow +\infty}\lim_{k\rightarrow +\infty} \int_{B_{R\,\rho^k}(x^k)}|F_{A^k}|^2\ dx^4\quad.
\ee
We have now to study the loss of Yang Mills energy in the so called {\it neck region} between the {\it first bubble} and the {\it macroscopic solution} $A^\infty$ to which $A^k$ converges
away from zero. Precisely we are studying
\be
\label{X.1-10}
\begin{array}{l}
\ds\lim_{k\rightarrow +\infty}\int_{B^4}|F_{A^k}|^2 dx^4-\int_{B^4}|F_{A^\infty}|^2 dx^4-\int_{S^4}|F_{D^1_1}|^2\ dvol_{S^4}\\[5mm]
\ds\quad\quad\quad=\lim_{R\rightarrow +\infty}\lim_{k\rightarrow +\infty} \int_{B_{R^{-1}}(x^k)\setminus B_{R\,\rho^k}(x^k)}|F_{A^k}|^2\ dx^4
\end{array}
\ee
For any $G$ there is a minimal Yang-Mills energy among all non-trivial Yang-Mills Fields. This can be proved easily observing that if the energy is less than 
the $\ep_{G,4}$ threshold, the connections can be represented by a global smooth  Yang-Mills 1-form on $S^4$ to which the $C^l$ estimates of theorem~\ref{th-X.1-a} apply. Hence since $A$ satisfies globally on $S^4$ the PDE (\ref{VI.3aaa-0}) and for small enough Yang-Mills energy this implies that $A$ is an
harmonic one form on $S^4$ which gives that it is a trivial Yang-Mills fields. Denote
\[
\mbox{YM}(G,S^4)=\min\lf\{\int_{S^4}|F_D|^2\ dvol_{S^4}\ ; \ D\mbox{ is a non zero }G\mbox{ Yang-Mills Field}\rg\}
\]
To simplify the presentation we assume that 
\be
\label{X.1-11}
\begin{array}{l}
\ds\lim_{k\rightarrow +\infty}\int_{B^4}|F_{A^k}|^2 dx^4-\int_{B^4}|F_{A^\infty}|^2 dx^4-\int_{S^4}|F_{D^1_1}|^2\ dvol_{S^4}\\[5mm]
\quad\quad\quad\quad\quad\quad\quad\quad<\mbox{YM}(G,S^4)\quad,
\end{array}
\ee
or in other words 
\be
\label{X.1-12}
\lim_{R\rightarrow +\infty}\lim_{k\rightarrow +\infty} \int_{B_{R^{-1}}(x^k)\setminus B_{R\,\rho^k}(x^k)}|F_{A^k}|^2\ dx^4<\mbox{YM}(G,S^4)\quad.
\ee
Alternatively we would have to go through some standard and fastidious induction procedure to remove all the bubbles one by one - each of them taking at least an amount of YM$(G,S^4)$ Yang-Mills energy - and we would be anyway reduced at the end to study
the loss of energy in annuli region where (\ref{X.1-12}) holds (such a procedure is described for instance in \cite{BR} proposition III.1 in the framework of
{\it Willmore surfaces}). 
Under the assumption (\ref{X.1-12}) the goal is ultimately to prove
\[
\lim_{R\rightarrow +\infty}\lim_{k\rightarrow +\infty} \int_{B_{R^{-1}}(x^k)\setminus B_{R\,\rho^k}(x^k)}|F_{A^k}|^2\ dx^4=0\quad,
\]
that will finish the proof of the theorem. We are now going to prove the following claim.

\medskip

\noindent{\bf Claim 1 :}
\[
\begin{array}{l}B_{R_\delta\,\rho^k}(x^k)
\ds\forall \, \delta >0\quad\exists\, R_\delta>1\quad \mbox{s. t. }\quad\forall\ r\in[R_\delta\,\rho^k,R_\delta^{-1}]\\[5mm]
\ds \limsup_{k\rightarrow +\infty}\int_{B_{2\,r}(x^k)\setminus B_{r}(x^k)}|F_{A^k}|^2\ dx^4<\delta
\end{array}
\]
\noindent{\bf Proof of claim 1.} We argue by contradiction. Assume there exists $\delta_0>0$ such that for all $R>1$ there exists $r^k\in[R\,\rho^k,R^{-1}]$
\[
\limsup_{k\rightarrow +\infty}\int_{B_{2\,r^k}(x^k)\setminus B_{r^k}(x^k)}|F_{A^k}|^2\ dx^4>\delta_0\quad.
\]
Since we can find a sequence $r_k$ for any $R>1$, using a diagonal argument and the extraction of a subsequence we can assume that
\[
\frac{r^k}{\rho^k}\rightarrow +\infty\quad\mbox{ and }\quad r^k\rightarrow 0\quad.
\]
Now we introduce
\[
s^k:=\inf\lf\{
\begin{array}{c}
\ds s\ ;\ \exists\, x\in B_{2\,r^k}(x^k)\setminus B_{r^k}(x^k)\quad\mbox{ s. t. }\\[5mm]
 \ds\int_{B_s(x)}|F_{A^k}|^2\ dx^4=\min\lf\{\frac{\ep_{G,4}}{2},\frac{\delta_0}{16}\rg\}
 \end{array}
 \rg\}
\]
Let $\ti{x}^k$ be a point in the dyadic annuls $B_{2\,r^k}(x^k)\setminus B_{r^k}(x^k)$ where this infimum is achieved. We clearly have 
$$B_{s^k}(\ti{x}^k)\subset B_{R^{-1}}(x^k)\setminus B_{R\,\rho^k}(x^k)$$
for any $R>1$ and $k$ large enough. Dilating the Yang-Mills connection $A^k$ about $\ti{x}^k$ at a rate $(s^k)^{-1}$ we again obtain
a limiting non trivial Yang-Mills field, a second {\it bubble}, $\ti{A}^\infty$ either on ${\R}^4$ or on ${\R}^4\setminus\{0\}$ depending whether $s^k/r^k$ tends to zero or not.
In any case we have
\[
\lim_{R\rightarrow +\infty}\lim_{k\rightarrow +\infty} \int_{B_{R^{-1}}(x^k)\setminus B_{R\,\rho^k}(x^k)}|F_{A^k}|^2\ dx^4\ge\int_{{\R}^4}|F_{\ti{A}^\infty}|^2\ dx^4\ge \min\lf\{\frac{\ep_{G,4}}{2},\frac{\delta_0}{16}\rg\}
\]
and the {\it point removability} theorem~\ref{th-V.6} for Yang-Mills fields implies that $\ti{A}^\infty$ extends to a non-trivial Yang-Mills Field on $S^4$
\[
\int_{{\R}^4}|F_{\ti{A}^\infty}|^2\ dx^4\ge\mbox{YM}(G,S^4)
\]
which contradicts the fact that we are working under the assumption that there is only one {\it bubble} (i.e. assumption  (\ref{X.1-12})). So we have proved claim 1.

\medskip

Combining claim 1 and the $\epsilon-$regularity  theorem~\ref{th-X.1-a} we obtain 
\[
\begin{array}{l}
\ds\forall \, \delta >0\quad\exists\, R_\delta>1\quad \mbox{s. t. }\quad\\[5mm]
\ds \forall\ x\in B_{R_\delta^{-1}}\setminus B_{R_\delta\,\rho^k}(x^k) \quad\quad |x|^2\ |F_{A^k}|^2(x)<\delta
\end{array}
\]
Consider an Uhlenbeck Coulomb gauge $(A^k)^{g^k}$ in the annulus $B_{2\,R_\delta\,\rho^k}(x^k)\setminus B_{R_\delta\,\rho^k}(x^k)$. Introduce
$\chi$ to be a cut-off function such that
\[
\lf\{
\begin{array}{l}
\ds\chi(x)\equiv 1\quad\quad\mbox{ in }{R}^4\setminus B^4_2(0)\\[5mm]
\ds\chi(x)\equiv 0\quad\quad\mbox{ in } B^4_1(0)
\end{array}
\rg.
\]
and let $\chi^k(x):=\chi(R_\delta\,\rho^k\ (x-x^k))$. Extend $(A^k)^{g^k}$ in $B_{2\,R_\delta\,\rho^k}(x^k)$ by taking
\[
\ti{A}^k:=\chi^k (A^k)^{g^k}
\]
Using again the $C^l$ estimates (\ref{X.1-0})  of $\epsilon-$regularity  theorem~\ref{th-X.1-a} for the Uhlenbeck Coulomb gauge $(A^k)^g$ in the annulus $B_{2\,R_\delta\,\rho^k}(x^k)\setminus B_{R_\delta\,\rho^k}(x^k)$ we obtain for any $x\in B_{2\,R_\delta\,\rho^k}(x^k)\setminus B_{R_\delta\,\rho^k}(x^k)$ 
\be
\label{X.1-13}
\begin{array}{l}
\ds|x|\ |\ti{A}^k|(x)+|x|^{2}\ |\nabla \ti{A}^k|(x)\\[5mm]
\ds\quad\le\, C\ \lf[\int_{B_{4\,R_\delta\,\rho^k}(x^k)\setminus B_{R_\delta\,\rho^k/2}(x^k)}|F_{A^k}|^2\ dx^4\rg]^{1/2}
\end{array}
\ee
We have then produced an extension $\ti{A}^k$ of $A^k$ inside
the ball $B_{2\,R_\delta\,\rho^k}(x^k)$ equal to $A^k$ in $B_{R_\delta^{-1}}\setminus B_{R_\delta\,\rho^k}(x^k)$ and satisfying
\[
\||x|^2\ |F_{\ti{A}^k}|\|_{L^\infty(B_{R_\delta^{-1}}(x^k))}<\sqrt{\delta}
\]
This implies in particular
\be
\label{X.1-13-aa}
\|F_{\ti{A}^k}\|_{L^{2,\infty}(B_{R_\delta^{-1}}(x^k))}\le C\ \sqrt{\delta}
\ee
where $C>0$ is a constant independent of $\delta$ and $k$. Taking $\delta$ small enough we can apply theorem~\ref{th-III.3} and find a gauge that we denote
simply $\ov{A}^k$ and which satisfies
\be
\label{X.1-13-a}
\lf\{
\begin{array}{l}
\ds\int_{B_{R_\delta^{-1}}(x^k)}|\ov{A}^k|^2+\sum_{i,j=1}^4|\p_{x_i}\ov{A}^k_j|^2\, dx^m\le C_G\, \int_{B_{R_\delta^{-1}}(x^k)}|F_{\ov{A}^k}|^2\, dx^4\\[8mm]
\ds d^\ast \ov{A}^k=0\quad\quad\mbox{ in }B_{R_\delta^{-1}}(x^k)\\[5mm]
\ds\iota_{\p B_{R_\delta^{-1}}(x^k) }^\ast(\ast \ov{A}^k)=0
\end{array}
\rg.
\ee
Using the gauge invariance of Yang-Mills integrant (\ref{VI.3-sst}) together with  the $C^l$ estimates (\ref{X.1-0})  of $\epsilon-$regularity  theorem~\ref{th-X.1-a} applied to the Uhlenbeck Coulomb gauge $(A^k)^g$ in the annulus $B_{2\,R_\delta\,\rho^k}(x^k)\setminus B_{R_\delta\,\rho^k}(x^k)$ we obtain
on $B_{R_\delta^{-1}}(x^k)$
\be
\label{X.1-14}
\begin{array}{l}
|d^\ast_{\ov{A}^k}F_{\ov{A}^k}|\le C\ |\nabla \chi^k|\ |(A^k)^{g^k}|^2+C\ |\chi^k|\ |(A^k)^{g^k}|^3\\[5mm]
\quad\le C\ (R_\delta\,\rho^k)^{-3}\ {\mathbf 1}_{k,\delta}\ \delta
\end{array}
\ee
where ${\mathbf 1}_{k,\delta}$ is the characteristic function of $B_{2\,R_\delta\,\rho^k}(x^k)\setminus B_{R_\delta\,\rho^k}(x^k)$. This implies the following estimate
\be
\label{X.1-15}
\|d^\ast_{\ov{A}^k}F_{\ov{A}^k}\|_{L^{(4/3,1)}(B_{R_\delta^{-1}}(x^k))}\le\ C\ \int_{ B_{2\,R_\delta\,\rho^k}(x^k)\setminus B_{R_\delta\,\rho^k}(x^k) }|F_{\ov{A}^k}|^2\, dx^4\quad.
\ee
where we recall that $L^{(4/3,1)}$ is the Lorentz space whose dual is the {\it Marcinkiewicz weak} $L^4$ space : $L^{4,\infty}$. Using the embeding
(\ref{III.19}) for $p=2$ and $m=4$ 
\[
W^{1,2}(B^4)\hookrightarrow L^{4,2}(B^4)\quad,
\]
we obtain from (\ref{X.1-13-a}) the estimate
\be
\label{X.1-16}
\|\ov{A}^k\|^2_{L^{4,2}(B_{R_\delta^{-1}}(x^k))}\le C_G\, \int_{B_{R_\delta^{-1}}(x^k)}|F_{\ov{A}^k}|^2\, dx^4\quad.
\ee
Using now one of  the embeddings (\ref{III.17}) :
\[
L^{2}(B^4)\cdot L^{4,2}(B^4)\hookrightarrow L^{4/3,1}(B^4)\quad,
\]
we obtain
\be
\label{X.1-17}
\|d^\ast d\ov{A}^k\|_{L^{4/3,1}(B_{R_\delta^{-1}}(x^k))}\le C_G\, \int_{B_{R_\delta^{-1}}(x^k)}|F_{\ov{A}^k}|^2\, dx^4\quad.
\ee
Combining this fact with the three lines of (\ref{X.1-13-a}) together with classical elliptic estimates in Lorentz spaces (see \cite{SW}) gives
\[
\|F_{\ov{A}^k}\|^2_{L^{2,1}(B_{R_\delta^{-1}}(x^k))}\le C_G\, \int_{B_{R_\delta^{-1}}(x^k)}|F_{\ov{A}^k}|^2\, dx^4\quad.
\]
Combining this inequality with the estimate (\ref{X.1-13-aa})  of the curvature in the dual space $L^{2,\infty}$ in the {\it neck region} we obtain
\[
\begin{array}{l}
\ds\forall\ \delta>0\quad\exists\, R_\delta>1\quad\mbox{ s. t. }\\[5mm]
\ds\limsup_{k\rightarrow+\infty}\int_{B_{R_\delta^{-1}}(x^k)\setminus  B_{2\,R_\delta\,\rho^k}(x^k) }|F_{{A}^k}|^2\, dx^4\le C\ \sqrt{\delta}
\end{array}
\]
from which we deduce
\[
\lim_{R\rightarrow +\infty}\lim_{k\rightarrow +\infty} \int_{B_{R^{-1}}(x^k)\setminus B_{R\,\rho^k}(x^k)}|F_{A^k}|^2\ dx^4=0
\]
This implies
\[
|F_{A^{k}}|^2\ dx^4\rightharpoonup \mu^\infty:=|F_{A^{\infty}}|^2\ dx^4+\int_{S^4}|F_{D^1_1}|^2\ dvol_{S^4}\ \delta_0\quad.
\]
This completes the proof of the theorem~\ref{th-X.1-c}.   \hfill $\Box$

\section{The resolution of the Yang-Mills Plateau problem in super-critical dimensions.}
\reset

\subsection{The absence of $W^{1,2}$ local gauges.}

We can reformulate the sequential weak closure of $W^{1,2}$ connections we proved in the previous sections for the dimensions up to 4 in the following way. Let $G$ be a compact Lie group and 
$(M^m,h)$  a compact riemanian manifold. Introduce the space of so called {\it Sobolev connections} defined by
\[
{\mathfrak A}_G(M^m):=\lf\{ 
\begin{array}{l}
A\in L^2(\wedge^1M^m,{\mathcal G})\ ; \ \int_{M^m}|dA+A\wedge A|_h^2\ dvol_h<+\infty\\[5mm]
\mbox{ locally }\exists\ g\in W^{1,2}\quad\mbox{ s.t. }\quad A^g\in W^{1,2}
\end{array} 
\rg\}
\]
then we have proved the following result.

\begin{Th}
\label{th-VII.1}
For $m\le 4$ the space ${\mathfrak A}_G(M^m)$ is weakly sequentially closed below any given Yang-Mills energy level : precisely
For any $A^k\in  {\mathfrak A}_G(M^m)$ satisfying
\[
\limsup_{k\rightarrow+\infty}YM(A^k)=\int_{M^m}|dA^k+A^k\wedge A^k|_h^2\ dvol_h<+\infty
\]
there exists a subsequence $A^{k'}$ and a Sobolev connection $A^\infty\in  {\mathfrak A}_G(M^m)$ such that
\[
d(A^{k'},A^\infty):=\inf_{g\in W^{1,2}(M^m,G)}\int_{M^m}|A^{k'}-(A^\infty)^g|_h^2\ dvol_h\longrightarrow 0
\]
moreover
\[
YM(A^\infty)\le\liminf_{k'\rightarrow 0} YM(A^{k'})\quad.
\]
\hfill $\Box$
\end{Th}
\begin{Rm}
\label{rm-VII.1}
Observe that the space ${\mathfrak A}_G(M^m)$ contains for instance global $L^2$ one forms taking values into the Lie algebra ${\mathcal G}$ that correspond to smooth connections of some
principal $G-$bundle over $M^m$. If the Yang-Mills energy of  a sequence of such smooth connections is uniformly bounded, we can extract a subsequence converging weakly to a Sobolev connection and  corresponding 
possibly to \underbar{another} $G-$bundle. This possibility of ''jumping'' from one bundle to another, as predicted for instance in the {\it concentration compactness} result theorem~\ref{th-X.1-b}, is encoded in the definition of ${\mathfrak A}_G(M^m)$. \hfill $\Box$
\end{Rm}
Because of this weak closure property the space ${\mathfrak A}_G(M^m)$ is the ad-hoc space for minimizing Yang-Mills energy in dimension less or equal than 4. This is however not the case in higher dimension. We have the following proposition.
\begin{Prop}
\label{pr-VII.1}
For $m> 4$ the space ${\mathfrak A}_{SU(2)}(M^m)$ is \underbar{ not} weakly sequentially closed below any given Yang-Mills energy level : precisely
there exists $A^k\in  {\mathfrak A}_{SU(2)}(M^m)$ satisfying
\[
\limsup_{k\rightarrow+\infty}YM(A^k)=\int_{M^m}|dA^k+A^k\wedge A^k|^2\ dx^m<+\infty
\]
 and a Sobolev connection $A^\infty\in  L^2$ such that
\[
d(A^{k'},A^\infty):=\inf_{g\in W^{1,2}(M^m,SU(2))}\int_{M^m}|A^{k'}-(A^\infty)^g|_h^2\ dvol_h\longrightarrow 0
\]
but in every neighborhood $U$ of every point of $M^m$ there is \underbar{no} $g$ such that $(A^\infty)^g\in W^{1,2}(U)$ .\hfill $\Box$
\end{Prop}
The proposition is not difficult to prove but we prefer to illustrate this fact by a small cartoon. We consider a sequence $A^k$ of smooth 1-forms on $B^5$, the unit 5 dimensional ball, into $su(2)$ such that $\limsup_{k\rightarrow+\infty}YM(A^k)<+\infty$. The drawing is representing the flow lines of the divergence free vector field associated to the closed {\it Chern 4-form} : $Tr(F_{A^k}\wedge F_{A^k})$.
In this cartoon $A^k$ is weakly converging in $L^2$ to some limit $su(2)-$valued 1-form $A^\infty$ such that $F_{A^\infty}\in L^2$ but this 1-form satisfies
\[
d\lf(Tr(F_{A^\infty}\wedge F_{A^\infty})\rg)=8\pi^2\ [\delta_P-\delta_N]\ne 0
\]
where $P$ and $N$ are two distinct points of $B^5$ - the red dots in the cartoon. Hence for almost every small radii $r>0$ we have for instance
\[
\int_{\p B^5_r(P)} Tr(F_{A^\infty}\wedge F_{A^\infty})=8\,\pi^2
\]
Assume there would exist $g\in W^{1,2}$ such that $(A^\infty)^g\in W^{1,2}$ in the neighborhood of $P$. The gauge invariance of the Chern form
gives that
\[
\int_{\p B^5_r(P)} Tr(F_{(A^\infty)^g}\wedge F_{(A^\infty)^g})=8\pi^2
\]
However we have seen in section III that the fact that $\iota^\ast_{\p B^5_r}(A^\infty)^g$ is in $W^{1,2}(\wedge^1 \p B^5_r,su(2))$ for almost every $r$ - due to Fubini theorem - imposes
\[
\int_{\p B^5_r(P)} Tr(F_{(A^\infty)^g}\wedge F_{(A^\infty)^g})=0
\]
which is a contradiction.
\begin{center}
\includegraphics[width=5cm,height=4cm]{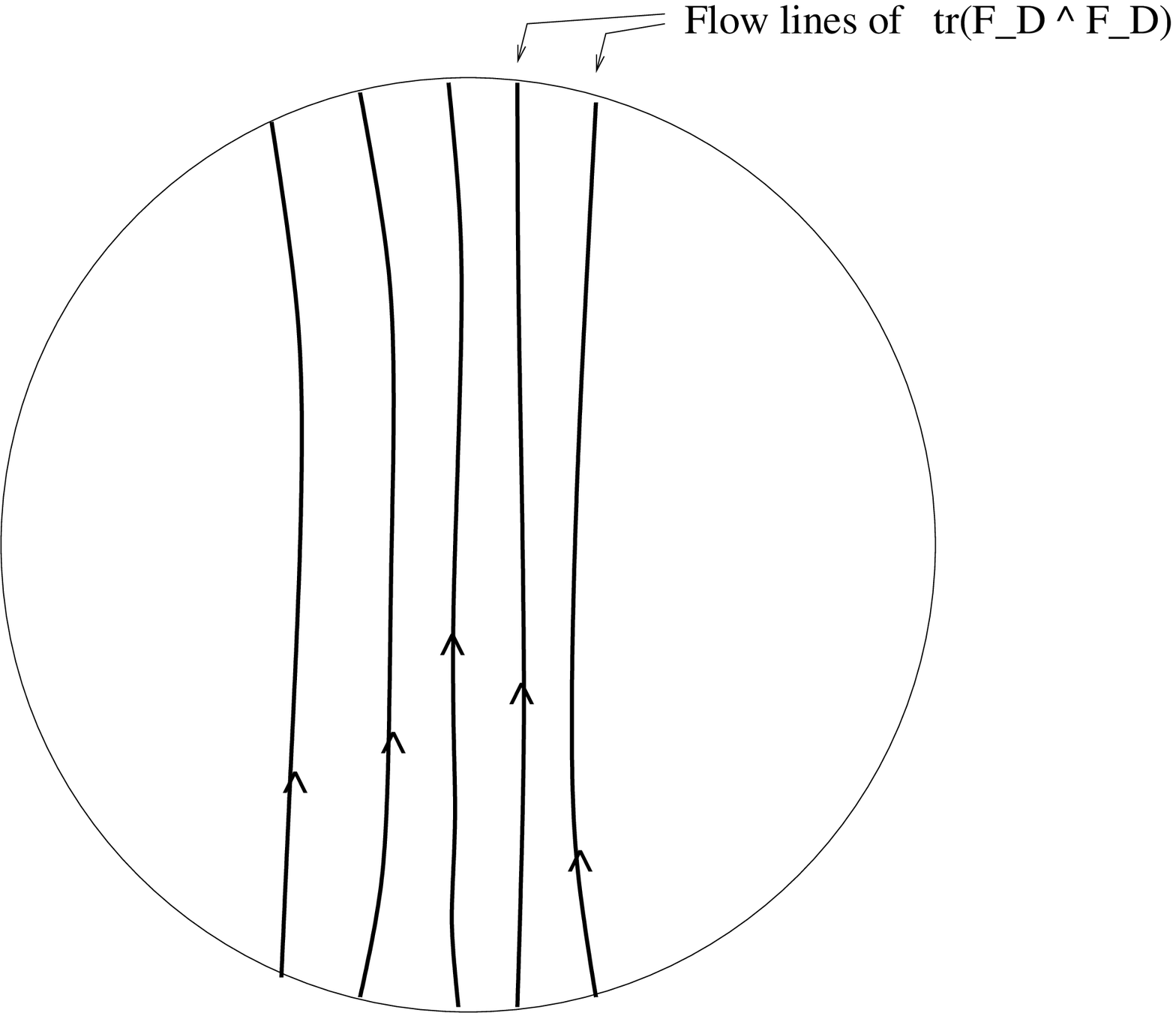}
\end{center}

\begin{center}
Fig.~1: Sequence of smooth connections time 1
\end{center}

\begin{center}
\includegraphics[width=5cm,height=4cm]{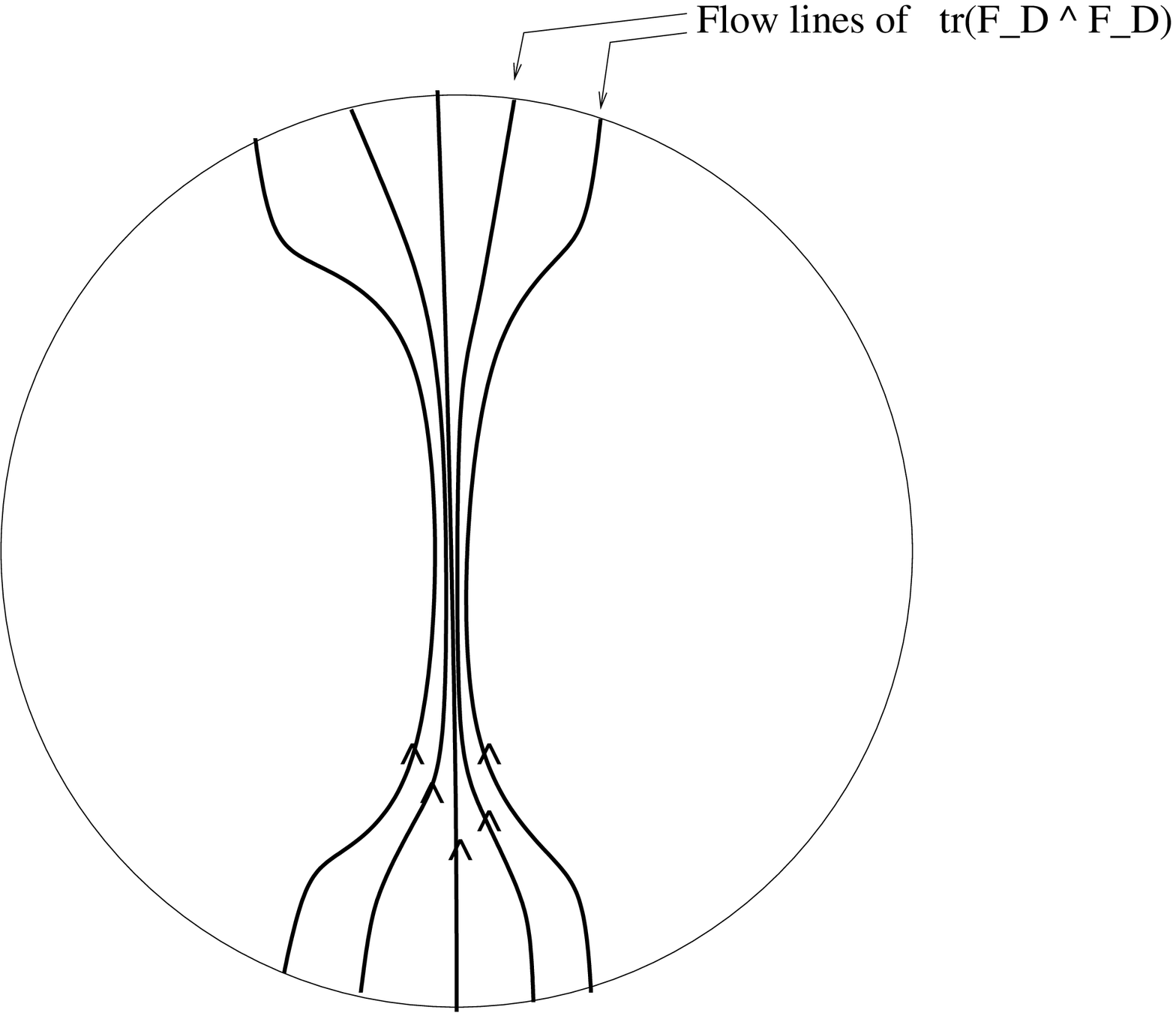}
\end{center}

\begin{center}
Fig.~2: Sequence of smooth connections time 2

\end{center}

\begin{center}
\includegraphics[width=5cm,height=4cm]{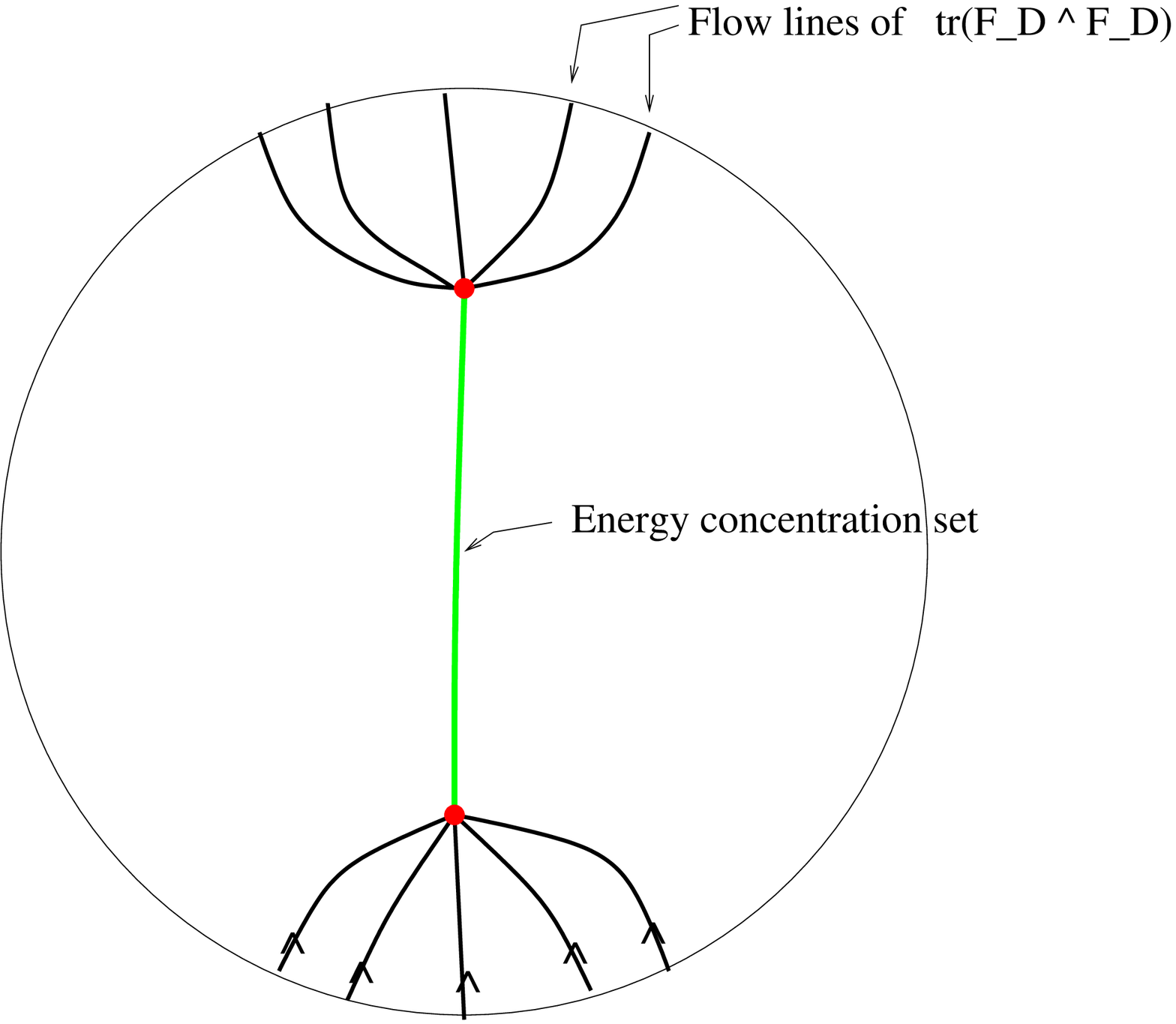}
\end{center}

\begin{center}
Fig.~3: Sequence of smooth connections - the limit.

\end{center}
The construction of the counter example of proposition $A^\infty$ for $M^m=B^5$ can be achieved by generating countably many dipoles of the form ($P$, $N$) at the limit which realize a dense subset of $B^5$ and then in such a way that
\[
supp\lf[d\lf(Tr(F_{A^\infty}\wedge F_{A^\infty})\rg)\rg]=B^5\quad.
\]
This can be done by controlling the Yang-Mills energy.

In the above cartoon the limiting 1-form $A^\infty$ is a 1-form of a smooth connection but on a $SU(2)-$bundle which is only defined over $B^5\setminus \{P\}\cup\{N\}$. What this example says is that, starting in 5 dimension, the Yang-Mills
energy is not coercive enough in such a way that it's control does not prevent the corresponding bundle to degenerate and to have local twists at the limit. In order to find an ad-hoc weakly sequentially closed space below any yang-Mills energy level, in such a way that the {\it Yang-Mills Plateau problem} can be solved we have then to relax the notion of {\it Sobolev connections},
which was implicitly assuming that the underlying bundle was smooth, by allowing the bundle, the carrier of the connection, to have singularities. This effort is similar to the one produced by Federer and Fleming while producing the class of {\it integer rectifiable currents} - i.e. sub manifolds with singularities in a way - in order to solve the {\it Classical Plateau Problem} in super critical dimension $m>2$. We are looking for a {\it Geometric measure theoretic} version of bundle and connections.

\subsection{Tian's results on the  compactification of the space of smooth Yang-Mills Fields in high dimensions.}

The need of developing a  {\it Geometric measure theoretic} version of bundle and connections beyond the too small class of {\it Sobolev connections} on smooth bundles has been already encountered in the study of the compactification of the
{\it moduli space } of smooth yang-Mills fields by G.Tian in \cite{Ti}.

\begin{Th}
\label{th-VII.2}
Let $G$ be a compact Lie group and $A^k$ be a sequence of ${\mathcal G}-$valued 1-forms in $B^m$. Assume $A^k$ are all smooth solutions to the Yang-Mills equation
\[
d_{A^k}^\ast F_{A^k}=0\quad.
\]
Assume
\[
\limsup_{k\rightarrow +\infty}\int_{B^m}|dA^k+A^k\wedge A^k|^2\ dx^m<+\infty
\]
Then there exists a subsequence $A^{k'}$ and a limiting  ${\mathcal G}-$valued 1-forms $A^\infty$
\[
d(A^{k'},A^\infty):=\inf_{g\in W^{1,2}(M^m,G)}\int_{M^m}|A^{k'}-(A^\infty)^g|_h^2\ dvol_h\longrightarrow 0
\]
Moreover there exists a $m-4$ rectifiable closed subset of $B^m$, $K$,  of finite $m-4$  Hausdorff measure, ${\mathcal H}^{m-4}(K)<+\infty$ such that
in 
\[
\begin{array}{l}
\forall\ B_r(x_0)\subset B^m\setminus K\quad \quad\exists\ g\in W^{1,2}(B_r(x_0),G)\quad\mbox{ s. t. }\\[5mm]
(A^\infty)^g\mbox{ is a smooth solution of Yang-Mills equation in }B_r(x_0) 
\end{array}
\]
and the following weak convergence as Radon measure holds
\be
\label{VII.1}
|F_{A^{k'}}|^2\ dx^m\rightharpoonup |F_{A^\infty}|^2\ dx^m+ f\ {\mathcal H}^{m-4}\res K
\ee
where ${\mathcal H}^{m-4}\res K$ is the restriction to $K$ of the $m-4$ Hausdorff measure and $f$ is an absolutely continuous function
with respect to this measure.
\hfill $\Box$
\end{Th}
This result is very close to a similar result proved in \cite{Li} by F.H.Lin for harmonic maps in super-critical dimension $m\ge 3$.

\medskip

\noindent{\bf Proof of theorem~\ref{th-VII.2}.}
The starting point of the proof of theorem~\ref{th-VII.2} is the following monotonicity formula computed first by P. Price in \cite{Pri}.
\begin{Prop}
\label{pr-VII.2}{\bf[Monotonicity formula]}
Let $m\ge 4$ and $A$ be a ${\mathcal G}-$valued 1-forms in $B^m_1(0)$ assume that $A$ is a smooth solution to the Yang-Mills equation
\[
d_{A}^\ast F_{A}=0\quad\quad\mbox{ in }B^m_1(0).
\]
then the following monotonicity formula holds
\be
\label{VII.2-aaa}
\forall\, p\in B^m_1(0)\quad\forall\ B^m_r(p)\subset B^m_1(0)\quad\frac{d}{dr}\lf[\frac{1}{r^{m-4}}\int_{B^m_r(p)}|F_A|^2\ dx^m\rg]\ge 0
\ee
\hfill $\Box$
\end{Prop}
The {\it monotonicity formula} is a direct consequence of the {\it stationarity condition} which is satisfied by any smooth critical point of the Yang-Mills
lagrangian. Precisely the {\it stationarity condition} says
\be
\label{VII.1-aaa}
\forall\ X\in C^\infty_0(B^m,{\R}^m)\quad\quad\lf.\frac{d}{dt}\int_{B^m}|\Psi_t^\ast F_A|^2\rg|_{t=0}=0
\ee
where $\Psi_t$ is the flow of $X$. The monotonicity formula is obtained by applying (\ref{VII.1-aaa}) to the each vector-field of the following form :  On the ball $B_r(p)$ the vector-field $X$  is equal to the radial one, $X=\p/\p r$ for canonical coordinates centered at $p$ which generates the dilations centered at $p$, and it realizes a  smooth interpolation to $0$ outside $B_{r+\delta}(p)$ for any $\delta>0$ (see \cite{Pri}). Once such a vector field is chosen one computes (\ref{VII.1-aaa}) and make $\delta$ tend to zero. This computation gives then (\ref{VII.2-aaa}).

\medskip

The second ingredient of the proof is the extension of the {\it Coulomb gauge extraction} in dimension larger than $4$ to the framework of the so called Morrey spaces
where the $m-4$ densities of Yang-Mills energy are assumed to be small everywhere and at any scale. The following result has been obtained independently by T.Tao and G.Tian in \cite{TT} and by Y.Meyer together with the author of the present notes in \cite{MR} .
\begin{Th}
\label{th-VII-2-b}{\bf [Coulomb Gauge extraction]}
Let $m\ge 4$ and $G$ be a compact Lie group, there exists $\ep_{m,G}>0$ such that for any \underbar{smooth} ${\mathcal G}-$valued 1-forms $A$ in $B^m_1(0)$  satisfying the small Morrey energy condition
\[
\|F_A\|^2_{M^0_{2,4}(B^m_1(0))}:=\sup_{ p\in B^m_{1}(0),\ r>0}\frac{1}{r^{m-4}}\int_{B^m_r(p)\cap B^m_1(0)}|F_A|^2\ dx^m<\ep_{m,G}
\]
then there exists a gauge  $g\in W^{2,2}(B_{1/2}^m(0),G)$ such that
\[
\sup_{ p\in B^m_{1}(0),\ r>0}\frac{1}{r^{m-4}}\int_{B^m_r(p)\cap B^m_1(0)}\sum_{i,j=1}^m|\p_{x_i}(A^g)_j|^2\le C\ \|F_A\|^2_{M^0_{2,4}(B^m_1(0))}
\]
and
\[
\lf\{
\begin{array}{l}
\ds d^\ast(A)^g=0\quad\quad\mbox{ in }B^m_1(0)\\[5mm]
\ds \iota_{\p B^m_1(0)}^\ast (A)^g=0
\end{array}
\rg.
\]
where the constant $C$ only depends on $m$ and $G$.\hfill $\Box$
\end{Th}
In this Coulomb gauge any Yang-Mills smooth connection one form   satisfies
\be
\label{VI.3aabd}
\Delta A^g=d^\ast(A^g\wedge A^g)+[A^g,\res dA^g]+[A^g,\res (A^g\wedge A^g)]\quad,
\ee
We shall make now use of the following generalization of theorem~\ref{th-VI.2} to Morrey spaces.
\begin{Th}
\label{th-VI.2-b}
Let $m\ge 4$ and $N\in {\N}^\ast$. Let $f\in C^\infty({\R}^N\times({\R}^m\otimes{\R}^N),{\R}^N)$ and let $g\in C^\infty({\R}^N,{\R}^N)$ such that there exists $C>0$
satisfying
\be
\label{VI.3aabb}
|f(\xi,\Xi)|\le C\ |\xi|\ |\Xi|\quad\quad\mbox{ and }\quad\quad|g(\xi)|\le C\ |\xi|^3\quad.
\ee
There exists $\ep>0$ such that for any $u$ in $L^4\cap W^{1,2}(B^m,{\R}^N)$ satisfying
\be
\label{VI.3aabc}
\sup_{ p\in B^m_{1}(0),\ r>0}\frac{1}{r^{m-4}}\int_{B^m_r(p)\cap B^m_1(0)}|\nabla u|^2\ dx^m<\ep_{m,G}
\ee
and
\be
\label{VI.4-b}
\Delta u=f(u,\nabla u)+g(u)\quad,
\ee
then we have for any $l\in{\N}$ the existence of $C_l>0$ such that
\be
\label{VI.4-c}
\|\nabla^l u\|^2_{L^\infty(B^m_{1/2}(0))}\le C_l\ \sup_{ p\in B^m_{1}(0),\ r>0}\frac{1}{r^{m-4}}\int_{B^m_r(p)\cap B^m_1(0)}|\nabla u|^2\ dx^m
\ee
.\hfill $\Box$
\end{Th}
The proof of this theorem is more or less identical to the one of  theorem~\ref{th-VI.2} replacing the different spaces  in $4$ dimensions by their {\it Morrey  counterparts} in higher dimension, bearing in mind that Calderon-Zygmund theory extends with the natural exponents to these spaces (see \cite{Mo}).

\medskip

Combining the {\it monotonicity formula},  theorem~\ref{th-VII-2-b} and theorem~\ref{th-VI.2-b} applied to the PDE (\ref{VI.3aabd}), adapting the arguments we followed for proving the corresponding result - theorem~\ref{th-X.1-a} - from the 4-D counterparts of theorem~\ref{th-VI.2-b} in the conformal dimension 4,  we obtain the following $\epsilon-$regularity result~\footnote{An $\epsilon-$regularity theorem for smooth Yang-Mills fields has first been obtained by H.Nakajima (see \cite{Na}). It is however a ''gauge invariant reslult'' which gives only an $L^\infty$ bound on the curvature under the small energy assumption but is not providing any control of the connection in some gauge. The proof of Nakajima $\epsilon$-regularity for smooth Yang-Mills fields 
is following the arguments originally introduced by R. Schoen in \cite{Sch} for proving the corresponding result for smooth harmonic map. It is using the {\it Bochner Formula}
as a starting point together with the {\it maximum principle} and the {\it Moser iteration technique}. }
\begin{Th}
\label{th-VII-2-a}{\bf[$\epsilon-$regularity for smooth Yang-Mills]}
Let $m\ge 4$ and $G$ be a compact Lie group, there exists $\ep_{m,G}>0$ such that for any \underbar{smooth} ${\mathcal G}-$valued 1-forms $A$ in $B^m_1(0)$   satisfying the Yang-Mills equation
\[
d_{A}^\ast F_{A}=0\quad\mbox{ in }B^m_1(0).
\]
and the small energy condition
\[
\int_{B^m_1(0)}|F_A|^2\ dx^m<\ep_{m,G}
\]
then there exists a gauge $g$ in which the following estimates holds : for any $l\in {\N}$ there exists $C_l>0$ such that
\[
\|\nabla^l(A)^g\|^2_{L^\infty(B_{1/2}(0))}\le C_l\ \int_{B^m_1(0)}|F_A|^2\ dx^m
\]
\hfill $\Box$
\end{Th}
\noindent{\bf Proof of theorem~\ref{th-VII.2} continued.}

\medskip

Let
\[
E^k_{r}:=\lf\{p\in B^m\ ;\ \frac{1}{r^{m-4}}\int_{B_r(p)}|F_{A^{k}}|\ dx^m\ge\ep_{m,G}\rg\}
\]
where $\ep_{m,G}$ is the epsilon in the $\epsilon-$regularity theorem~\ref{th-VII-2-a}.The {\it monotonicity formula} implies
\[
\forall k'\in {\N}\quad \forall \ r<\rho\ \quad E^k_{r}\subset E^k_{\rho}\quad .
\]
Hence, by a standard diagonal argument we can extract a subsequence such that $E_{k',2^{-j}}$ converges to a limiting closed set $E_{\infty,2^{-j}}$ which of course satisfy
\[
E^{\infty}_{2^{-j-1}}\subset E^{\infty}_{ 2^{-j}}\quad.
\]
Let
\[
K:=\bigcap_{j\in{\N}}E^{\infty}_{ 2^{-j}}
\]
A classical Ferderer-Zimmer covering argument gives
\[
{\mathcal H}^{m-4}(K)<+\infty\quad.
\]
With the $\epsilon-$regularity theorem at hand, extracting possibly a further subsequence following a diagonal argument, we can ensure that  $A^{k'}$ converges locally away from $K$ in every $C^l-$norm  in the  Coulomb gauges constructed in theorem~\ref{th-VII-2-b} and we have
\[
\mu^{k'}:=|F_{A^{k'}}|^2\ dx^m\rightharpoonup \mu^\infty:=|F_{A^\infty}|^2\ dx^m+\nu
\]
where $\nu$ is a Radon measure supported in the closed set $K$. Because of the Radon measure convergence, the monotonicity formula (\ref{VII.2-aaa}) satisfied by $A^k$ can be transferred to the measure $\mu^\infty$
\[
\forall\, p\in B^m_1(0)\quad\forall\ B^m_r(p)\subset B^m_1(0)\quad\frac{d}{dr}\lf[\frac{\mu^\infty(B^m_r(p))}{r^{m-4}}\ dx^m\rg]\ge 0
\]
from which we deduce
\[
\theta^{m-4}(\mu^\infty,p):=\lim_{r\rightarrow 0}\frac{\mu^\infty(B^m_r(p))}{r^{m-4}}\ge 0\quad\mbox{ exists for every }p\in B^m_1(0)\quad.
\]
Observe that
\be
\label{VII.xx-1}
K=\lf\{p\in B^m_1(0)\ ;\ \theta^{m-4}(\mu^\infty,p)>0\rg\}\quad.
\ee
Using the monotonicity we have that for any $\rho\in (0,1)$ and any $p\in K\cap B_\rho(0)$ 
\[
\nu(B^m_r(p))\le \frac{r^{m-4}}{(1-\rho)^{m-4}}\lim_{k'\rightarrow +\infty}\int_{B^m_1(0)}|F_{A^{k'}}|^2\ dx^m
\]
We deduce from this inequality that $\nu$ is {\it absolutely continuous} with respect to ${\mathcal H}^{m-4}\res K$, the restriction to $K$ of the $m-4-$Hausdorff measure. Let $\delta>0$, define
\[
G_\delta:=\lf\{p\in B^m_1(0)\ ;\ \delta<\limsup_{r\rightarrow 0}\frac{1}{r^{m-4}}\int_{B_r^m(p)}|F_{A^\infty}|^2\ dx^m\rg\}\quad.
\]
Considering for any $\delta>0$ for any $\eta>0$ and any $p\in G_\delta$ a radius $0<r_p^\eta<\eta$ such that $$\int_{B_{r_p^\eta}^m(p)}|F_{A^\infty}|^2\ dx^m\ge \frac{\delta}{2}\ (r_p^\eta)^{m-4}$$
For any $\eta>0$ we extract from the covering $(B_{r_p^\eta}^m(p))_{p\in G_\delta}$ a Besicovitch subcovering $(B_{r_{i}}^m(p))_{i\in I}$ of $G_\delta$
in such a way that there exists an integer $N_m>0$ depending only on $m$ such that each point  of $B^m$ is covered by at most $N_m$ balls of this sub-covering. We then have
\[
\begin{array}{l}
\ds{\mathcal H}^{m-4}(G_\delta)\le\frac{2}{\delta}\limsup_{\eta\rightarrow 0} \sum_{i\in I}\int_{B_{r_{i}}^m(p)}|F_{A^\infty}|^2\ dx^m\\[5mm]
\ds\quad\quad\quad\quad\le \frac{2\, N_m}{\delta}\limsup_{\eta\rightarrow 0}\int_{dist(x,K)<\eta}|F_{A^\infty}|^2\ dx^m
\end{array}
\]
Since $K$ is closed and has Lebesgue measure zero we deduce that
\[
{\mathcal H}^{m-4}\lf(\bigcap_{\delta>0} G_\delta\rg)=0
\]
or in other words
\[
\mbox{ for }{\mathcal H}^{m-4}\mbox{ almost every }p\in B^m\quad \limsup_{r\rightarrow 0}\frac{1}{r^{m-4}}\int_{B_r^m(p)}|F_{A^\infty}|^2\ dx^m=0
\]
Using the characterization of $K$ given by (\ref{VII.xx-1}) and the fact that $\nu$ is {\it absolutely continuous} with respect to ${\mathcal H}^{m-4}\res K$ we deduce from the previous fact that
\[
\mbox{ for }\nu\mbox{ almost every }p\in B^m\quad\lim_{r\rightarrow 0}\frac{\nu(B_r(p))}{r^{m-4}}\quad\mbox{ exists and is positive}
\]
The following result, which is an important contribution to {\it Geometric measure theory} was proved by D.Preiss in \cite{Pre} and answered
positively to a conjecture posed by Besicovitch. It permits to conclude the proof of theorem~\ref{th-VII.2}.\hfill $\Box$
\begin{Th}
\label{th-VII-preiss}
Let $\nu$ be a Borel non-negative measure in $B^m$. Assume that $\nu$ almost everywhere the $n-$dimensional density of $\nu$ exists and is positive
then $\nu$ is supported by a $n-$dimensional rectifiable subset in $B^m$.\hfill $\Box$
\end{Th}

\medskip

In order to have a more complete description of the {\it compactification of the moduli space} of smooth Yang-Mills fields in supercritical dimension two main questions are left open in this theorem 

\begin{itemize}
\item[i)]
What are the special geometric properties satisfied by the set $K$ ? 
\item[ii)] What is the energy defect $f(x)$ ?
\item[ii)]
What is the regularity of $A^\infty$ - modulo gauge  - throughout $K$ ?
\end{itemize}

In subcritical dimension $m<4$, due to the analysis we have exposed in the previous sections, we have that $K=\emptyset$, $f=0$ and modulo gauge transformation the limiting connection extends to a smooth Yang-Mills field over the whole ball.

\medskip

In critical $4$ dimensions these questions are answered by the point removability theorem~\ref{th-V.6} : the set $K$ is made of isolated points,
the function $f$ is the sum of the Yang-Mills energies of Yang-Mills fields over $S^4$ (see \cite{Ri3}) and, modulo gauge transformation the limiting connection extends to a smooth Yang-Mills field over the whole ball.

\medskip

In super-critical dimension it is expected that $K$ with the multiplicity $f$ defines a a so called {\it stationary varifold} (see \cite{Si}). This belief comes from the fact that smooth Yang-Mills fields satisfy the {\it stationarity condition}  (\ref{VII.1-aaa}) and it is expected that this condition should still be satisfied by the weak limit $A^\infty$ itself, and hence, due to the Radon measure convergence (\ref{VII.1}) , it would be ''transfered'' to the measure $f\ {\mathcal H}^{m-4}\res K$. This  last fact is equivalent to the stationarity of $(K,f)$.

\medskip

Regarding the regularity of $A^\infty$ a result of T.Tao and G.Tian \cite{TT} asserts that $A^\infty$ is a smooth Yang-Mills connection of a smooth  bundle defined away of a closed subset $L\subset K$ satisfying ${\mathcal H}^{m-4}(L)=0$. This partial regularity result is probably not optimal 
but this optimality or non-optimality is an open problem (a similar open question exists for {\it stationary harmonic maps} - see \cite{RS})

\subsection{The $\Omega-$anti-self-dual instantons.}

In 4 dimension a special  class of solutions to the Yang-Mills equation, the {\it anti-self-dual instantons} have been considered by S.Donaldson in the early 80's to produce invariants of differential 4 dimensional manifolds. On a $4$-dmensionnal riemannian manifold $(M^4,h)$ for some given $SU(n)-$bundle over $M^m$ we consider connections $A$ solutions to the equation
\be
\label{VII.2}
\ast_h F_{A}=-F_A\quad.
\ee
This equation is issued from an elliptic complex (see \cite{DK}) and is the natural generalization in 4 dimension of the {\it flat connection equation} $F_A=0$
on riemann surfaces considered to classify holomorphic complex structures over such a surface. It defines special solutions to Yang-Mills equation. Indeed taking the covariant exterior derivative with respect to $A$, using the Bianchi identity $d_AF_A=0$ one obtains $d^{\ast_h}F_A=0$. The {\it anti-self-dual equation} (\ref{VII.2}) is generalized in higher dimension 
\be
\label{VII.3}
\ast_h F_A=-\ \Omega\wedge F_A\quad,
\ee
where $\Omega$ is a \underbar{closed} $m-4$ form. Again, due to {\it Bianchi identity}, by taking the covariant exterior derivative with respect to the connection $A$, using {\it Leibnitz identity} on 1-derivations and the fact that $\Omega$ is closed, we obtain the Yang-Mills equation $d^{\ast_h}F_A=0$. The $\Omega-$anti-self-dual equation is not elliptic in general. There are however special situations of geometric interest when the base manifold has a restricted holonomy.

\medskip

\noindent $\bullet$ {\it Hermitian Yang-Mills fields.}

\medskip

Let $(M^{2n},h)$ be an even dimensional  riemannian manifold. We assume that $M^{2n}$ is equipped with an integrable complex structure  $J_M$ - i.e.
the brackets operation leaves the space of $1-0$ vector fields of $TM\otimes{\C}$ invariant
\[
\forall\ X,Y\mbox{ vector fields }\quad J\,[X-i\,J\, X,Y-i\, J\, Y]= i\, [X-i\,J\, X,Y-i\, J\, Y]
\]
Finally we assume that $(M^{2m},h,J_M)$ is K\"ahler : $\omega(\cdot,\cdot):=g(\cdot,J_M\cdot)$ is a closed 2-form. It defines a non degenerate 2-form and $\om^m/m !=dvol_g$. Let
\[
\Om:=\om^{m-2}/(m-2) !
\]
Consider an hermitian vector bundle $E$ associated to a principal  $SU(n)$ bundle over $M^m$ with projection map $\pi\ : E\rightarrow M^{2m}$. A connection $\nabla$ is $\Omega$ anti self-dual if and only if it satisfies the Hermitian Yang-Mills equations :
\[
\ast F_\nabla=-\ \Omega\wedge F_\nabla\quad\Longleftrightarrow\quad
\lf\{
\begin{array}{l}
\ds F_\nabla^{0,2}=0\\[3mm]
\ds \omega\cdot F_\nabla^{1,1}=0
\end{array}
\rg.
\]
where $F^{0,2}_{\nabla}$  (resp. $F^{1,1}_\nabla$ ) is the $0-2$ (resp. $1-1$) part of the curvature (the space $T^\ast M\otimes{\C}$  is decomposed according to the eigen-spaces of $J_M$ for the eigenvalues $i$ and $-i$) so the Hermitian Yang-Mills equation 
\[
\left\{
\begin{array}{l}
\displaystyle F_\nabla^{0,2}=0\quad\Leftrightarrow\quad\forall\ X,Y\quad F_\nabla(X-iJX,Y-iJY)=0\\[5mm]
\displaystyle\omega\cdot F_\nabla^{1,1}=0\quad\Leftrightarrow\quad\sum_{l=1}^m(F_\nabla)_{\epsilon_l,J\,\epsilon_l}=0
\quad\mbox{ where }\omega=\sum_{l=1}^m\epsilon_l\wedge J\epsilon_l
\end{array}
\right.
\]
and $\ep_l$ denotes an orthonormal basis of $(T^\ast M^{2m},h)$. 

\medskip

The {\it Hermitian Yang-Mills} equation can be interpreted as follows. The {\it equivariant horizontal distribution of plane} $H$ associated to $\nabla$ defines an almost complex structure $J_\nabla$ on $E$ in the following way
\[
\forall \, \xi\in E\ \forall X\in T_\xi E\quad J_\nabla(X):= J_E(X^V)+(J_M(\pi_\ast X))^H
\]
we recall that  $X^V$ is the projection onto the tangent vertical space (the kernel of the projection $\pi_\ast$) with respect to the horizontal plane $H$, and $J_E$ is the complex structure on the tangent vertical space defined by the $SU(n)$ structure group of the bundle.
The first part of the equation {\it Hermitian Yang-Mills} can be reformulated as follows
\[
\begin{array}{l}
\forall\ X,Y\mbox{ vector fields in }M^{2m}\quad F^{0,2}(X,Y)=0\Longleftrightarrow\\[5mm]
 J_\nabla\,[(X-i\,J_M\, X)^H,(Y-i\, J_M\, Y)^H]= i\, [(X-i\,J_M\, X)^H,(Y-i\, J_M\, Y)^H]
\end{array}
\]
which is equivalent to say that $J_\nabla$ is \underbar{integrable} and that the hermitian bundle is holomorphic.

\medskip

 If $A$ is a  $su(n)-$valued 1-form representing $\nabla$ in an orthonormal trivialization, this integrability condition implies, by switching to a local holomorphic trivialization, that there is a gauge change $g$ (non-unitary anymore but taking values into $Gl(n,{\C})$) such that locally
\[
g^{-1}\ov{\p}g+g^{-1}A^{0,1} g=(A^g)^{0,1}=0
\]
and the $0-1$ part of the connection $\nabla$ in this holomorphic trivialization coincide with  $\ov{\p}$. Since $A$ is taking values into $su(n)$, we have that
$A^{1,0}=-\ov{A^{0,1}}^t$. Thus we obtain the fact that 
\[
(A^{1,0})^g=h^{-1}\,\p h\quad,
\]
where $h:=\ov{g}^t\, g$ is taking values into invertible self-dual matrices. So the so called ''{Einstein part}'' of the equation $\om\cdot F^{1,1}_\nabla=0$ becomes
equivalent to the non-linear elliptic equation\footnote{ Recall that in a K\"ahler manifold $\om\cdot\ov{\p}\p=\Delta$}
\[
\om\cdot \ov{\p}[h^{-1}\,\p h]=0\quad.
\]
\medskip

\noindent $\bullet$ {\it $SU(4)-$instantons in Calaby-Yau 4-folds.}

\medskip
In high energy physics and later in geometry (see the PhD thesis of C.Lewis , \cite{DT}, \cite{Ti}) the following generalization of {\it instantons} has been introduced. Consider an $SU(n)$  principal bundle $P$ over $(M^8,g)$ a {\it Calabi-Yau manifold} of complex dimension 4 - i.e. $(M^8,g)$ has holonomy $SU(4)$. Such a manifold 
posses an integrable complex structure $J_M$ for which $(M^8,g,J_M)$ is K\"ahler and it  posses in addition a global holomorphic $4-0$ form $\theta$ of unit norm - unique modulo unit complex number multiplication -  and satisfying
\[
\theta\wedge\ov{\theta}=\frac{\om^4}{4!}
\]
The $4-$form $\Om$ defined by
\[
\Omega:=4\,\Re(\theta)+\frac{\om^2}{2}
\]
is closed, parallel and of unit {\it co-mass} ( see \cite{Ti} lemma 4.4.1) : $\forall\ x\in M$
\[
1=\|\Om\|_\ast(x)=\sup\lf\{\frac{<\Om, u_1\wedge u_2\wedge u_3\wedge u_4 >}{\prod_{i=1}^4 |u_i|}\quad;\quad u_i\in T_xM\setminus \{0\}\rg\}
\]
The holomorphic $4-0$ form $\theta$ defines an isometry of the space of $0-2$ forms on $M^{8}$ such that
\[
\forall\,\al\,,\, \beta\in \wedge^{0,2}M^8\quad\quad \al\wedge\ast_\theta\beta=\al\cdot\beta\ \ov{\theta}
\]
Some basic computation gives the the $\Omega-$anti self dual equation in this case is equivalent to
\[
\ast F_\nabla=-\ \Omega\wedge F_\nabla\quad\Longleftrightarrow\quad
\left\{
\begin{array}{l}
\displaystyle (1+\ast_\theta)\,F_\nabla^{0,2}=0\\[5mm]
\displaystyle\omega\cdot F_\nabla^{1,1}=0
\end{array}
\right.
\]
which is also known as the {\it $SU(4)$ instanton equation}.

\subsection{Tian's regularity conjecture on $\Omega-$anti-self-dual instantons.}

Tian's result in the case of $\Omega$ anti-self dual instantons for a closed $m-4-$form $\Omega$ of co-mass less than 1 is the following.
\begin{Th}
\label{th-VII.3}
Let $(M^m,h)$ be a compact riemannian manifold. Let $\Omega$ be a smooth $m-4$ closed form in $M^m$. Assume $\Omega$ has co-mass less than 1
\[
\||\Om|_\ast\|_{L^\infty(M^m)}\le 1
\]
where
\[
|\Om|_\ast(x):=\sup\lf\{\frac{<\Om, u_1\wedge \cdots\wedge u_{m-4} >}{\prod_{i=1}^{m-4} |u_i|}\quad;\quad u_i\in T_xM\setminus \{0\}\rg\}
\]
Let $E$ be an hermitian vector bundle associated to an $SU(n)$-bundle over $M^m$ and let $\nabla^k$ be a sequence of smooth $SU(n)-$connections satisfying the $\Omega$ anti self dual instantons equation
\[
\ast_h F_{\nabla^k}=-\ \Omega\wedge F_{\nabla^k}\quad\mbox{ in }M^m
\]
Then, modulo extraction of a subsequence,  there exists a $m-4$ rectifiable closed subset of $M^m$, $K$,  of finite $m-4$  Hausdorff measure, ${\mathcal H}^{m-4}(K)<+\infty$, an hermitian bundle $E_0$ defined over $M^m\setminus K$ and a smooth $SU(n)$ connection $\nabla^\infty$ of $E_0$ such that
\be
\label{VII.4}
\ast_h F_{\nabla^\infty}=-\ \Omega\wedge F_{\nabla^\infty}\quad\mbox{ in }M^m\setminus K
\ee
moreover
\[
\begin{array}{l}
\forall\ B_r(x_0)\subset M^m\setminus K\quad \quad\exists\  \mbox{ a sequence of trivializations  s. t. }\\[5mm]
A^{k'}\rightarrow A^\infty\quad \mbox{ strongly in } C^l(B^m)\quad\forall l\in {\N}
\end{array}
\]
where $\nabla^{k'}\simeq d+A^{k'}$ in these trivializations and $\nabla^\infty\simeq d+A^\infty$ in a trivialization of $E_0$ over $B^m$.
Moreover there exists an integer rectifiable current $C$ such that for any smooth $m-4$ form $\varphi$ on $M^m$
\[
\int_{M^m}Tr(F_{\nabla^{k'}}\wedge F_{\nabla^{k'}})\wedge\varphi \rightarrow \int_{M^m}Tr(F_{\nabla^{\infty}}\wedge F_{\nabla^{\infty}})\wedge\varphi +8\pi^2\,C(\varphi)
\]
and the current $C$ is calibrated by $\Om$ 
\[
C(\Om)={\mathbf M}(C)=\sup\{C(\varphi)\ ;\ \||\varphi|_\ast\|_{L^\infty(M^m)}\le 1\}
\]
where ${\mathbf M}$ is the mass of the current $C$. Finally the following convergence holds weakly as Radon measures
\[
|F_{\nabla^{k'}}|_h^2\ dvol_{h}\rightharpoonup|F_{\nabla^{\infty}}|_h^2\ dvol_{h}+8\pi^2\ \Theta(C)\ {\mathcal H}^{m-4}\res K
\]
where $\Theta(C)$ is the integer valued $L^1$ function with respect to the restriction of the $m-4$ Hausdorff measure to $K$ and which is giving the multiplicity of the current $C$ at each point.\hfill $\Box$
\end{Th}
Observe that no bound is a-priori needed for the Yang-Mills energy of the sequence. We have indeed
\[
\begin{array}{l}
\ds YM(\nabla^{k})=-\int_{M^m}Tr(F_{\nabla^{k}}\wedge\ast_hF_{\nabla^{k}})\ dvol_h\\[5mm]
\ds\quad\quad=\int_{M^m}Tr(F_{\nabla^{k}}\wedge F_{\nabla^{k}})\wedge\Om
\end{array}
\]
Since $\Om$ is closed this integral only depends on the cohomology class of the other closed form $Tr(F_{\nabla^{k}}\wedge F_{\nabla^{k}})$ which is the second Chern class
of $E$ and which is independent of $\nabla^k$.

\medskip
\begin{Rm}
\label{rm-VII.3}
Regarding the limiting bundle $E_0$, it is important to insist on the fact that there is no reason for $E_0$ to extend through $K$ over the whole manifold $M^m$
as a smooth bundle. Hence there is a-priori no meaning to give to $\nabla^\infty$ over the whole manifold and therefore the $\Omega-$anti-self dual equation
(\ref{VII.4}) cannot even hold in a distributional sense throughout $K$ if we do not relax the notion of bundle and connections. \hfill $\Box$
\end{Rm}
\medskip

An important question directly related to the regularity of the limiting configuration $(E_0,\nabla^\infty,C)$ is the following open problem.

\medskip

\noindent{\bf Open problem.} {\it Show that the limiting current $C$ has no boundary :
\[
\p C=0\quad.
\]}
\hfill $\Box$

\medskip 

\noindent Observe that this open question is equivalent to 
\[
d\lf(Tr(F_{\nabla^{\infty}}\wedge F_{\nabla^{\infty}})\rg)=0\quad.
\] 
This last question should be equivalent to the following 

\medskip

\noindent{\bf strong approximation property} : Does there exist a sequence $D^k$ of smooth  $SU(n)-$connection of smooth hermitian bundles  over $M^m$ such that
\[
\lim_{k\rightarrow +\infty}\inf_{g\in W^{1,2}} \int_{M^m}|\nabla^\infty-(D^k)^g|_h^2+|F_{\nabla^\infty}-g^{-1}F_{D^k} g|^2_h\ dvol_h= 0\quad ?
\]
\medskip

If we would know that $\p C=0$ then $C$ defines a {\it calibrated integral cycle}. The optimal regularity for {\it calibrated} or {\it semi-calibrated} integral cycles of dimension 2 has been proved in \cite{RT} and \cite{Be}. Such cycles have at most isolated point singularities. More generally, the calibrated condition implies that such a cycle  is homologically {\it mass minimizing} (see \cite{HL}). From this later fact, using  Almgren regularity result \cite{Al} proved also by C.De Lellis and E.Spadaro (\cite{DS1}, \cite{DS2} and \cite{DS3}), we would obtain that $C$ is the integration along a rectifiable set which is a smooth dimension $m-4$ sub-manifold away from a co-dimension 2 singular set with smooth integer multiplicity away from that set. This result is optimal : integration along holomorphic curves in ${\C}{P}^n$ is a calibrated integral cycle for the Fubini Study K\"ahler form and can have isolated singularities which are of co-dimension 2 within the curve.

\medskip

In his paper Tian made the following conjecture.

\medskip

\noindent{\bf Tian regularity conjecture} {\it Let $(E_0,\nabla^\infty)$ be the weak limit of smooth $\Omega$-anti-self dual instantons on a bundle $E$. Then the limiting bundle $E_0$ and the limiting connections $\nabla^\infty$ extend to smooth bundle, resp. smooth connection, away from a closed co-dimension 6
set $L$ in $M^m$}.\hfill $\Box$.

\medskip

There is one case which has been completely settled and where the conjecture has been proved. This is the case of {\it Hermitian Yang-Mills} fields.
In that case the currents defined by $$\varphi\longrightarrow\int_{M^{2m}}Tr(F_{\nabla^k}\wedge F_{\nabla^k})\wedge\varphi$$
is a $(m-2)-(m-2)$ positive cycle (i.e. calibrated by $\om^{m-2}/(m-2)!$), This condition is of course preserved at the limit and then $$\varphi\longrightarrow\int_{M^m}Tr(F_{\nabla^{\infty}}\wedge F_{\nabla^{\infty}})\wedge\varphi +8\pi^2\,C(\varphi)$$
is also a  $(m-2)-(m-2)$ positive cycle (i.e. calibrated by $\om^{m-2}/(m-2)!$). The points of non zero density correspond to the support of $C$ and using an important result by Y.T. Siu \cite{Siu} we obtain that  $C$ is the integration along an holomorphic sub-variety of complex co-dimension $2$.  Using now
a result by S.Bando and Y.T.Siu \cite{BS} we obtain that $E_0$ extends to an analytic reflexive sheaf over the whole $M^{2m}$ it is then locally free (i.e. it realizes a smooth bundle) away from a closed complex co-dimension 3 subset of $M^{2m}$ which is included in $K$ (see for instance \cite{K}). They also prove a point removability asserting that $\nabla^\infty$ defines a smooth connection on the part where the sheaf is free which proves the {\it Tian regularity conjecture}
in the special case of {\it Hermitian Yang-Mills} fields. 

\subsection{The space of weak connections.}

As we saw in the previous section, in dimension larger than 4, ''bundles with singularities'' arise naturally as ''carriers'' of limits of smooth Yang-Mills fields with uniformly bounded energy over smooth bundles. If now we remove the assumption to be Yang-Mills and just follow sequences of connections with uniformly bounded Yang-Mills energy over smooth bundles we have seen in the beginning of this section that the limiting carrying bundle can have twists everywhere
on the base !
Similarly, taking a sequence of closed sub-manifolds with uniformly bounded volume the limit ''escape'' from the space of smooth sub-manifold and can be singular. The main achievement of the work of Federer and Fleming has been to introduce a class of objects, the {\it integral cycles} which complete the space of closed sub-manifold
with uniformly bounded volume and which was suitable to solve the {\it Plateau problem} in a general framework. The purpose of the work \cite{PR3} is to define 
a class of {\it weak bundles} and {\it weak connections} satisfying a closure property under uniformly bounded Yang-Mills energy and suitable to solve the {\it Yang-Mills Plateau problem}. 

\medskip

We introduce the following stratified definition.

\begin{Dfi}
\label{df-VII.1}
Let $G$ be a compact Lie group and 
$(M^m,h)$  a compact riemanian manifold. For $m\le 4$ the space of {\it weak connections} ${\mathcal A}_G(M^m)$ is defined to coincide with the space of  {\it Sobolev connections} defined by
\[
{\mathfrak A}_G(M^m):=\lf\{ 
\begin{array}{l}
A\in L^2(\wedge^1M^m,{\mathcal G})\ ; \ \int_{M^m}|dA+A\wedge A|_h^2\ dvol_h<+\infty\\[5mm]
\mbox{ locally }\exists\ g\in W^{1,2}\quad\mbox{ s.t. }\quad A^g\in W^{1,2}
\end{array} 
\rg\}
\]
For $m>4$ we define the space of {\it weak connections} ${\mathcal A}_G(M^m)$ to be
\[
{\mathcal A}( M^m):=\lf\{ 
\begin{array}{l}
A\in L^2(\wedge^1M^m,{\mathcal G})\ ; \ \int_{M^m}|dA+A\wedge A|_h^2\ dvol_h<+\infty\\[5mm]
\forall\ p\in M^m\quad\mbox{ for a.e. }r>0\quad\iota_{\p B_r(p)}^\ast A\in {\mathcal A}_G(\p B^{m}_r(p))
\end{array} 
\rg\}
\]
where $B^m_r(p)$ denotes  the geodesic ball of center $p$ and radius $r>0$ and $\iota_{\p B_r(p)}^\ast A$ is the restriction of the 1-form $A$ on the boundary
of $B^m_r(p)$.
\hfill $\Box$
\end{Dfi}
As an illustration of this space, It is not difficult to check that the limiting connections $\nabla^\infty$ from Tian's closure  theorem are in ${\mathcal A}_{G}$. We have the following result which justifies the previous definition.
\begin{Th}
\label{th-VII.4}
The space ${\mathcal A}_{G}(B^5)$ is weakly sequentially closed bellow any Yang-Mills energy level. Precisely, let $A^k\in {\mathcal A}_{G}(B^5)$ such that
\[
\limsup_{k\rightarrow +\infty}\int_{B^5}|dA^k+A^k\wedge A^k|^2\ dx^5<+\infty
\]
then there exists a subsequence $k'$ and  $A^\infty\in {\mathcal A}_{G}(B^5)$ such that
\[
d(A^{k'},A^\infty):=\inf_{g\in W^{1,2}(B^5,G)}\int_{B^5}|A^{k'}-(A^\infty)^g|^2\ dx^5\longrightarrow 0
\]
\hfill $\Box$
\end{Th}
We conjecture in \cite{PR3} that this result extends to dimension $m\ge 6$. A proof of theorem~\ref{th-VII.4} is given in \cite{PR3} and is using the following strong approximation theorem whose proof is rather technical.
\begin{Th}
\label{th-VII.5}
Let $A\in {\mathcal A}_{G}(B^5)$ then there exists $A^k$ which are the connection ${\mathcal G}$ 1-forms on $B^5$ associated to smooth connections of a sequence of smooth bundles over $B^5$ minus finitely many points such that
\[
\lim_{k\rightarrow +\infty}\inf_{g\in W^{1,2}(B^5,G)} \int_{M^m}|A-(A^k)^g|_h^2+|F_{A}-g^{-1}F_{A^k} g|^2\ dx^5= 0
\]
\hfill $\Box$
\end{Th}
These two last results can be seen as a non-abelian version of the following weak sequential closure, respectively strong approximation, results proved in \cite{PR1}.
\begin{Th}
\label{th-VII.6}
Let $p>1$ and let $F_k$ be a sequence of $L^p$ two forms in $B^3$ such that 
\[
\forall\ p\in B^3\quad\mbox{ for a.e. }B_r(p)\subset B^3\quad\int_{\p B_r(p)} F_k\in 2\pi\,{\Z}\quad.
\]
Assume $F_k\rightharpoonup F_\infty$ \underbar{weakly} in $L^p$ then
\[
\forall\ p\in B^3\quad\mbox{ for a.e. }B_r(p)\subset B^3\quad\int_{\p B_r(p)} F_\infty\in 2\pi\,{\Z}\quad.
\]
\hfill $\Box$
\end{Th}
The strong approximation result says the following.
\begin{Th}
\label{th-VII.7}
Let $p>1$ and let $F$ be an $L^p$ two forms in $B^3$ such that 
\[
\forall\ p\in B^3\quad\mbox{ for a.e. }B_r(p)\subset B^3\quad\int_{\p B_r(p)} F\in 2\pi\, {\Z}\quad.
\]
then there exists a sequence $A^k$ of smooth connections of a sequence of smooth $U(1)$ bundles over $B^3$ minus finitely many points such that
\[
\lim_{k\rightarrow +\infty}\int_{B^3}|F_{A^k}-F|^p\ dx^3=0
\]
\hfill$\Box$
\end{Th}

\subsection{The resolution of the Yang-Mills Plateau problem in 5 dimensions.}

In order to solve the Yang-Mills Plateau problem in the 5-dimensional ball $B^5$, we introduce a sub-space of ${\mathcal A}_{SU(n)}(B^5)$ of weak
connections admitting a trace.

\medskip

Let $\eta$ be a weak connection 1-form of ${\mathcal A}_{G}(S^4)$ we denote by ${\mathcal A}_{SU(n)}^\eta(B^5)$  the subspace of 
${\mathcal A}_{SU(n)}(B^5)$ made of weak connection 1-forms $A$ such that
\[
\lim_{r\rightarrow 1^{-}}\inf_{g\in W^{1,2}(S^4,G)}\int_{S^4}| D_r^\ast A-g^{-1}\,dg-g^{-1}\eta\, g|^2\ dvol_{S^4}
\]
where $D_r^\ast A$ is the pull-back on $S^4\simeq \p B^5_1(0)$ of the restriction of $A$ to $\p B^5_r(0)$ by the dilation of ratio $r$. 
It is proved in \cite{PR3} that for any $\eta\in {\mathcal A}_{G}$ the space ${\mathcal A}_{G}^\eta$ is weakly sequentially closed under any Yang-Mills energy level. Hence we have the following theorem.

\begin{Th}
\label{th-VII.8}{\bf [Existence for YM minimizers]}
Let $G$ be a compact Lie group and let $\eta$ be an arbitrary Sobolev  connection one form in ${\mathcal A}_G(S^4)$ then the folowing
minimization problem is achieved
\be
\label{VII.6}
\inf\lf\{\int_{B^5}|dA+A\wedge A|^2\ dx^5\ ;\ A\in {\mathcal A}_G^\eta(B^5)\rg\}\quad.
\ee\hfill $\Box$
\end{Th}
\noindent Solutions to this minimization problem will be called {\it solutions to the Plateau problem} for the boundary connection $\eta$. 

Once the existence of the minimizer is established it is legitimate to ask about the regularity for these solutions to the {\it Yang-Mills Plateau problem}.
Before to give the optimal regularity we give an intermediate result which holds in the general class of {\it stationary Yang-Mills} in the space
of {\it weak connections} ${\mathcal A}_{G}(B^5)$.

\begin{Th}
\label{th-VII-2-al}{\bf[$\epsilon-$regularity for weak stationary YM fields in 5D]}
Let  $G$ be a compact Lie group, there exists $\ep_{G}>0$ such that for any weak connection 1-forms $A$ in ${\mathcal A}_G(B^5)$   satisfying \underbar{weakly} the Yang-Mills equation
\[
d^\ast F_A+[A,\res F_{A}]=0\quad\mbox{ in }{\mathcal D}'(B^5_1(0)).
\]
assume also that it satisfies the stationarity condition
\[
\forall\ X\in C^\infty_0(B^5,{\R}^5)\quad\quad\lf.\frac{d}{dt}\int_{B^5}|\Psi_t^\ast F_A|^2\rg|_{t=0}=0
\]
where $\Psi_t$ is the flow of $X$, and the small energy condition
\[
\int_{B^5_1(0)}|F_A|^2\ dx^5<\ep_{G}
\]
then there exists a gauge $g\in W^{1,2}(B^5_{1/2}(0),G)$ in which the following estimates hold : for any $l\in {\N}$ there exists $C_l>0$ such that
\[
\|\nabla^l(A)^g\|^2_{L^\infty(B^5_{1/2}(0))}\le C_l\ \int_{B^5_1(0)}|F_A|^2\ dx^5
\]
\hfill $\Box$
\end{Th}
This theorem is a consequence the {\it monotonicity formula} deduced from the {\it stationarity condition} and a weak version of the {\it Coulomb gauge }extraction theorem~\ref{th-VII-2-b} which generalizes the work of T.Tao- G.Tian \cite{TT} for {\it admissible connections} or the work of Y.Meyer and the author
\cite{MR} for {\it approximable connections} to the  general framework of {\it weak connections} in ${\mathcal A}_G(B^5)$.

\medskip

With the {\it $\epsilon-$regularity } result at hand, using a {\it Luckhaus lemma} together with a {\it Federer dimension reduction method} following the proof of the regularity result of R.Schoen and K.Uhlenbeck for minimizing harmonic map \cite{SU1} we obtain the following theorem.
\begin{Th}{\bf [Regularity for minimizers of the Yang-Mills Plateau problem.]}
\label{th-VII-2-alm}
Let $G$ be a compact Lie group, and let $\eta$ be the connection 1-form associated to a smooth connection of some $G-$bundle over $\p B^5$ then the minimizers of the Yang-Mills Plateau problem (\ref{VII.6}) are smooth connections of a smooth bundle defined on $B^5$ minus finitely many points.\hfill $\Box$
\end{Th}
It is proved in \cite{PR3} that the theorem is optimal in the sense that there exists a smooth $su(2)-$valued 1-form $\eta$ - defining then a connection of the trivial bundle over $S^4$ - such that any solution to the 
Plateau problem (\ref{VII.6}) has isolated singularities.

\end{document}